\documentclass[a4paper]{article}

\usepackage[all]{xy}
\usepackage[latin1]{inputenc}        
\usepackage[dvips]{graphics,graphicx}
\usepackage{amsfonts,amssymb,amsmath,color,mathrsfs, amstext, esint}
\usepackage{amsbsy, amsopn, amscd, amsxtra, amsthm, verbatim, yhmath}
\usepackage{algorithmicx}
\usepackage{mathtools}
\usepackage{algorithm}
\usepackage{upref}
\usepackage{enumitem}
\usepackage{placeins}
\usepackage{appendix}

\usepackage[colorlinks,
linkcolor=red,
anchorcolor=red,
citecolor=red
]{hyperref}
\newtheorem{theorem}{Theorem}[section]
\newtheorem{lemma}{Lemma}[section]
\newtheorem{proposition}[lemma]{Proposition}
\newtheorem*{proposition*}{Proposition}
\newtheorem{corollary}[lemma]{Corollary}
\newtheorem*{corollary*}{Corollary}
\newtheorem{remark}[lemma]{Remark}
\newtheorem{definition}[lemma]{Definition}
\usepackage{geometry}
\usepackage{float}
\allowdisplaybreaks
\newcommand{\vertiii}[1]{{\left\vert\kern-0.25ex\left\vert\kern-0.25ex\left\vert #1 
		\right\vert\kern-0.25ex\right\vert\kern-0.25ex\right\vert}}
\usepackage{yhmath}
\definecolor{DarkGreen}{rgb}{0,0.6,0}
\newcommand{\lw}[1]{\textcolor{DarkGreen}{[LW: #1]}}

\newcommand{\tdet}{t_3}

\newcommand{\ud}{\,\mathrm{d}}
\newcommand{\rr}{r_*}
\newcommand{\rss}{r_{**}}
\newcommand{\N}{\mathbb{N}}
\newcommand{\R}{\mathbb{R}}
\newcommand{\opS}{\mathcal{S}}
\newcommand{\Z}{\mathbb{Z}}

\newlength{\leftstackrelawd}
\newlength{\leftstackrelbwd}
\def\leftstackrel#1#2{\settowidth{\leftstackrelawd}%
{${{}^{#1}}$}\settowidth{\leftstackrelbwd}{$#2$}%
\addtolength{\leftstackrelawd}{-\leftstackrelbwd}%
\leavevmode\ifthenelse{\lengthtest{\leftstackrelawd>0pt}}%
{\kern-.5\leftstackrelawd}{}\mathrel{\mathop{#2}\limits^{#1}}}

\title{Optimal artificial boundary conditions based on second-order correctors for three dimensional random elliptic media}
\author{Jianfeng Lu \thanks{Mathematics Department, Duke University, Box 90320, Durham, NC, 27708, USA (jianfeng@math.duke.edu)} \and Felix Otto\thanks{Max Planck Institute for Mathematics in the Sciences, Inselstr. 22, 04103, Leipzig, Germany (otto@mis.mpg.de)} \and Lihan Wang\thanks{Department of Mathematical Sciences, Carnegie Mellon University, 311 Hamerschlag Drive, Pittsburgh, PA, 15213, USA (lihanw@andrew.cmu.edu)}}
\date{\today}

\begin{document}
\maketitle

\begin{abstract}

We are interested in numerical algorithms for computing the electrical field generated
by a charge distribution localized on scale $\ell$ in an infinite heterogeneous medium, 
in a situation where the medium is only known in a box of diameter $L\gg\ell$
around the support of the charge. We propose a boundary condition that with overwhelming
probability is (near) optimal with respect to scaling in terms of $\ell$ and $L$,
in the setting where the medium is a sample from a stationary ensemble with a finite range of dependence
(set to be unity and with the assumption that $\ell \gg 1$). 
The boundary condition is motivated by quantitative stochastic homogenization 
that allows for a multipole expansion \cite{bella2017effective}.

This work extends \cite{lu2018optimal}, the algorithm in which is optimal in two dimension, and thus we need to take
quadrupoles, next to dipoles, into account. This in turn relies on stochastic estimates
of second-order, next to first-order, correctors. These estimates are provided for finite range ensembles
under consideration, based on an extension of the semi-group approach of \cite{GO2015most}.
\end{abstract}
\tableofcontents

\section{Introduction and Main Results}

Consider a conducting medium as described by a symmetric $\lambda$-uniform coefficient field
$a=a(x)$ in $d$-dimensional space, that is, for any $x, \xi \in \R^d$
\begin{align}\label{eqn:intrunifell}
\lambda|\xi|^2\le \xi\cdot a(x)\xi\le|\xi|^2.
\end{align}
Consider a localized charge distribution that is overall neutral, as described
by a compactly supported dipole density $g$. Let us give a sense to its characteristic scale
$\ell$ by assuming that it is of the form
\begin{align}\label{eqn:conditionrhs}
g(x)=\hat g(\frac{x}{\ell})
\end{align}
for some sufficiently smooth $\hat g$ supported in the unit ball. We are interested
in the field $\nabla u$ the charge generates, which is the decaying solution of
the elliptic divergence-form equation 
\begin{align}\label{eqn:intrbaseq}
\nabla\cdot(a\nabla u+g)=0\quad\mbox{in}\;\mathbb{R}^d.
\end{align}
%
%

In this paper, we address the following question: Suppose we only know the medium
$a$ in some box $Q_{2L}:=(-2L,2L)^d$, to what precision may we infer the value of
$\nabla u$? What is a practical algorithm to retrieve it?
Heuristically, one expects $\nabla u(x)$ to decay as a dipole, i.e., like $(\frac{\ell}{|x|})^d$, thus we expect that changing the coefficient field $a$ outside the box $Q_{2L}$ will affect $\nabla u$ to order $(\frac{\ell}{L})^d$,
and imposing homogeneous Dirichlet conditions on $\partial Q_L$ would do no worse
-- and this would be the end of the story and the paper. 

In this paper, however, we consider a more specific situation, namely when $a$
is sampled from a stationary ensemble $\langle\cdot\rangle$, which puts us
into the context of stochastic homogenization. More precisely, 
we shall assume
that $\langle\cdot\rangle$ is of finite range, which we set to be unity without loss of generality.
This means that for two sets $D$ and $D'\subset\mathbb{R}^d$ with distance larger than $1$, the restrictions $a_{|D}$ and $a_{|D'}$ are independent. What information
may we retrieve in this case? Consider again changing the coefficient field outside $Q_{2L}$ (now with $L\gg 1$) but keeping the statistical ensemble, heuristically, we expect that due to stochastic cancellations, the impact on $\nabla u$ reduces to $\frac{1}{\sqrt{L^d}} (\frac{\ell}{L})^d$.
Indeed, the additional attenuation factor $\frac{1}{\sqrt{L^d}}$ comes from the Central Limit Theorem (CLT) scaling involving the square root
of the relevant volume, non-dimensionalized by the correlation length. In fact,
this precision on inferring the value of $\nabla u$ cannot be improved, as the following
lower bound on the variance of $\nabla u$ from previous work \cite{lu2018optimal}, conditioned on the 
restriction $a_{|Q_{2L}}$, shows, which we expect to hold for generic ensembles.

\begin{theorem}\label{thm:luottooptimal} \cite[Theorem 2]{lu2018optimal}
There exists a stationary, unit-range ensemble
$\langle\cdot\rangle$ supported on $a$'s satisfying \eqref{eqn:intrunifell} with the following property: Consider the solution $u$ of \eqref{eqn:intrbaseq}, where $g$ is of the form \eqref{eqn:conditionrhs} for some $\ell$ and $\hat{g}$,
then there exists a radius $R$ such that for any $\omega=\frac{1}{R^d}\hat\omega(\frac{x}{R})$
for some sufficiently smooth $\hat\omega$ supported in $B_1$ with $\int\hat\omega=1$,
\begin{align*}
\big\langle\big|\int\omega\nabla u-\langle\int\omega\nabla u \mid a_{|Q_{2L}}\rangle\big|^2\big\rangle^\frac{1}{2}
\ge \frac{1}{C}(\frac{\ell}{L})^d(\frac{1}{L})^\frac{d}{2}\quad\mbox{ provided } \; \
\frac{L}{C}\ge \ell \ge C.
\end{align*}
Here the radius $R$ and the constant $C$ depend only on the ensemble, $\hat g$, and $\hat\omega$.
\end{theorem}

In \cite{lu2018optimal}, a practical algorithm, that saturates this scaling for $d=2$, was proposed and analyzed. The error of such an algorithm is $O((\frac{\ell}{L})^d(\frac{1}{L})^{1-})$ in any dimension. In this paper, we tackle the more physically relevant case of $d=3$, which requires a
substantial modification of the algorithm and its analysis.

\medskip

We will propose a deterministic algorithm, Algorithm \ref{alg:truealg}, that involves the realization $a$ only in terms of its restriction $a_{|Q_{2L}}$. The algorithm saturates the lower bound of Theorem \ref{thm:luottooptimal} in terms of scaling. More specifically, by solving a couple of auxiliary boundary value problems with homogeneous Dirichlet boundary conditions on $Q_{2L}$, $Q_{\frac{7}{4}L}$, and $Q_{\frac{3}{2}L}$, this algorithm constructs Dirichlet boundary data $u_L$ on $\partial Q_L$, which in turn defines $u^{(L)}$, the output of the algorithm, by solving \begin{equation}\label{eqn:diropt}
-\nabla \cdot a \nabla u^{(L)}=\nabla \cdot g\text{ in }Q_L,\hspace{0.3in} u^{(L)}=u_L\text{ on }\partial Q_L.
\end{equation}For $L,\ell,R \gg 1$ and with overwhelming probability, 
this algorithm saturates the lower bound of Theorem \ref{thm:luottooptimal}, with the little caveat that the CLT exponent, $\frac{3}{2}$ for $d=3$, has to be replaced by 
$\beta<\frac{3}{2}$. Part of the probabilistic nature of the statement is
contained in the random radius $r_{**}$, which can be interpreted as the scale
from which onwards stochastic homogenization is effective.

\begin{theorem}\label{thm:mainthm} Let $d= 3$ and $\langle\cdot\rangle$ be a stationary, unit-range ensemble 
supported on
$a$'s satisfying \eqref{eqn:intrunifell}. Let $g$ be of the form \eqref{eqn:conditionrhs} for some $\ell$ and $\hat{g}$, let $u$ denote the solution of \eqref{eqn:intrbaseq}, and $u^{(L)}$ be the output of Algorithm \ref{alg:truealg} for some $L\ge \ell$. Then for any $\beta<\frac{3}{2}$, there exists a random radius $r_{**}$ such that conditioning on $\ell\ge \rss$, with probability at least $1 - \exp(-L^{1/C})$, we have for any $R\in [\rss,L]$,\footnote{Here $\fint$ denotes spatial average} 
\begin{equation*}
\Big(\fint_{B_R}|\nabla(u^{(L)}-u)|^2\Big)^\frac{1}{2}
\le C(\frac{\ell}{L})^d(\frac{r_{**}}{L})^\beta \quad\mbox{ provided } \;
\frac{L}{C}\ge  \ell \ge C.
\end{equation*}
Moreover, the radius $r_{**}$ satisfies
\begin{align}\label{eqn:thm1rss}
\bigl\langle\exp(r_{**}^\beta)\bigr\rangle\le C.
\end{align}
Here $C$ denotes a constant that depends only on $\lambda$, $\hat g$, $\beta$ and $\hat\omega$ appearing in Theorem \ref{thm:luottooptimal}.
\end{theorem}

We believe that with some additional work, it is possible to derive an ``a posteriori'' style result similar to \cite{lu2018optimal}, that is, we could define some computable $\rss^{(L)}$ that plays the same role as $\rss$ in Theorem~\ref{thm:mainthm}. We also comment here that Theorem \ref{thm:mainthm} holds with any $\beta <2$ for $d\ge 4$, and the algorithm is thus also near-optimal when $d=4$. Obtaining the optimal algorithm for $d\ge 5$ requires computing correctors of order three or higher and we do not discuss it here.

\smallskip

Before we further discuss ideas of Algorithm \ref{alg:truealg} and the proof of the theorem, let us first compare our result with the previous work \cite{lu2018optimal}. A main difference lies in the introduction of the functions $\psi_T^{(L)}$ and the corresponding coefficients $c_T^{(L)}$, which are approximations of second-order correctors and quadrupoles that will be introduced below. These are available for $d>2$ and indeed necessary for the algorithm to (almost) reach the CLT-scaling $\beta<\frac{3}{2}$. This aspect of difference will be discussed in more details in Section \ref{subsec:multipole}. A more technical difference is that Algorithm~\ref{alg:truealg} uses a different approximation $\phi_T^{(L)}$ of the first-order corrector $\phi$, compare \eqref{eqn:phiTL} to \cite[(18)]{lu2018optimal}. The difference lies in the massive term $\frac{1}{T}\phi_T^{(L)}$. The reason of such change is discussed in Section \ref{subsec:mascor}.

\begin{algorithm}
\caption{Optimal algorithm for the approximate solution $u^{(L)}$ in $Q_L$ \label{alg:truealg}}
\begin{algorithmic}[1]
\State For $\beta\in (1,\frac{3}{2})$ set $\varepsilon = \frac{3}{2}-\beta$ and $T = L^{2(1-\varepsilon)}$. For $i=1,\cdots,d$, solve for the approximate first-order corrector $\phi_{i,T}^{(L)}$:
\begin{equation}\label{eqn:phiTL}
 \dfrac{1}{T}\phi_{i,T}^{(L)}-\nabla  \cdot a \nabla \phi_{i,T}^{(L)} =\nabla\cdot ae_i \,  \mbox{ in }Q_{2L}, \hspace{0.3in} \phi_{i,T}^{(L)}=0 \, \mbox{ on }\partial Q_{2L}. 
\end{equation}
\State Calculate the approximate homogenized coefficients via \begin{equation}\label{eqn:algahL}
    a_h^{(L)}e_i=\int \omega q_{i,T}^{(L)},
\end{equation} where \begin{equation}\label{eqn:defqiTL} q_{i,T}^{(L)}:=a(e_i+\nabla\phi_{i,T}^{(L)})\end{equation} and $\omega(x)=\frac{1}{L^d}\hat{\omega}(\frac{x}{L})$ with $\hat{\omega}$ as in Theorem \ref{thm:luottooptimal}. 
\State Find $\tilde{u}_h^{(L)}$ on $\partial Q_L$: 
\begin{equation}\label{eqn:alguhtildeL}
    \tilde{u}_h^{(L)} =\int G_h^{(L)}* (\nabla\cdot g),
\end{equation}where $G_h^{(L)}(x) := \frac{1}{4\pi\left|(a_h^{(L)})^{-1/2}x\right|}$ is the Green function for the constant-coefficient operator $-\nabla\cdot a_h^{(L)} \nabla$.
\State Solve for approximate first-order flux correctors $\sigma_{i,T}^{(L)}=\{\sigma_{ijk,T}^{(L)}\}_{j,k}$:  \begin{equation}\label{eqn:algsigma}
    \dfrac{1}{T}\sigma_{ijk,T}^{(L)}-\Delta \sigma_{ijk,T}^{(L)} =\partial_j q_{ik,T}^{(L)}-\partial_k q_{ij,T}^{(L)} \, \mbox{ in }Q_{\frac{7}{4}L}, \hspace{0.3in} \sigma_{ijk,T}^{(L)}=0 \, \mbox{ on }\partial Q_{\frac{7}{4}L}. 
\end{equation}
\State Solve for approximate second-order correctors $\psi_{ij,T}^{(L)}$:  \begin{equation}\label{eqn:2ndcorapprox} 
     \dfrac{1}{T}\psi_{ij,T}^{(L)} - \nabla\cdot a \nabla \psi_{ij,T}^{(L)} = \nabla\cdot (\phi_{i,T}^{(L)}a-\sigma_{i,T}^{(L)})e_j \, \mbox{ in }Q_{\frac{3}{2}L}, \hspace{0.3in} \psi_{ij,T}^{(L)}=0 \,\mbox{ on }\partial Q_{\frac{3}{2}L}. 
\end{equation}
\State For the indices \begin{equation}\label{eqn:calJ}(i,j)\in \mathcal{J}=\{(1,2),(1,3),(2,3),(2,2),(3,3)\},\end{equation} calculate  \begin{equation}\label{eqn:cijlt}
    c_{ij,T}^{(L)}=-\int g\cdot  \nabla \Bigl(\sum_{k=1}^3\phi_{k,T}^{(L)}\partial_k v_{h,ij}^{(L)}+(2-\delta_{ij})(\psi_{ij,T}^{(L)}-\dfrac{a_{hij}^{(L)}}{a_{h11}^{(L)}}\psi_{11,T}^{(L)})\Bigr) ,
\end{equation} where $v_{h,ij}^{(L)}$ denote the $a_h^{(L)}$-harmonic polynomials \begin{equation}\label{eqn:harmpolL}
    v_{h,ij}^{(L)}=(1-\dfrac{1}{2}\delta_{ij})(x_ix_j-\dfrac{a_{hij}^{(L)}}{a_{h11}^{(L)}}x_1^2).
\end{equation}
\State Obtain $u_h^{(L)}$ as \begin{equation}\label{eqn:algapproxbdry}
    u_h^{(L)}=\tilde{u}_h^{(L)}+ \sum_{i=1}^3(\int  g \cdot\nabla \phi_{i,T}^{(L)})\partial_i G_h^{(L)} +\sum_{(i,j)\in\mathcal{J}}c_{ij,T}^{(L)}\partial_{ij} G_h^{(L)}.
\end{equation} 
\State Solve for $u^{(L)}$ (here and for the rest of the paper we adopt Einstein's summation convention for repeated indices): \begin{equation}\label{eqn:finalapprox}
   -\nabla \cdot a \nabla u^{(L)}=\nabla \cdot g\text{ in }Q_L,\hspace{0.3in} u^{(L)}=(1+\phi_{i,T}^{(L)}\partial_i+\psi_{ij,T}^{(L)}\partial_{ij}) u_h^{(L)}\text{ on }\partial Q_L.
\end{equation} 
\end{algorithmic}
\end{algorithm}

\subsection*{Related Works}\label{subsec:reference} 

Quantitative stochastic homogenization, which dates back to Yurinskii \cite{yurinskii1986averaging}, has been intensively studied during the past decade. Naddaf and Spencer \cite{naddaf1997homogenization,naddaf1998estimates} introduced the notion of spectral-gap inequality and captured the CLT-type optimal scaling in stochastic homogenization under the condition of small ellipticity contrast, which is removed in \cite{gloria2011optimal, gloria2012optimal,gloria2015quantification} for discrete coefficients. The spectral-gap inequality is then refined to logarithmic Sobolev inequality in \cite{marahrens2015annealed} and further generalized to ensembles with potentially thick correlation tails in \cite{gloria2019quantitative,gloria2014regularity}. Another approach by Armstrong and Smart \cite{armstrong2016quantitative}, Armstrong, Kuusi and Mourrat \cite{armstrong2016mesoscopic,armstrong2017additive} uses a variational approach and obtain similar results. See \cite{armstrong2017quantitative} for a nice monograph. The Gaussianity of the energy of the solution was first identified by Nolen \cite{nolen2014normal}, and has been extended to homogenized coefficients by Gloria and Nolen \cite{gloria2016quantitative}, both for the representative volume element method. The covariance structure of correctors have been studied in \cite{mourrat2016correlation}, based on the annealed Green's function estimates in \cite{marahrens2015annealed}. This is extended to study the fluctuation of homogenization error in \cite{gu2016scaling}. Scaling limit of the correctors have been studied in \cite{mourrat2017scaling}. Second-order correctors have been studied in \cite{fischer2016higher, bella2017stochastic}. Correctors of higher order have been constructed in \cite{gu2017high} up to order of $\lceil\frac{d}{2}\rceil$ with suboptimal error estimate in two-scale expansion, which is improved to be optimal in \cite{duerinckx2020higher}. Higher-order correctors up to order $d$ were recently constructed in \cite{duerinckx2021non} using a distributional formulation. In a perturbative regime, it is possible to define homogenized coefficients up to order $2d$ \cite{duerinckx2019remark} by using a breakthrough result of Bourgain \cite{bourgain2018homogenization}, in its refined form established in \cite{kim2019averaged}.

\smallskip

Developing efficient numerical algorithms has been a major motivation behind the study of quantitative stochastic homogenization, see e.g., the review \cite{anantharaman2012introduction}. Let us just mention a few directions here: Quantitative error estimates for computation of effective coefficients in stochastic homogenization have been studied in \cite{bourgeat2004approximations,gloria2012numerical,gloria2012spectral,egloffe2015random},  where strategies using different boundary conditions or massive terms have been studied. Representative volume method, a popular approach used by engineers, are systematically analyzed in \cite{fischer2019choice, khoromskaia2019numerical}. Iterative multigrid methods have been studied in \cite{mourrat2019efficient, hannukainen2019computing,armstrong2018iterative}. There have also been abundant research in numerical homogenization where the coefficients do not necessarily arise from random setting. An approach using an embedded corrector problem for approximating homogenized coefficients has been considered in \cite{cances2015embedded}. Variance reduction methods have been developed in \cite{le2016special,blanc2016some}. Multiscale finite element methods have been developed in \cite{babuvska1983generalized, babuvska1994special} and extended in \cite{hou1997multiscale}, see \cite{efendiev2009multiscale} for a review. Heterogeneous multiscale methods have been developed and studied in \cite{weinan2003heterogeneous, ming2005analysis} and many other works, see \cite{weinan2007heterogeneous,abdulle2012heterogeneous} for reviews. In \cite{oden2000estimation} they proposed a method that aims at recovering local oscillations by solving a local problem using the approximate solution of the homogenized equation as its boundary condition. Localized orthogonal decomposition approaches have been studied in \cite{maalqvist2014localization, fischer2019priori}. Relationship between numerical homogenization and Bayesian inference have been investigated in \cite{owhadi2015bayesian,owhadi2017multigrid}. 

\section{Ideas Behind the Algorithm}
\subsection{Correctors, Homogenized Coefficients and Fluxes, Effective Multipoles}\label{subsec:multipole}
Let us recall the standard first-order correctors $\phi_i$, which play a central role in homogenization theory \cite{kozlov1979averaging, papanicolaou1979boundary}. For each direction $i=1,\cdots,d$, the first-order corrector $\phi_i$ is defined such that $x_i+\phi_i$ is $a$-harmonic\footnote{See below for a uniqueness argument}: 
\begin{equation}\label{intrphi}
-\nabla\cdot a(e_i+\nabla\phi_i)=0.
\end{equation}
Provided the ensemble is stationary and ergodic, the operator $-\nabla \cdot a \nabla \ $ homogenizes on large scale to $-\nabla \cdot a_h \nabla$, with the constant and deterministic homogenized coefficient $a_h$ given by \begin{equation}\label{intrhomcoeff}
a_he_i:=\langle q_i \rangle=\lim\limits_{L\uparrow\infty}\fint_{B_L}q_i\hspace{0.1in}\text{ where }\hspace{0.1in}q_i:=a(e_i+\nabla\phi_i).
\end{equation}
Given $a_h$, we define $\tilde{u}_h$ to be the solution of the homogenized equation  
\begin{equation}\label{intruhtilde}
-\nabla\cdot  a_h\nabla \tilde{u}_h=\nabla\cdot g.
\end{equation} 

\smallskip

The flux corrector $\sigma$, introduced in \cite{gloria2014regularity} in the setting of stochastic homogenization (see earlier ideas in periodic homogenization in \cite[Proposition 7.2]{jikov2012homogenization}), plays a convenient role in quantitative analysis. Since $q_i-a_he_i$, which can be viewed as a $(d-1)$-form, has zero expectation and is divergence free in view of \eqref{intrphi} and \eqref{intrhomcoeff}, there exists a $(d-2)$-form, which can be expressed as a skew-symmetric tensor field $\sigma_i$, such that 
\begin{equation}\label{eqn:divsigma}
    q_i-a_he_i=\nabla\cdot\sigma_i,
\end{equation} 
where we use the convention $(\nabla\cdot\sigma_i)_j= \partial_k \sigma_{ijk}$.
Clearly $\sigma_i$ is only determined up to a $(d-3)$-form, so that a gauge has to be chosen to make its construction unique. It is convenient to choose \begin{equation}\label{eqn:intrsig}
-\Delta\sigma_{ijk}= \partial_j q_{ik}-\partial_k q_{ij}= \nabla \cdot (q_{ik}e_j-q_{ij}e_k).
\end{equation}
With the help of the $\sigma_i$'s one can express the homogenization error in divergence form: for any $a_h$-harmonic function $u_h$, the two-scale expansion $(1+\phi_i\partial_i)u_h$ is close to being $a$-harmonic in the sense of \begin{equation}\label{eqn:1sthomerr}
    -\nabla\cdot a \nabla (1+\phi_i\partial_i) u_h = \nabla \cdot (\phi_i a-\sigma_i)\nabla\partial_i u_h.
\end{equation}
The functions $(\phi,\sigma)$ are uniquely determined up to a random constant by requiring $(\nabla \phi, \nabla \sigma)$ to be stationary fields (i.e. $\nabla \phi(a(\cdot+z),x)=\nabla \phi(a,x+z)$ for any $x,z\in \R^d$) with finite second moments and zero mean \cite[Lemma 1]{gloria2014regularity}. In dimension $3$ under the unit-range assumption, $(\phi,\sigma)$ themselves are stationary, and have finite stretched exponential moments \cite[Theorem 4.1, Proposition 6.2]{armstrong2017quantitative}, \cite[Corollary 2]{GO2015most}. Hence $(\phi,\sigma)$ are uniquely determined by requiring $\langle \phi \rangle = \langle \sigma \rangle =0$ (see discussions in \cite{gloria2014regularity}). 
\smallskip

In dimension $3$, it is well-known that enriching the two-scale expansion by second-order correctors $\psi_{ij}$ leads to a better approximation \cite{bella2017stochastic}. More precisely, given an $a_h$-harmonic function $u_h$, we may improve our two-scale expansion to $(1+\phi_i\partial_i+\psi_{ij}\partial_{ij})u_h$, which is a second-order approximation of an $a$-harmonic function. The characterizing property of second-order correctors is the following equation \footnote{It might be possible to construct $\psi$ via $-\nabla\cdot a (\nabla \psi_{ij}+\phi_ie_j) = e_j\cdot (q_i-a_he_i)$ in order to completely avoid computing $\sigma_i$. We will not discuss the details here.} \begin{equation}\label{eqn:2ndcordef}
-\nabla \cdot a \nabla\psi_{ij}=\nabla\cdot (\phi_i a-\sigma_i)e_j. 
\end{equation} Equation \eqref{eqn:2ndcordef} encodes the following property of the $\psi_{ij}$'s: for an $a_h$-harmonic quadratic polynomial $u_h$, $(1+\phi_i\partial_i+\psi_{ij}\partial_{ij})u_h$ is exactly $a$-harmonic. In practice we only need $\{\psi_{ij}\}_{i,j}$ in form of linear combinations $E_{ij}\psi_{ij}$ with coefficients $\{E_{ij}\}_{i,j}$ that are symmetric and satisfy the trace condition $a_{hij}E_{ij}=0$, which is a five dimensional space for $d=3$. Hence, it suffices to restrict to $\{\psi_{ij}\}_{(i,j)\in \mathcal{J}}$, where $\mathcal{J}$ is the index set defined in Algorithm \ref{alg:truealg}. 
\smallskip

Analogously to first-order corrector case, we need second-order flux correctors, which were first introduced in \cite{bella2017stochastic} for quantitative analysis. Since by \eqref{eqn:2ndcordef}, \begin{equation} \label{intrpsidivfr}
p_{ij}:=a\nabla \psi_{ij}+(\phi_i a-\sigma_i)e_j
\end{equation}
is divergence free, there exist $\Psi_{ij}=\{\Psi_{ijkn} \}_{k,n=1,\cdots, d}$, skew-symmetric with respect to $k$ and $n$, such that \begin{equation}\label{intrpijdiv}
p_{ij}=\nabla\cdot \Psi_{ij}.
\end{equation}
Since $a$ is symmetric, the second-order homogenized coefficient vanishes (see discussions in \cite[calculations on page 5]{bella2017stochastic} or \cite[Corollary 1]{bella2017effective}). Similar to the argument for $\sigma_i$, $\Psi_{ij}$ can be viewed as a $(d-2)$-form and is determined only up to a $(d-3)$-form, and the following choice of gauge is convenient:
\begin{equation}\label{intrPsi}
-\Delta \Psi_{ijkn}=\partial_k p_{ijn}-\partial_n p_{ijk}.
\end{equation}
We learn from \eqref{eqn:2ndcordef} that $\nabla \psi_{ij}$, and thus $p_{ij}$ by \eqref{intrpsidivfr}, and $\nabla \Psi_{ijkn}$ by \eqref{intrPsi}, can be constructed to be stationary if and only if $\phi_i$ and $\sigma_i$ are, which is the case only for $d>2$. Although $\Psi$ is not needed to formulate our algorithm, it will be used in our analysis to upgrade the homogenization error \eqref{eqn:1sthomerr} from first to second-order: for any $a_h$-harmonic function $u_h$,\begin{equation*}
    -\nabla\cdot a \nabla (1+\phi_i\partial_i +\psi_{ij}\partial_{ij} )u_h= -\nabla \cdot (\psi_{ij}a-\Psi_{ij})\nabla \partial_{ij} u_h.
\end{equation*}

 The second-order correctors $\psi$ and $\Psi$ are not expected to be stationary, but rather to grow at a rate a bit worse that $2-\frac{d}{2}=\frac{1}{2}$ away from the origin, which we capture through an exponent $\beta<\frac{d}{2}=\frac{3}{2}$ that measures the amount the growth rate stays below quadratic. The first-order correctors $\phi$ and $\sigma$ are stationary for $d>2$, but they are of course not bounded on $\R^d$; we capture this through an exponent $\alpha \in (\frac{1}{2},1)$ that measures the amount the growth rate below linear. Following \cite{bella2017stochastic, bella2017effective}, we introduce the random radius $\rss\ge 1$ starting from which we have the desired growth rate:
\begin{equation}\label{intrpsiorig}
\dfrac{1}{r^2}\Bigl(\fint_{B_r}|(\psi,\Psi)-\fint_{B_r} (\psi,\Psi)|^2\Bigr)^{\frac{1}{2}} \le (\frac{\rss}{r})^\beta\hspace{0.1in}\text{ for all }r\ge \rss,
\end{equation}
 and \begin{equation}\label{intrphigr}
\dfrac{1}{r}(\fint_{B_r}|(\phi,\sigma)|^2)^{\frac{1}{2}} \le (\frac{\rss}{r})^\alpha\hspace{0.1in}\text{ for all }r\ge \rss.
\end{equation} For convenience we take $\alpha=1-\varepsilon$ and $\beta=\frac{3}{2}-\varepsilon$ for the $\varepsilon>0$ fixed in Algorithm \ref{alg:truealg}. This $\rss$ is the one that appears in the statement of Theorem \ref{thm:mainthm}.
\smallskip

As observed in \cite{bella2015quantitative, bella2017effective}, if we solve the boundary value problem \begin{equation*}
    -\nabla \cdot a \nabla u_1=\nabla \cdot g\text{ in }Q_L,\hspace{0.3in} u_1=(1+\phi_i\partial_i+\psi_{ij}\partial_{ij})\tilde{u}_h\text{ on }\partial Q_L,
\end{equation*}
the solution has error $O((\frac{\ell}{L})^d)$, as it fails to capture the correct multipole behavior, which is the far-field behavior generated from the intrinsic moments of the localized r.h.s.\ $\nabla\cdot g$. We now recall the discussions in \cite{bella2017effective} and describe the far-field behavior of $u$ in order to design the correct boundary condition on $\partial Q_L$ to solve \eqref{eqn:diropt}.
\smallskip

Let $X_k$ be the space of $a$-harmonic functions of growth rate $\le k$, more precisely\footnote{Here $B_r^c:= \mathbb{R}^d \backslash B_r$.} \[X_k:=\Big\{u: -\nabla \cdot a \nabla u=0, ~~~ \limsup_{R\to\infty} R^{1-k}\Bigl(\fint_{B_R}|\nabla u|^2 \Bigr)^\frac{1}{2}<\infty\Big\},\] and $Y_k(r)$ be $a$-harmonic functions outside $B_r$ with decay rate $\ge \, k+d-2$, more precisely, \[Y_k(r):=\Big\{u: -\nabla \cdot a \nabla u=0 \text{ in }B_r^c, ~~~ \limsup_{R\to\infty} R^{k+(d-2)+1}\Bigl(\dfrac{1}{R^d}\int_{B_R^c}|\nabla u|^2 \Bigr)^\frac{1}{2}<\infty\Big\}.\] We also use $X_k^h$ and $Y_k^h(r)$ to denote similar spaces with $a_h$ in place of $a$ in the definition of $X_k$ and $Y_k(r)$ respectively. The spaces $X_k$ and $Y_k(r)$ are related through the bilinear form: for $u\in Y_k(r)$ and $v\in X_k$, \begin{equation}\label{eqn:bilform}
    (u,v)_a:=-\int \nabla\eta\cdot (va\nabla u-ua\nabla v),
\end{equation} where the cutoff function $\eta\equiv 1$ in $B_r$ and $\eta\equiv 0$ in $B_{2r}^c$. Note that the value of the integral does not depend on the choice of $\eta$ since $va\nabla u-ua\nabla v$ is divergence-free outside of $B_r$. We use $(u_h,v_h)_{a_h}$ to denote the bilinear form similarly defined as \eqref{eqn:bilform} with $a_h$ in place of $a$.  In the case of a constant coefficient $a_h$, the spaces $X_k^h$ and $Y_k^h(r)$ are well-understood: $X_k^h$ is the finite-dimensional linear space spanned by $a_h$-harmonic polynomials of degree at most $k$, while $Y_k^h(r)$ is the infinite-dimensional linear space spanned by $k$-th or higher derivatives of $G_h$, the Green's function of $-\nabla\cdot a_h\nabla$. In the language of electrostatics, the quotient spaces $Y_1^h(r)/Y_2^h(r)$ and $Y_2^h(r)/Y_3^h(r)$ are spanned by first and second derivatives of $G_h$ respectively, and thus represent dipoles and quadrupoles. As shown in \cite[Proposition 5]{bella2017effective}, similar to the pairing $(\cdot,\cdot)_{a_h}$ providing an isomorphism between $Y_1^h(r)/Y_3^h(r)$ and $(X_2^h/X_0^h)^*$, $(\cdot,\cdot)_a$ provides an isomorphism between $Y_1(r)/Y_3(r)$ and $(X_2/X_0)^*$. For Liouville principles which describe the equivalence between spaces $X_k^h$ and $X_k$, see also \cite{armstrong2016mesoscopic, fischer2016higher}.

By \cite[Lemma 4]{bella2017effective}, we know that $u\in Y_1(\ell)$. Therefore, by \cite[Theorem 2]{bella2017effective}, under the assumptions \eqref{intrpsiorig} and \eqref{intrphigr}, for $\ell\ge \rss$, there exists a $u_h \in Y_1^h(\ell)$, uniquely defined up to an element in $Y_3^h(\ell)$ by \begin{equation}\label{eqn:lm1homphsm}
(u,v)_a=(u_h,v_h)_{a_h}, \ \text{ for any } v_h\in X_2^h \ \text{ where }\ v:=(1+\phi_i\partial_i+\psi_{ij}\partial_{ij})v_h\in X_2,
\end{equation} that captures the effective multipole behavior of $u$. We state a modification of this result in the following Proposition \ref{prop:effectivemultipole}. \begin{proposition}\label{prop:effectivemultipole}
	Consider a coefficient field $a$ on $\R^d$ satisfying \eqref{eqn:intrunifell}. Suppose that there exists a tensor $a_h$ and, for $i=1,2,3$, a scalar field $\phi_i$ and a skew-symmetric tensor field $\sigma_i$ such that \eqref{intrphi}, \eqref{intrhomcoeff} and \eqref{eqn:divsigma} hold, and for $i,j=1,2,3$, a scalar field $\psi_{ij}$ and a skew-symmetric tensor field $\Psi_{ij}$ such that \eqref{eqn:2ndcordef}, \eqref{intrpsidivfr} and \eqref{intrpijdiv} hold. For fixed $\varepsilon>0$ and $\alpha=1-\varepsilon, \beta=\frac{3}{2}-\varepsilon$, suppose that there exists a radius $\rss$ such that \eqref{intrpsiorig} and \eqref{intrphigr} hold. Moreover, let us assume the following convergences in distribution as $R\to \infty$:
	\begin{align*}
	    \phi_i(R\cdot) & \rightharpoonup 0 \\ \sigma_{i}(R\cdot) & \rightharpoonup 0 \\ \big(e_j\cdot \sigma_i \nabla \phi_k -e_j\cdot \sigma_k \nabla \phi_i\big)(R\cdot) & \rightharpoonup 0. \stepcounter{equation} \tag{\theequation} \label{eqn:qualtlgscconv}
	\end{align*}
	
    Let $g$ be of the form \eqref{eqn:conditionrhs} for some $\ell\ge \rss$ and $\hat{g}$, let $u$ be the solution of \eqref{eqn:intrbaseq}, and let $u_h \in Y_1^h(\ell)$ satisfy \eqref{eqn:lm1homphsm}. Then we have for any $R\ge \rss$, \begin{equation*}
	\Bigl(\frac{1}{R^d}\int_{B_R^c} \Bigl|\nabla \Bigl(u-(1+\phi_i\partial_i+\psi_{ij}\partial_{ij})u_h\Bigr)\Bigr|^2\Bigr)^{\frac{1}{2}} \le C(\frac{\ell}{R})^d(\frac{\rss}{R})^\beta,
	\end{equation*} 
	where $C$ is a constant depending only on $d$, $\lambda$, $\varepsilon$ and $\hat{g}$. 
\end{proposition}\noindent Solving \eqref{eqn:diropt} with boundary condition $u_L=(1+\phi_i\partial_i+\psi_{ij}\partial_{ij})u_h$ will indeed provide the desired level of approximation, which is the statement of Corollary \ref{cor:multipole} below. 
\begin{corollary}\label{cor:multipole}
	Under the same assumptions as in Proposition \ref{prop:effectivemultipole}, for $L\ge \ell$ let $\hat{u}$ be the solution of the equation \begin{equation}\label{eqn:coruhat} -\nabla \cdot a \nabla \hat{u}=\nabla\cdot g \text{ in } Q_L, \hspace{0.3in} \hat{u}=(1+\phi_i\partial_i+\psi_{ij}\partial_{ij})u_h \text{ on } \partial Q_L,\end{equation} where $u_h\in Y_1^h(\ell)$ satisfies \eqref{eqn:lm1homphsm}. Then we have for any $R\in[\rss,L]$, \begin{equation*}
\Bigl(\fint_{B_R} |\nabla(\hat{u}-u)|^2\Bigr)^{\frac{1}{2}} \le C(\frac{\ell}{L})^d(\frac{\rss}{L})^\beta,
	\end{equation*} where $C$ is a constant depending only on $d$, $\lambda$, $\varepsilon$ and $\hat{g}$.
\end{corollary}
\smallskip

 For the purpose of Algorithm \ref{alg:truealg} we now derive the explicit expression for $u_h$ from \eqref{eqn:lm1homphsm}. We start with $\tilde{u}_h$ defined in \eqref{intruhtilde}, which is an element of $Y_1^h(\ell)$. We make the ansatz \begin{equation}\label{eqn:multipole}
    u_h=\tilde{u}_h+\xi_i\partial_i G_h +\sum_{(i,j)\in \mathcal{J}} c_{ij}\partial_{ij}G_h \ \mod{Y_3^h(\ell)},
\end{equation}
where we recall the index set $\mathcal{J}$ is defined in \eqref{eqn:calJ}. This ansatz is motivated by the fact that $\{\partial_i G_h\}_{i=1,2,3}$ is a basis of the $3$-dimensional space $Y_1^h(\ell)/Y_2^h(\ell)$ and $\{\partial_{ij} G_h\}_{(i,j)\in\mathcal{J}}$ a basis of the $5$-dimensional space $Y_2^h(\ell)/Y_3^h(\ell)$, so that $\xi_i\partial_i G_h +\sum_{(i,j)\in \mathcal{J}} c_{ij}\partial_{ij}G_h$ is a general element of $Y_1^h(\ell)/Y_3^h(\ell)$.

 By equations \eqref{eqn:intrbaseq} and \eqref{intruhtilde} and integration by parts, we have for all $v_h\in X_2^h$,
\begin{equation}\label{eqn:biluva}
(u,v)_a=\int g \cdot \nabla (1+\phi_i\partial_i+\psi_{ij}\partial_{ij})v_h, \ \text{ and } \ (\tilde{u}_h,v_h)_{a_h}=\int g \cdot \nabla v_h.
\end{equation}
If $v_h\in X_0^h$ (i.e. $v_h$ is a constant) then both expressions vanish. Moreover, comparing \eqref{eqn:biluva} with \eqref{eqn:lm1homphsm}, we obtain the following identity for the multipole correction of $\tilde{u}_h$:  \begin{equation}\label{eqn:deltauh}
    (\xi_i\partial_i G_h +\sum_{(i,j)\in \mathcal{J}} c_{ij}\partial_{ij}G_h,v_h)_{a_h}= \int g\cdot  \nabla (\phi_i\partial_i+\psi_{ij}\partial_{ij})v_h \ \mbox{ for }v_h\in X_2^h / X_0^h.
\end{equation} 
Note that for $v_h\in X_2^h$, \begin{equation*}
    (\partial_i G_h, v_h)_{a_h}=-\int \nabla \eta \cdot (v_h a_h \nabla \partial_i G_h-\partial_i G_h a_h \nabla v_h) = \partial_i v_h(0),
\end{equation*}and similarly $(\partial_{ij} G_h,v_h)_{a_h} =-\partial_{ij} v_h$. Therefore, substituting this into \eqref{eqn:deltauh}, we obtain \begin{equation}\label{eqn:impvh}
    \xi_i \partial_i v_h(0) - \sum_{(i,j)\in \mathcal{J}} c_{ij} \partial_{ij} v_h = \int g\cdot  \nabla (\phi_i\partial_i+\psi_{ij}\partial_{ij})v_h.
\end{equation} Thus, choosing  $v_h=x_i$ for $i=1,2,3$ in \eqref{eqn:impvh}  we obtain \[\xi_i=\int g\cdot \nabla\phi_i,\] which is consistent with the expression obtained in \cite{lu2018optimal}. To determine the coefficients $c_{ij}$, we test \eqref{eqn:impvh} with $5$ linearly independent $a_h$-harmonic polynomials $v_{h,ij}$ that are homogeneous of degree $2$, and obtain linear equations of $c_{ij}$. Choosing the basis \begin{equation}\label{eqn:harmpol} v_{h,ij}=(1-\frac{1}{2}\delta_{ij})(x_ix_j-\frac{a_{hij}}{a_{h11}}x_1^2) \ \mbox{ for } (i,j)\in \mathcal{J},\end{equation} leads to the explicit formula which highly resembles \eqref{eqn:cijlt} in Algorithm \ref{alg:truealg}\begin{equation}\label{eqn:cij}
    c_{ij}=-\int g\cdot \nabla \Bigl(\phi_k\partial_k v_{h,ij}+(2-\delta_{ij})(\psi_{ij}-\dfrac{a_{hij}}{a_{h11}}\psi_{11})\Bigr) .
\end{equation}
Summing up, the function $u_h$ that captures the correct multipole behavior is given by \begin{equation} \label{eqn:effectivequadp}
u_h=\tilde{u}_h+ (\int g\cdot \nabla \phi_i)\partial_i G_h +c_{ij}\partial_{ij} G_h,
\end{equation} with $c_{ij}$ given by \eqref{eqn:cij}. This motivates Algorithm~\ref{alg:fake} to obtain an approximation $\hat{u}$ to the solution of \eqref{eqn:intrbaseq}.
\begin{algorithm}
\caption{Idealized algorithm \label{alg:fake}}
\begin{algorithmic}[1]
\State Solve \eqref{intrphi} for first-order correctors $\phi_i$.
\State Determine the homogenized coefficients $a_h$ via \eqref{intrhomcoeff}.
\State Solve \eqref{intruhtilde} for $\tilde{u}_h$ on $\partial Q_L$ by $\tilde{u}_h = \int G_h*(\nabla\cdot g)$.
\State Solve \eqref{eqn:intrsig} for first-order flux correctors $\sigma_{ijk}$ and \eqref{eqn:2ndcordef} for second-order correctors $\psi_{ij}$.
\State Obtain $u_h$ via \eqref{eqn:effectivequadp}.
\State Solve \eqref{eqn:coruhat} for $\hat{u}$, which is the approximation we desire.
\end{algorithmic}
\end{algorithm}

 Algorithm \ref{alg:fake} is, however, not computationally practical since several quantities like $\phi_i,\sigma_i,\psi_{ij}$ still require solving a whole-space problem like \eqref{intrphi} and thus knowledge of a realization of $a$ outside of $Q_{2L}$. Fortunately, we can replace each of these quantities in Algorithm \ref{alg:fake} by a computable surrogate with a small approximation error. This leads to Algorithm \ref{alg:truealg} and the error is only affected by a multiplicative constant. This is a consequence of the following Proposition \ref{prop:apriori}, which allows us to pass from Corollary \ref{cor:multipole} to the error estimates of Algorithm \ref{alg:truealg} in Theorem \ref{thm:mainthm}:
\begin{proposition}\label{prop:apriori}
	Under the same assumptions as in Proposition \ref{prop:effectivemultipole}, for $\sqrt{T}=L^{1-\varepsilon}\ge \ell \ge \rss$, define $\phi_{T}^{(L)},\sigma_{T}^{(L)}$ and $\psi_{T}^{(L)}$ as in \eqref{eqn:phiTL}, \eqref{eqn:algsigma} and \eqref{eqn:2ndcorapprox}. We assume they are good approximations of $\phi,\sigma,\psi$, in the sense of \footnote{Here and for the rest of the paper the notation ``$(f_1,f_2)\lesssim X$" (when obviously $X\ge 0$) means $f_1 \lesssim X$ and $f_2 \lesssim X$.}
	\begin{equation}\label{eqn:pp1assump}
	    \Bigl(\fint_{B_R} \bigl|\bigl(\sqrt{T}\nabla(\phi_{T}^{(L)}-\phi),\phi_T^{(L)}-\phi,\nabla(\psi_{T}^{(L)}-\psi)\bigr)\bigr|^2\Bigr)^\frac{1}{2} \le  \sqrt{T}(\dfrac{\rss}{L})^\beta \ \mbox{ for }  R\in \bigl\{\ell, \dfrac{5}{4}L\bigr\},
	\end{equation}  and \begin{equation}\label{eqn:pp1psiTgrowth}
	    (\fint_{B_L} |\psi_{T}^{(L)}-\psi|^2)^\frac{1}{2} \le  L^2(\dfrac{\rss}{L})^{\beta}.
	\end{equation}Define $a_h^{(L)}$ as in \eqref{eqn:algahL}, and we assume it is a good approximation of $a_h$ in the sense of \begin{equation}\label{eqn:pp1approxah}
	    |a_h^{(L)}-a_h| \le  (\dfrac{\rss}{L})^\beta.
	\end{equation} Given $\hat{g}$, let $u$ be defined as in \eqref{eqn:intrbaseq} and $u^{(L)}$ in \eqref{eqn:finalapprox}. Then we have for any $R\in [\rss,L]$,\begin{equation*}
	\Bigl(\fint_{B_R} |\nabla(u^{(L)}-u)|^2\Bigr)^{\frac{1}{2}} \le C(\frac{\ell}{L})^d(\frac{\rss}{L})^\beta,
	\end{equation*} with a constant $C$ of the same type as in Theorem \ref{thm:mainthm}.
\end{proposition}

\subsection{Massive Approximation of Correctors}\label{subsec:mascor}
We start with introducing the notations for this section and beyond. Given a length scale $R>0$, we define the exponential averaging function \[\eta_{R}(x):=\frac{c_d}{R^d}\exp(-\frac{|x|}{R})\] with the constant $c_d$ such that $\int_{\R^d} \eta_R(x)\ud x=1$. We also define the Gaussian \[G_R(x):=\dfrac{1}{(2\pi R^2)^\frac{d}{2}}\exp(-\frac{|x|^2}{2R^2}).\] For any function $f$, we use 
\begin{equation}\label{eqn:gaussconv}
f_R:=f*G_R
\end{equation}
to denote the convolution of $f$ with $G_R$. For any $s>0$, we define the following norm for a random variable $F$ that quantifies its tail:  
\begin{equation}\label{eqn:defsnorm}
\lVert F\rVert_s:=\inf\Bigl\{M\ge 0: \bigl\langle \exp\bigl((\dfrac{|F|}{M}+c)^s\bigr) \bigr\rangle-\exp(c^s)\le 1\Bigr\} \qquad \text{ with } \, c=\left\{\begin{aligned} & \Big(\frac{1-s}{s}\Big)^\frac{1}{s} & s\in(0,1), \\ & 0 & \text{otherwise.} \end{aligned} \right.
\end{equation}
Here the constant $c$ is chosen such that the function $[0,\infty) \ni x\mapsto \exp\bigl((x+c)^s\bigr)-\exp(c^s)$ is convex, which by Jensen's inequality makes $\lVert \cdot \rVert_s$ a norm. Therefore, if $\|F\|_s<\infty$, then $\langle \exp(r|F|^{\tilde{s}})\rangle<\infty$ for $r<\frac{1}{\|F\|_s}$ and $\tilde{s} \in (0,s]$.

We now illustrate why we may expect $\phi_T^{(L)}$ and $\psi_T^{(L)}$ to be good approximations of $\phi$ and $\psi$ in the sense of \eqref{eqn:pp1assump}. Let $\phi_T$ satisfy the equation \begin{equation}\label{eqn:intrphiTeqn}
\dfrac{1}{T}\phi_T-\nabla\cdot a(e+ \nabla\phi_T )=0.
\end{equation}
The massive corrector $\phi_T$, which was first considered in the early works of \cite{papanicolaou1979boundary,yurinskii1986averaging} and then \cite{gloria2011optimal,gloria2012optimal,gloria2017quantitative}, is an approximation of $\phi$ that has the advantage of being defined deterministically and it is automatically stationary. Indeed, in the class \begin{equation*}
    \sup_x \fint_{B_1(x)} (\phi_T^2 + |\nabla \phi_T|^2)<\infty
\end{equation*}
there exists a unique solution to \eqref{eqn:intrphiTeqn}.  The massive corrector serves as the bridge between $\phi$ and $\phi_T^{(L)}$, and we will show its closeness to both of them. The following Proposition \ref{prop:phismallscale} shows that if we choose the length scale $\sqrt{T}$ to be close to $L$, then the estimates \eqref{eqn:pp1assump} on $\phi-\phi_T$ hold with high probability.
\begin{proposition}\label{prop:phismallscale} For $d\ge 3$ and $\sqrt{T}\ge 1$, we have\footnote{We use ``$\lVert F\rVert_{s-} \lesssim X$" hereafter to denote for any $s'\in (0,s)$, $\lVert F\rVert _{s'} \lesssim_{s'} X$.}
 \begin{equation} \label{eqn:pp4eq2}
    \bigl\lVert(\int \eta_{\sqrt{T}} |\sqrt{T}\nabla(\phi_T-\phi)|^2)^\frac{1}{2} \bigr\rVert_{2-} \lesssim \sqrt{T}^{-\frac{1}{2}}.
\end{equation}
Moreover, there exists a random radius $\rr$ with \begin{equation}\label{eqn:estrstar}
	\lVert \rr\rVert_d \lesssim 1,
	\end{equation} and such that for all $R>0$, we have\footnote{Hereafter we use $I$ to denote indicator (characteristic) functions.}
\begin{equation} \label{eqn:pp4eq1}
\Bigl\lVert I(R \ge \rr)(\int \eta_{R} |\bigl(\sqrt{T}\nabla(\phi_T-\phi),\phi_T-\phi \bigr)|^2)^\frac{1}{2} \Bigr\rVert_{2-} \lesssim \sqrt{T}^{-\frac{1}{2}}.
\end{equation}
\end{proposition} 
\begin{remark}
    For $d>4$, the approximation error \eqref{eqn:pp4eq1} saturates at $\sqrt{T}^{-1}$. A better approximation which has error $\sqrt{T}^{1-\frac{d}{2}}$ is given in \cite[Theorem 3]{GO2015most} and arises from iterated Richardson extrapolation of $\phi_T$. We expect a similar strategy to work in the optimal approximation of higher order correctors.
\end{remark}
\begin{remark}
    The quantity $\rr$, also known as the minimal radius, is the smallest scale on which the elliptic and parabolic $C^{0,1}$-estimates (See Lemma \ref{lem:parabolicmvp}) hold. Large scale regularity was first considered in \cite{avellaneda1987compactness} for periodic homogenization (see also the monograph \cite{shen2018periodic}), and \cite{armstrong2016quantitative, gloria2014regularity} introduced the random variable $\rr$ in the stochastic setting. For any fixed $\delta \in (0,1)$ we define $\rr$ as in \cite{GO2015most}\footnote{Note that the definitions in \cite{armstrong2016quantitative, gloria2014regularity} are slightly different}: \begin{equation}\label{eqn:regminrad}
    \rr=\rr(0):=\inf\Bigl\{ r\ge 1  \mid \forall R\ge r,  \ \dfrac{1}{R}\Bigl(\fint_{B_R} |(\phi,\sigma)-\fint_{B_R} (\phi,\sigma)|^2\Bigr)^\frac{1}{2} \le \delta \Bigr\}.
\end{equation} 
Comparing the above \eqref{eqn:regminrad} with \eqref{intrphigr}, we observe $\rss\ge \rr$. The stochastic estimate \eqref{eqn:estrstar} is proved in \cite[Corollary 6]{GO2015most}. For slightly different definitions of $\rr$ under various probabilistic settings, we refer to \cite[Theorem 1.1]{armstrong2016quantitative} and \cite[Theorem 1]{gloria2014regularity}  for corresponding stochastic estimates.
\end{remark}
On the other hand, $\phi_T$ near the origin can be well-approximated by a function that only depends on $a$ through its restriction to the finite domain $Q_L$, which is exactly achieved by the function $\phi_{T}^{(L)}$ defined in \eqref{eqn:phiTL}. Though the Dirichlet boundary conditions break the stationarity of $\phi_{T}^{(L)}$, deterministic methods are sufficient to prove that the difference  $\phi_T-\phi_{T}^{(L)}$ is sub-algebraically small in $\frac{\sqrt{T}}{L}$ for $\sqrt{T}\ll L$. The exact statement of this is deferred to Proposition \ref{prop:truncation}. 
\smallskip

As a comparison, if we consider the Dirichlet approximation as \cite[(18)]{lu2018optimal}
\begin{equation*}
 -\nabla  \cdot a \nabla \phi_{i}^{(L)} =\nabla\cdot ae_i \,  \mbox{ in }Q_{2L}, \hspace{0.3in} \phi_{i}^{(L)}=0 \, \mbox{ on }\partial Q_{2L}, 
\end{equation*}
then it is unclear if one could prove anything stronger than what is proved in \cite{lu2018optimal} 
\begin{equation*}
    \Big(\fint_{Q_L}|\nabla(\phi^{(L)}-\phi)|^2\Big)^\frac{1}{2} \lesssim L^{-1},
\end{equation*}
which is the desired scaling for $d=2$ but insufficient for the CLT scaling required in \eqref{eqn:pp1assump} for $d=3$. The bottleneck is that $\phi^{(L)}$ is not stationary and it is thus unclear how probabilistic tools can be applied here. While one may use appropriate oversampling techniques to obtain desirable approximations using $\phi^{(L)}$, in this work we opt for our approximation $\phi_T^{(L)}$ in \eqref{eqn:phiTL} since we could prove the required approximation bounds and it is not more numerically difficult to compute than $\phi^{(L)}$. 

Equipped with estimates on $\phi$ we may derive estimates for the approximation of homogenized coefficients $a_h$. In analogy to \eqref{intrhomcoeff} we introduce the modified flux \begin{equation}\label{eqn:defqsubT} q_T:=a(e+\nabla \phi_T).\end{equation} For any smooth weight function $\omega(x)$ supported in the unit ball satisfying $\int_{\R^d} \omega=1$, and rescaled as $\omega_L(x)=\frac{1}{L^d}\omega(\frac{x}{L})$, $\int \omega_L q_T$ is a good approximation of $a_h$, which establishes that the assumption \eqref{eqn:pp1approxah} of Proposition \ref{prop:apriori} holds with high probability: 
\begin{lemma}\label{lem:ahahL} Let $L \ge \sqrt{T} \ge 1$. Then for $d\ge 3$ we have
\begin{equation}\label{eqn:ppahahL}
     \lVert a_he_i-\int \omega_L q_{i,T} \rVert_{2-}  \lesssim \sqrt{T}^{-\frac{3}{2}}.
 \end{equation}
\end{lemma}
Using the same procedure, we can also approximate the first-order flux corrector $\sigma$ and the second-order corrector $\psi$ using their massive counterparts, denoted as $\sigma_T$ and $\psi_T$ respectively, which are defined by the equations \begin{align}
& \dfrac{1}{T} \sigma_{ijk,T}-\Delta\sigma_{ijk, T}= \nabla \cdot (q_{ik,T}e_j-q_{ij,T}e_k),  \label{eqn:moifiedflux} \\
& \dfrac{1}{T}\psi_{ij,T}-\nabla\cdot a \nabla\psi_{ij,T}=\nabla\cdot (\phi_{i,T} a-\sigma_{i,T})e_j.  \label{eqn:2ndcorinitapprox}
\end{align}
Similar to $\phi_T$, both $\sigma_T$ and $\psi_T$ are well-defined stationary fields, and can be approximated by functions that only depend on $a$ through $a|_{Q_{2L}}$, which are $\sigma_T^{(L)}$ and $\psi_T^{(L)}$ defined in \eqref{eqn:algsigma} and \eqref{eqn:2ndcorapprox} respectively. The following proposition gives estimates on $\psi_T$, establishing that the assumptions \eqref{eqn:pp1assump} and \eqref{eqn:pp1psiTgrowth} hold with high probability. \begin{proposition}\label{prop:psibounds} Let $\rr$ be the same random radius as defined in Proposition \ref{prop:phismallscale}. Then for $d\ge 3$ and $\sqrt{T}\ge 1$, $R>0$, \begin{align}
        & \lVert I(R\ge \rr)(\int \eta_{R} |\psi_T|^2)^\frac{1}{2}\rVert_{1-}  \lesssim \sqrt{T}^\frac{1}{2}, \label{eqn:pp6eq2}\\
	& \bigl\lVert I(R\ge \rr)\bigl(\int \eta_{R}  |\nabla (\psi-\psi_{T})|^2\bigr)^\frac{1}{2}  \bigr\rVert_{1-}  \lesssim \sqrt{T}^{-\frac{1}{2}}. \label{eqn:pp6eq4}
\end{align}
\end{proposition}
\noindent The massive correctors $\phi_T$ and $\psi_T$ can both be approximated by functions that depend on $a$ only through its restriction to a finite domain, with an error smaller than any power of $\frac{\sqrt{T}}{L}$ for $\sqrt{T}\ll L$, which is the result of the following Proposition \ref{prop:truncation}:
\begin{proposition}\label{prop:truncation}
    Let $L\ge \sqrt{T}\ge 1$. Let $\phi_T^{(L)},\sigma_T^{(L)}$ and $\psi_T^{(L)}$ be defined through \eqref{eqn:phiTL}, \eqref{eqn:algsigma} and \eqref{eqn:2ndcorapprox}, then for any $p<\infty$,  \begin{align}
    & \Bigl(\fint_{Q_{\frac{7}{4} L}} \bigl|\bigl(\sqrt{T}\nabla (\phi_T-\phi_{T}^{(L)}), \phi_T-\phi_T^{(L)}\bigr)\bigr|^2\Bigr)^\frac{1}{2}\lesssim_p (\frac{\sqrt{T}}{L})^p, \label{eqn:phitruncation} \\ & \Bigl(\fint_{Q_{\frac{5}{4}L}}\bigl|\bigl(\sqrt{T}\nabla (\psi_T-\psi_{T}^{(L)}),\psi_T-\psi_T^{(L)}\bigr)\bigr|^2\Bigr)^\frac{1}{2}\lesssim_p \sqrt{T}(\frac{\sqrt{T}}{L})^p.  \label{eqn:psiapprox}
    \end{align} 
\end{proposition}
 The reason we cannot prove \eqref{eqn:phitruncation} with $Q_{2L}$ is that the Dirichlet boundary conditions generate a boundary layer where $\phi_T$ may be dramatically different from $\phi_T^{(L)}$. Therefore near $\partial Q_{2L}$, $\phi_T^{(L)}$ is not trustworthy and should not be used for the computation of $\sigma_T^{(L)}$, which is the reason why in Algorithm \ref{alg:truealg} the domain of computation shrinks when computing $\sigma_T^{(L)}$ and, for the same reason, shrinks further for $\psi_T^{(L)}$. 

We finally present Proposition \ref{prop:psirss}, which is the main ingredient of the proof for the stochastic estimate on $\rss$ in \eqref{eqn:thm1rss}, and, together with the estimate \eqref{eqn:estrstar} on $\rr$, bounds the probability of the event $\rss \ge R$ for any large $R$: \begin{proposition}\label{prop:psirss}
    Let $\rr$ be the same random radius as defined in Proposition \ref{prop:phismallscale}, and denote $(\nabla\psi,\nabla\Psi)_R$ as the convolution of the two functions $(\nabla \psi, \nabla \Psi)$ with the Gaussian kernel $G_R$ as in \eqref{eqn:gaussconv}, then for $d\ge 3$ and $R\ge 1$, \begin{align} & \bigl\lVert I(R\ge \rr)\bigl(\fint_{B_R} |\nabla (\psi, \Psi)|^2\bigr)^\frac{1}{2}\bigr\rVert_{1-}  \lesssim 1, \label{eqn:pp6eq1}\\ & \lVert I(R\ge \rr)\bigl(\nabla\psi,\nabla\Psi\bigr)_R\rVert_{1-}  \lesssim R^{1-\frac{d}{2}}. \label{eqn:pp6eq3} \end{align}
\end{proposition}

\section{Numerical Example}

For the numerical test, we will consider a discrete elliptic equation on $\Z^3$, so that we do not need to worry about error due to discretization. To set up the elliptic problem, we say the points $x,y\in\Z^3$ are neighbors if $\lVert x-y\rVert_{\ell_1}=1$, and draw an edge between $x$ and $y$ if they are neighbors. Denote $\mathbb{B}$ as the set of (undirected) edges, $\{\textbf{e}_1,\textbf{e}_2,\textbf{e}_3\}$ as the canonical basis in $\Z^3$, and $(a_e)_{e\in \mathbb{B}}$ as the random field. The discrete gradient is defined as $$\nabla f(x)=(f(x+\textbf{e}_1)-f(x),f(x+\textbf{e}_2)-f(x),f(x+\textbf{e}_3)-f(x))$$ and the divergence of $F=(F_1,F_2,F_3)$ is defined as \[\nabla\cdot F(x)=\sum_{i=1}^{3}(F_i(x)-F_i(x-\textbf{e}_i)).\] 

We consider the discrete equation, with the above notations, still in the form
\[-\nabla\cdot a\nabla u=\nabla\cdot g.\]
The coefficient field $a=a_e$ on edges $e \in \mathbb{B}$ are i.i.d.{} random matrices with values $\mathrm{Id}$ and $9\,\mathrm{Id}$ with probability $1/2$ each. 
For the right hand side, we take some function $f(x)$ compactly supported in the box $\{-1, 0, 1\}^3$ with average zero, so that there exists some vector valued function $g(x)$ such that $f=\nabla \cdot g$ for a function $g$ supported in the slightly larger box $Q_2$. 

We compare Algorithm~\ref{alg:truealg} with three other algorithms: 
\begin{enumerate} 
\item Solving the equation \eqref{eqn:diropt} with zero Dirichlet boundary condition.
\item Solving the equation \eqref{eqn:diropt} with modified correctors but without dipole or quadruple corrections, i.e., the boundary condition given by \begin{equation}\label{eqn:numnopole}
u_{nc}^{(L)}=(1+\phi_{i,T}^{(L)}\partial_i+\psi_{ij,T}^{(L)}\partial_{ij})\tilde{u}_h^{(L)} \ \mbox{ on } \partial Q_L.
\end{equation}
\item Solving the equation \eqref{eqn:diropt} with boundary condition corrected up to first-order correctors and dipoles, which is the algorithm proposed in \cite{lu2018optimal}: \begin{equation}\label{eqn:numdipole}
u_{dp}^{(L)}=(1+\phi_{k,T}^{(L)}\partial_k)\bigl(\tilde{u}_h^{(L)}+(\int g \cdot \nabla \phi_{i,T}^{(L)}) \partial_i G_h^{(L)}\bigr) \ \mbox{ on } \partial Q_L.
\end{equation}\end{enumerate} 

\begin{figure}[!ht]
    \centering
    \includegraphics[width=0.55\textwidth]{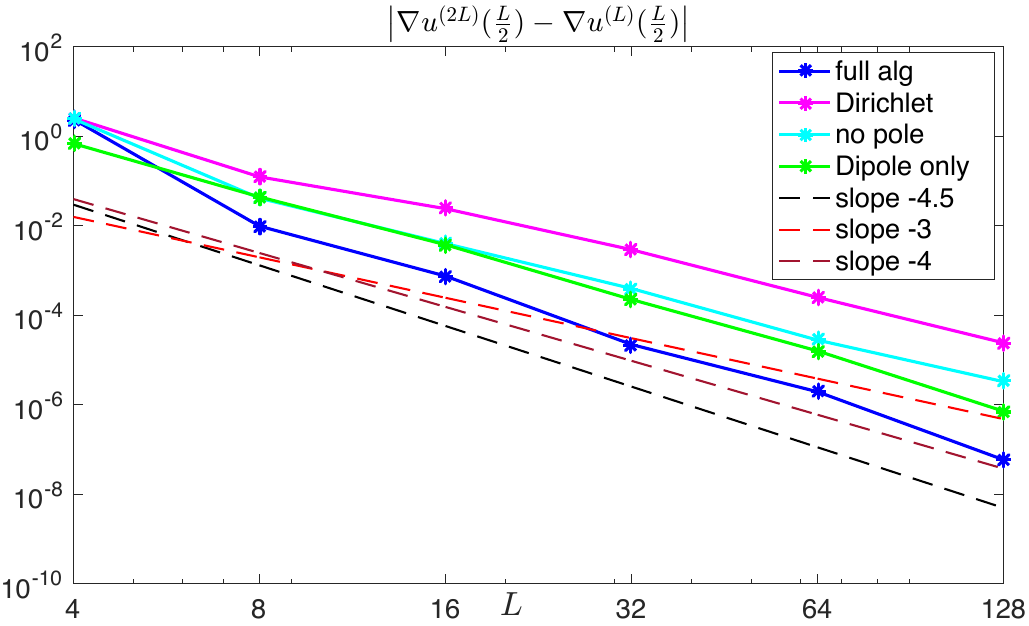}
    
    \includegraphics[width=0.55\textwidth]{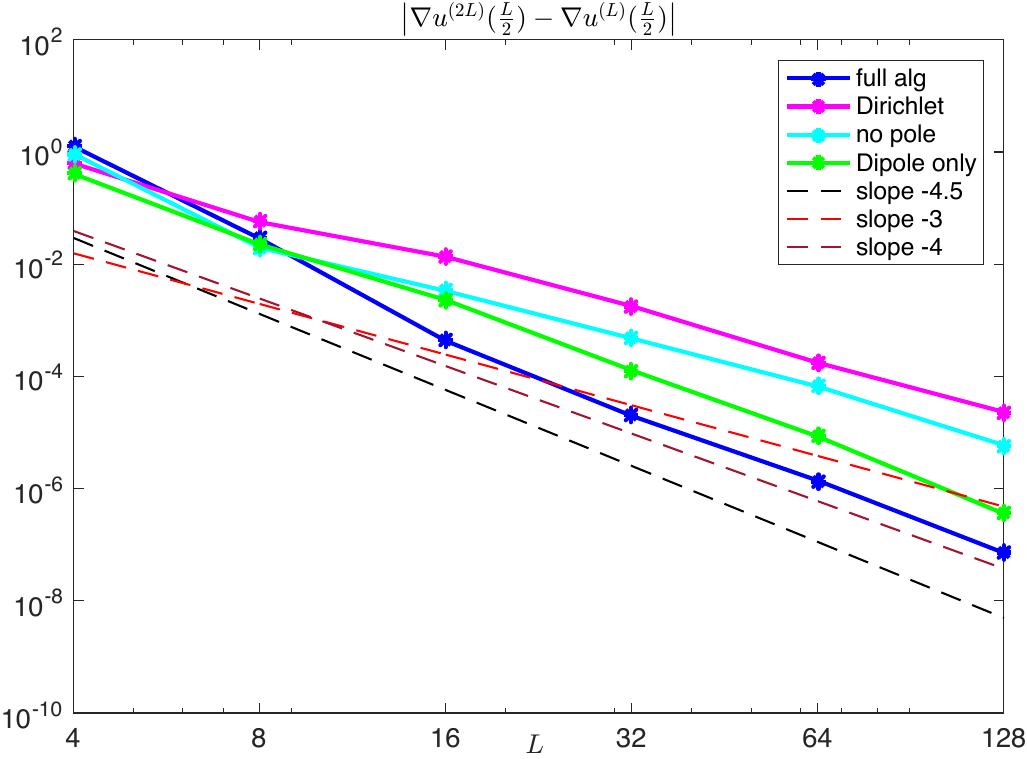}
    \caption{Numerical convergence rate of $\lvert \nabla u^{(2L)}(\frac{L}{2},\frac{L}{2},\frac{L}{2})-\nabla u^{(L)}(\frac{L}{2},\frac{L}{2},\frac{L}{2}) \rvert$ for the proposed Algorithm \ref{alg:truealg}, an algorithm with zero Dirichlet boundary condition, an algorithm without multipole corrections (defined in \eqref{eqn:numnopole}), and an algorithm with only dipole corrections (defined in \eqref{eqn:numdipole}).  The two figures correspond to two independent realizations of the random media and the same r.h.s. 
    \label{fig:fig1}}
\end{figure}

We compare the numerical rate of $\lvert \nabla (u^{(2L)}- u^{(L)})(\frac{L}{2},\frac{L}{2},\frac{L}{2})\rvert$ and plot it for various $L$ for all four algorithms. We would comment here that while our analysis is for the gradient averaged over a region, for simplicity we only compare the gradient at a single point. From Figure \ref{fig:fig1} we can observe that the Dirichlet algorithm and the no multipole algorithm both have convergence rates of $O(L^{-3})$, the dipole algorithm has a convergence rate of $O(L^{-4})$ while the proposed algorithm achieves  $O(L^{-4.5})$ convergence rate, which are consistent with our theoretical results. 

Let us remark that the reason we compare solutions at $(\frac{L}{2}, \frac{L}{2}, \frac{L}{2})$ is that, if the observation point is at the origin, then the dipole algorithm will achieve a better error rate of $O(L^{-4.5})$. The reason is that, if we compare the dipole algorithm and the full algorithm, \begin{equation*}
    u^{(L)}-u_{dp}^{(L)}=(1+\phi_{i,T}^{(L)}\partial_i+\psi_{ij,T}^{(L)}\partial_{ij})c_{kn}^{(L)}\partial_{kn}G_h^{(L)}\\+\psi_{ij,T}^{(L)}\partial_{ij}\bigl(\tilde{u}_h^{(L)}  +(\int  g \cdot\nabla \phi_{k,T}^{(L)})\partial_k G_h^{(L)}\bigr) \ \mbox{ on }\partial Q_L,
\end{equation*} then besides the term $c_{kn}^{(L)}\partial_{kn}G_h^{(L)}$, all other terms will result in a contribution of $O(L^{-4.5})$ or higher order at the origin. Moreover, in dimension $3$ the term $c_{kn}^{(L)}\partial_{kn}G_h^{(L)}$ is an odd function so its contributions will vanish near the origin, making the dipole algorithm behave similarly as the full algorithm.

\begin{figure}[ht]
    \centering
    \includegraphics[width=0.6\textwidth]{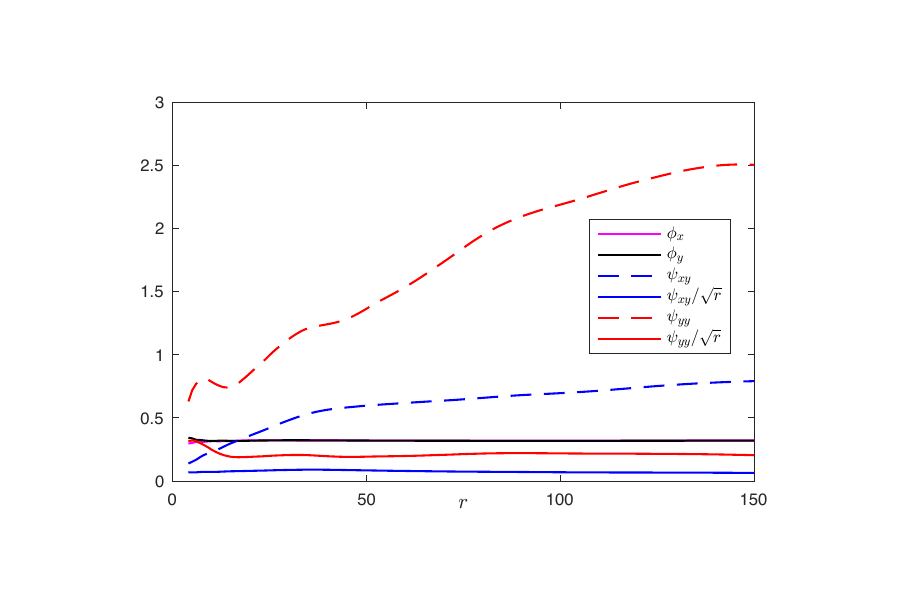}
    \includegraphics[width=0.6\textwidth]{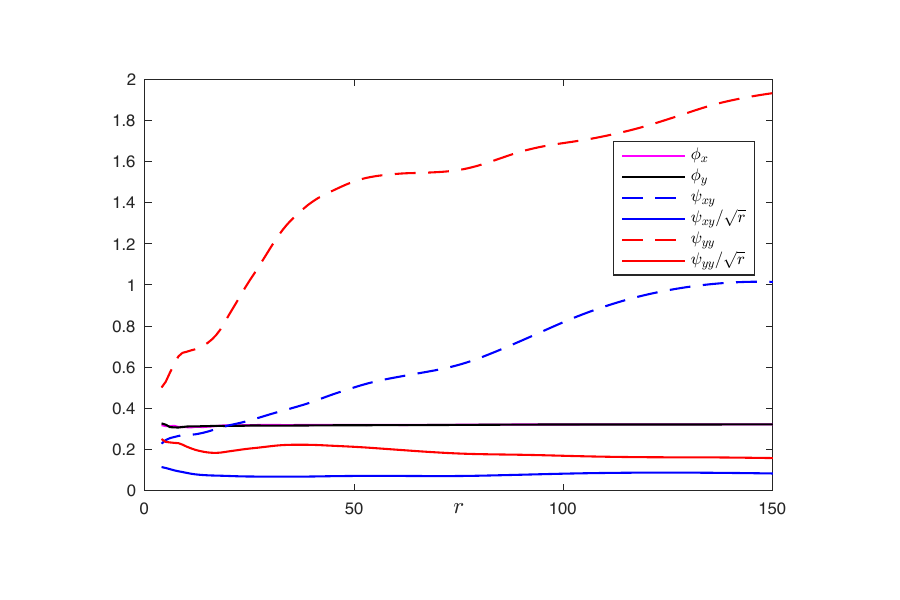}
    \caption{$L^2$-averages of $\phi$ and fluctuations of $\psi$. The two figures correspond to independent realizations of the random media. }
    \label{fig:growth}
\end{figure}

We also numerically compute $(\fint_{Q_r}|\phi_T^{(L)}|^2)^\frac{1}{2}$ and $(\fint_{Q_r}|\psi_T^{(L)}-\fint_{Q_r} \psi_T^{(L)}|^2)^\frac{1}{2}$ for a variety of $r$. Figure \ref{fig:growth} indicates that the quantities $(\fint_{Q_r} |\phi_T^{(L)}|^2)^\frac{1}{2}$ and $\frac{1}{\sqrt{r}}(\fint_{Q_r}|\psi_T^{(L)}-\fint_{Q_r} \psi_T^{(L)}|^2)^\frac{1}{2}$ are almost constants for all $r$, which is consistent with their growth estimates \eqref{intrphigr} and \eqref{intrpsiorig}. The figure also indicates numerically that $\rss$ should be of order $1$.

\section{Strategy of the Proof}\label{sec:lemmas}
\subsection{Parabolic Semigroup Representation of Correctors and Fluxes}
Our approach is based on the semigroup framework developed in \cite{GO2015most}, see also \cite[Chapter 9]{armstrong2017quantitative}. The central objects we study are the operators $S(t)$ and $\opS(t)$ related to $-\nabla\cdot a \nabla$. \begin{definition} The linear operators $S(t)$ and $\opS(t)$ are defined as follows: for an arbitrary vector field $g$, let $v$ solve the parabolic equation \begin{equation*}
         \partial_t v-\nabla\cdot a \nabla v=0 \ \mbox{ for }t>0, \ \ v(t = 0) =\nabla\cdot g,
\end{equation*} 
then \begin{equation}\label{eqn:cals}
        S(t)g :=v(t), \ \text{ and } \ \opS(t)g := g+a\int_0^t \ud\tau \nabla v(\tau). \end{equation} We also define the operator $\bar{S}(t)$ similarly to $S(t)$ with $a$ replaced by constant coefficient $\mathrm{Id}$, and $\opS^h(t)$ similarly to $\opS(t)$ with $a$ replaced by $a_h$. \end{definition} 
    Roughly speaking, $\opS(t)g$ is the flux accumulated from time $0$ to $t$ with initial condition $\nabla \cdot g$. As opposed to $S(t)$, $\opS(t)$, which is the same as the operator $S_{0\to t}$ defined in \cite{GO2015most}, does satisfy the semigroup property \cite[Lemma 2]{GO2015most} \begin{equation}\label{eqn:GOlm2}
    \opS(t_1)\opS(t_2) g=\opS(t_1+t_2)g.
\end{equation}Since the solution of the standard heat equation is the convolution of the initial condition with the heat kernel, which coincides with $G_{\sqrt{t}}$, $\bar{S}$ can be rewritten as \begin{equation}\label{eqn:barSgaussian}
        \bar{S}(t)g=(\nabla\cdot g)_{\sqrt{t}} \leftstackrel{\eqref{eqn:gaussconv}}{=} G_{\sqrt{t}}* (\nabla \cdot g).
    \end{equation}

The semigroup operators $S(t)$ and $\bar{S}(t)$ are essential since all correctors and their massive approximations can be represented using these operators, which is the building block for the estimates of correctors and fluxes in Propositions \ref{prop:phismallscale}--\ref{prop:truncation}. More precisely, these operators provide a resolution of the correctors in terms of quantities of controlled locality and amplitude. Below in Remark \ref{lem:psidecom} we collect the semigroup representation formulas for all correctors and fluxes. Here for simplicity of notation, we omit the indices, and introduce the vector product notation $\times$ so that $\sigma = \int_0^\infty \ud t_0 \bar{S}(t_0)\times q$ stands for $\sigma _{ijk}= \int_0^\infty \ud t_0 \bar{S}(t_0)(q_{ik}e_j-q_{ij}e_k)$.
\begin{remark}\label{lem:psidecom}
    The correctors $\phi, \sigma, \psi, \Psi$ and their massive approximations $\phi_T, \psi_T$ can be represented as follows 
    \begin{align}
 \phi & =\int_0^\infty \ud tS(t)ae, \label{eqn:expphi} \\
 \phi_T &= \int_0^\infty \ud t\exp(-\dfrac{t}{T}) S(t)ae, \label{eqn:expphiT}
    \\  \sigma  &=\int_0^\infty \ud t_0\bar{S}(t_0)\times q,\label{eqn:easydecomsig} \\
     \nabla \psi & =\int_0^\infty \ud t_0\int_0^\infty \ud t_1 \nabla S(t_0)(aS(t_1)ae-\bar{S}(t_1)\times ae) \nonumber \\ & \qquad - \int_0^\infty\ud t_0\int_0^\infty \ud t_1 \int_0^\infty \ud t_2 \nabla S(t_0)\bar{S}(t_1)\times a\nabla S(t_2)ae,\label{eqn:psidecom} \\
     \psi_T & =\int_0^\infty\ud t_0 \exp(-\dfrac{t_0}{T})\int_0^\infty \ud t_1\exp(-\dfrac{t_1}{T}) S(t_0)(aS(t_1)ae-\bar{S}(t_1)\times ae) \nonumber  \\ & \qquad - \int_0^\infty\ud t_0 \exp(-\dfrac{t_0}{T})\int_0^\infty \ud t_1 \exp(-\dfrac{t_1}{T}) \int_0^\infty \ud t_2 \exp(-\dfrac{t_2}{T}) S(t_0)\bar{S}(t_1)\times a\nabla S(t_2)ae,\label{eqn:psiTdecom} \\ \nabla \Psi &= \int_0^\infty \ud t \; \nabla \bar{S}(t) \times (a\nabla \psi+(a\phi-\sigma)e).\label{eqn:convPsi}
\end{align}
\end{remark}
We would like to comment here that \eqref{eqn:psidecom} and \eqref{eqn:convPsi} are formulated on the level of gradients, since only $\nabla (\psi,\Psi)$ are well-defined stationary random fields. To show the algebraic equivalence of the quantities, it suffices to show the r.h.s.\ of the equations satisfy the same elliptic equations as their counterparts on the l.h.s. The core argument we need is, suppose a function $w_T$ has the representation \begin{equation}\label{eqn:masparadecom} w_T=\int_0^\infty\ud t \exp(-\dfrac{t}{T})S(t)g,\end{equation} then it satisfies the massive equation $\dfrac{1}{T}w_T-\nabla\cdot a \nabla w_T=\nabla \cdot g.$ Indeed, \begin{align*}-\nabla\cdot a \nabla w_T &=-\int_0^\infty\ud t \exp(-\dfrac{t}{T}) \nabla \cdot a \nabla S(t)g = -\int_0^\infty \ud t\exp(-\dfrac{t}{T}) \partial_t S(t)g \\ &= S(0)g+\int_0^\infty \ud t \bigl( \partial_t \exp(-\dfrac{t}{T}) \bigr) S(t)g = \nabla \cdot g-\dfrac{1}{T}w_T.\end{align*} Similarly, if $w$ has the representation \begin{equation}\label{eqn:nomasparadecom}w=\int_0^\infty S(t)g\ud t,\end{equation} then $w$ satisfies $-\nabla\cdot a \nabla w=\nabla \cdot g.$ In particular, we obtain for $a=\mathrm{Id}$, the function  \begin{equation}\label{eqn:decomlap}
    w=\int_0^\infty \ud t(\nabla \cdot g)_{\sqrt{t}}.
\end{equation} satisfies $-\Delta w = \nabla \cdot g.$ Choosing $g=ae$ in \eqref{eqn:masparadecom} and \eqref{eqn:nomasparadecom} and we obtain the representations \eqref{eqn:expphiT} for $\phi_T$ and \eqref{eqn:expphi} for $\phi$ respectively, and choosing $g=q$ in \eqref{eqn:decomlap} yields \eqref{eqn:easydecomsig} for $\sigma$. The representations for $(\nabla\psi, \psi_T, \nabla \Psi)$ can be derived by applying the above arguments iteratively and we will not write the details here. To rigorously obtain the equivalences, in view of uniqueness of stationary correctors, it suffices to show that the r.h.s.\ integrals converge and therefore represent stationary functions, which is done in the proofs of \eqref{eqn:lm4eq1} for $(\phi,\sigma)$ and \eqref{eqn:pp6eq1}, \eqref{eqn:pp6eq2} for $(\nabla \psi, \nabla \Psi, \psi_T)$ (such result for $\phi_T$ is classic). 

\subsection{CLT-Cancellations and Propagation of Locality}
The next important notion is the so-called ``relative approximate locality". \begin{definition} Let $p>\frac{d}{2}$ be fixed, and let $g$ and $\bar{g}$ denote two stationary random fields. We say $g$ is approximately local on scale $r\ge 1$ relative to $\bar{g}$ if, for two realizations $a$ and $\tilde{a}$ satisfying $a=\tilde{a}$ in $B_{2R}$ for some $R\ge r$,  \begin{equation}\label{eqn:apploc}
    (\fint_{B_R} |g(a)-g(\tilde{a})|^2)^\frac{1}{2} \le (\dfrac{r}{R})^p \int \eta_{R} \big(\bar{g}(a)+\bar{g}(\tilde{a})\big).
\end{equation} \end{definition} Roughly speaking, the random field $g|_{B_R}$ ``essentially depends on $a$ only through $a|_{B_r}$ up to an error term $(\frac{r}{R})^p\bar{g}$''. This locality is at the basis of a CLT scaling, see Lemma \ref{lem:GOlm13}.

In order to estimate higher-order correctors and their fluxes, it is important to estimate the operator $S(T)$ acting on arbitrary $g$. More precisely, we extend the results in \cite{GO2015most} on approximate locality from the special and obviously local initial data $ae$ to more general initial data $g$ that are only approximately local. Our goal is to derive stochastic bounds and approximate locality properties for $(S(T)g, \nabla S(T)g,\opS(T)g)$, given $g$ approximately local on scale $r_0\ge 1$ relative to $\bar{g}$. In addition, while $ae$ is obviously bounded, we only assume stochastic integrability of $g$ and $\bar{g}$, in the sense of controlling $\lVert g\rVert_{s}$ and $\lVert \bar{g} \rVert_{s}$ for some $0<s\le 2$.

In the range $\sqrt{T}\le r_0$, $S(T)g$ does not benefit from stochastic cancellations through parabolic propagation of $g$. Therefore we can do no better than the following consequence of a deterministic estimate (see Lemma \ref{lem:GOlm1} below) \begin{equation}\label{eqn:sc5arg}\Bigl\lVert\Bigl(\int \eta_{r_0}\bigl|\bigl(T\nabla S(T)g,\sqrt{T}S(T)g,\opS(T)g\bigr)\bigr|^2\Bigr)^\frac{1}{2}\Bigr\rVert _{s} \lesssim \lVert(\int \eta_{r_0} |g|^2)^\frac{1}{2}\rVert _{s}.\end{equation}
The scenario is more subtle in the regime $\sqrt{T}\ge r_0$. While $S(T)g$ now benefits from stochastic cancellations, for which the ratio of the locality scale $r_0$ of $g$ and the parabolic scale $\sqrt{T}$ matters, it also suffers an increase of the locality scale to $\sqrt{T}$, as well as a loss of stochastic integrability to $\frac{2s}{s+2}-$. The first and third aspects are captured by Proposition \ref{prop:applocqt}, the second by Proposition \ref{prop:Stglocality}.
\begin{proposition}\label{prop:applocqt}
    Suppose $g$ and $\bar{g}$ satisfy \eqref{eqn:apploc} for $r_0\ge 1$. Then for all $R\ge 1, \ T\ge 0$, we have\footnote{Here we use $r_1 \vee r_2$ and $r_1 \wedge r_2$ to denote $\max\{r_1,r_2\}$ and $\min\{r_1,r_2\}$, respectively.} \begin{equation}\label{eqn:pp5eq1}
        \Bigl\lVert \Bigl(\opS(T)g-\langle \opS(T)g \rangle,\int_0^T \ud t \nabla S(t)g \Bigr)_R\Bigr\rVert_{\frac{2s}{s+2}-} \lesssim (\frac{r_0}{R})^\frac{d}{2}\Bigl(\lVert(\int \eta_{r_0} |g|^2)^\frac{1}{2}\rVert _s+\lVert\bar{g}\rVert _s\Bigr)
    \end{equation}
    and 
    \begin{equation}\label{eqn:pp5eq2}
        \Bigl\lVert\Bigl(\int \eta_{\sqrt{T}}\bigl|\bigl(T\nabla S(T)g,\sqrt{T}S(T)g\bigr)\bigr|^2\Bigr)^\frac{1}{2}\Bigr\rVert _{\frac{2s}{s+2}-} \lesssim (1\wedge\dfrac{r_0}{\sqrt{T}})^\frac{d}{2}\Bigl(\lVert(\int \eta_{r_0} |g|^2)^\frac{1}{2}\rVert _s+\lVert\bar{g}\rVert _s\Bigr).
    \end{equation}  
\end{proposition}
\begin{remark}\label{rmk:barSintg}
    The same estimates hold with $S$ replaced by $\bar{S}$, with stochastic integrability improved from $\frac{2s}{s+2}-$ to $s$, see \cite[Lemma 16]{GO2015most}.
\end{remark}
\begin{proposition}\label{prop:Stglocality}
      Suppose $g$ and $\bar{g}$ satisfy \eqref{eqn:apploc} for $r_0\ge 1$. For all $T\ge 0$, $(T\nabla S(T)g, \sqrt{T}S(T)g, \opS(T)g)$ is approximately local on scale $r_0\vee \sqrt{T}$ relative to some stationary $\bar{G}$ with \begin{equation}\label{eqn:bargbound} \lVert\bar{G}\rVert _{s_0} \lesssim (1\wedge \dfrac{r_0 }{\sqrt{T}})^\frac{d}{2}\Bigl(\lVert(\int \eta_{r_0} |g|^2)^\frac{1}{2}\rVert _s+\lVert\bar{g}\rVert _s\Bigr).
    \end{equation}Here $s_0=\frac{2s}{s+2}-$ when $\sqrt{T}> r_0$ and $s_0=s$ when $\sqrt{T}\le r_0$.
\end{proposition}
    Propositions \ref{prop:applocqt} and \ref{prop:Stglocality} tell us that if initially $ \lVert(\int \eta_{r_0} |g|^2)^\frac{1}{2}\rVert _s$ and $\lVert\bar{g}\rVert _s$ have the same upper bound, then so will $\lVert(\int\eta_{(r_0\vee \sqrt{T})}|(T\nabla S(T)g,\sqrt{T} S(T)g,\opS(T)g)|^2)^\frac{1}{2}\rVert _{\frac{2s}{s+2}-}$ and $\lVert\bar{G}\rVert _{\frac{2s}{s+2}-}$, so that we can essentially ``sweep under the rug'' the quantity quantifying the locality since it has the same size as the original quantity. This is convenient for estimating multiple time integrals. 
     
In particular, for $g=ae$, in which case we may set $r_0=1, \, \bar{g}=0$ and $s=\infty$, we recover the following stochastic estimates on $S(t)ae$ proved in \cite{GO2015most}, which are special cases of Propositions \ref{prop:applocqt} and \ref{prop:Stglocality}.  
\begin{lemma}\label{cor:GOCor4}  (\cite[Corollary 4]{GO2015most})
    For all $T\ge 0$, 
    \begin{equation}\label{eqn:semigroup}
\Bigl\lVert\bigl(\int \eta_{\sqrt{T}} |\bigl(T\nabla S(T)ae,\sqrt{T}S(T)ae\bigr)|^2\bigr)^\frac{1}{2}\Bigr\rVert_{2-} \lesssim (1\wedge \frac{1}{\sqrt{T}})^{\frac{d}{2}}.
\end{equation}
Moreover, $(T\nabla S(T)ae, Ta\nabla S(T)ae, \sqrt{T}S(T)ae,\opS(T)ae)$ is approximately local on scale $1 \vee \sqrt{T}$ relative to some stationary $\bar{g}$ with \begin{equation}\label{eqn:expressionbarg}
\lVert \bar{g}\rVert_{2-}\lesssim (1\wedge \frac{1}{\sqrt{T}})^{\frac{d}{2}}.
\end{equation}
\end{lemma}

With Lemma \ref{cor:GOCor4} and the representations \eqref{eqn:expphi}, \eqref{eqn:expphiT} of $\phi$ and $\phi_T$, we see why $\phi_T$ is a good approximation of $\phi$: in the range $t\gg T$, the contribution of $S(t)ae$ to $\phi_T$ is exponentially small, and the locality of $S(t)ae$ in conjunction with the finite range condition on $a$ make the contribution to $\phi$ small as well, which (almost) matches the bound \eqref{eqn:pp1assump}; in the range $t\lesssim T$, we can use $1-\exp(-\frac{t}{T})\lesssim \frac{t}{T}$ and the stochastic decay estimates of $S(t)ae$ to control $\phi_T-\phi$, which matches \eqref{eqn:pp1assump} in spatial dimension $3$ as well. This is the main intuition behind the proof of Proposition \ref{prop:phismallscale}. We use the same strategy to prove that $\psi_T$ is a good approximation of $\psi$, where we again divide the representation of $\psi_T-\psi$ into several regimes depending on the relationship between the $t_i$'s.
\smallskip

Another motivation of introducing the semigroup framework comes from deriving stochastic estimates for $\psi,\Psi$ and $\rss$ in our finite range setting. For ensembles that satisfy a logarithmic Sobolev inequality, these estimates are already established in \cite{bella2017stochastic}, while similar estimates are absent in our setting of finite range dependence (as a side note, contrary to intuition, ensembles that satisfy finite range dependence may not satisfy an LSI, see \cite[example after Theorem 6]{bella2016corrector}). The semigroup framework will be our main tool for establishing stochastic bounds on $\psi,\Psi$ (see Proposition \ref{prop:psibounds} for its precise statement), which eventually leads to the stochastic estimates on $\rss$, c.f. \eqref{eqn:thm1rss}.
\smallskip

We would also comment here that the representation of $(\psi,\Psi)$ involves two instances of the random $S(t)$, losing stochastic integrability twice, which is the reason why we only obtain the stochastic integrability of $\frac{2s}{s+2}\big|_{s=2-}=1-$. For the same reason, we would expect the $k$-th order corrector (when dimension $d\ge 2k-1$) to have stochastic integrability $\frac{2}{k}-$.

\subsection{Other Technical Lemmas}
We continue with presenting auxiliary lemmas for the proof of Proposition \ref{prop:applocqt}.  In order to capitalize in stochastic cancellations as in the proof of \cite[Theorem 1]{GO2015most}, for which we divide into dyadic scales and apply the CLT-estimate Lemma \ref{lem:GOlm13} to $\bigl(\opS(T)-\opS^h(\frac{T}{2})\opS(\frac{T}{2})\bigr)g$, we need the following approximate locality result. 
\begin{lemma}\label{lem:buckling}
	For $\sqrt{T}\ge r_0$, $\bigl(\opS(T)-\opS^h(\frac{T}{2})\opS(\frac{T}{2})\bigr)g$ is approximately local on scale $\sqrt{T}$ relative to $\bar{F}+(\frac{r_0}{\sqrt{T}})^\frac{d}{2}\bar{g}$, where \begin{equation}\label{eqn:barF}
	\bar{F}:=\fint_0^T \ud t(\dfrac{\sqrt{t}}{\sqrt{T}})^\frac{d}{2}\fint_0^{\sqrt{t}}\ud r(\dfrac{r}{\sqrt{t}})^\frac{d}{2}|(\opS(t)g-\langle \opS(t)g\rangle)_r|.
	\end{equation}
\end{lemma}

We also need some auxiliary estimates of $(\phi,\sigma,q)$ to prove Proposition \ref{prop:applocqt}, which are listed in the following Lemma \ref{lem:phisigb1}. The first two results are proven in \cite{GO2015most} while the latter three are not explicitly stated in \cite{GO2015most} since they involve stationary $(\phi,\sigma,q)$, so we provide a proof for them. \begin{lemma}\label{lem:phisigb1} Let $d\ge 3$, $\rr$ be the random radius defined in Proposition \ref{prop:phismallscale}, then for any $R\ge 1$,  
\begin{flalign}  \text{\cite[Corollary 1]{GO2015most}}  \hspace{6em} & 
     \lVert\bigl(\nabla \phi, \nabla \sigma, q-\langle q \rangle , \opS(t)ae-\langle \opS(t)ae \rangle\bigr)_R \rVert_{2-}\lesssim R^{-\frac{d}{2}},  \label{eqn:psqeq1} &\\   
     \text{\cite[Corollary 4]{GO2015most}} \hspace{6em} & \lVert (\int \eta_R |\nabla\phi|^2)^\frac{1}{2}\rVert_{2-}  \lesssim 1, \label{eqn:psqeq2} \\  & \lVert (\phi,\sigma)_R \rVert_{2-} \lesssim R^{1-\frac{d}{2}}, \label{eqn:lm4eq2}
    \\ & \bigl \lVert I(R\ge \rr) (a\phi-\langle a\phi \rangle)_R \bigr\rVert_{2-} \lesssim R^{1-\frac{d}{2}}. \label{eqn:aphiR}\\  & \bigl\lVert \big(\int \eta_R| (\phi,\sigma)|^2\bigr)^\frac{1}{2}\bigr\rVert_{2-} \lesssim 1, \label{eqn:lm4eq1}  
    \end{flalign}
\end{lemma}
\medskip




We conclude with two technical lemmas which facilitate the proof of Proposition \ref{prop:psibounds}. The first lemma acts like a combination of Lemma \ref{lem:GOlm1} (with $g$ replaced by $aS(t_1)g$) and Corollary \ref{cor:upgrademvp}.
\begin{lemma}\label{lem:cheating}
    For any $ R>0$, $t_1\ge t_0>0$,
\begin{align*} & \Bigl\|I(R \ge \rr) \Bigl(\int \eta_{R} \bigl|\int_0^{t_0} \ud \tau \nabla S(\tau)  aS(t_1)g\bigr|^2\Bigr)^\frac{1}{2} \Bigr\|_{2-} \\ & \qquad \lesssim \Bigl\|I(R \vee \sqrt{t_1}\ge \rr)\Bigl(\int \eta_{R \vee \sqrt{t_1}} \bigl(|\sqrt{t_1}\nabla S(\frac{t_1}{2})g|^2 + |S(\frac{t_1}{2})g|^2\bigr)\Bigr)^\frac{1}{2}\Bigr\|_{2-}. \stepcounter{equation} \tag{\theequation} \label{eqn:nobarcheating} \end{align*}
\end{lemma}
\noindent The second lemma deals with estimating the triple integral term of $\psi_T$, see \eqref{eqn:psiTdecom}, in the regime $t_2<t_1$. The goal is to utilize the bounds on $\opS(t)ae$, i.e. Lemma \ref{lem:phisigb1}. \begin{lemma}\label{lem:smalltibp}For $T \ge 1$, $\tdet \ge t_1>0$, \begin{align*}
         \bigl\lVert  \int_0^{t_1}\ud t_2 \bigl(1-\exp(-\dfrac{\tdet +t_2}{T})\bigr)\bigl(\sqrt{t_1}\nabla \bar{S}(t_1)\times a\nabla S(t_2)ae, \, \bar{S}(t_1)\times a \nabla & S(t_2)ae\bigr) \bigr\rVert_{2-}  \\ & \lesssim (1 \wedge \dfrac{\tdet}{T}) \frac{1}{\sqrt{t_1}} (1 \wedge \frac{1}{\sqrt{t_1}})^\frac{d}{2}, \stepcounter{equation} \tag{\theequation}\label{eqn:smalltibp} \\  \bigl\lVert  \int_0^{t_1}\ud t_2 \exp(-\dfrac{t_2}{T})(\sqrt{t_1}\nabla \bar{S}(t_1)\times a\nabla S(t_2)ae,\bar{S}(t_1)\times a\nabla S(t_2)ae )\bigr\rVert_{2-} &\lesssim  \frac{1}{\sqrt{t_1}} (1 \wedge \frac{1}{\sqrt{t_1}})^\frac{d}{2}. \stepcounter{equation}\tag{\theequation}\label{eqn:smalltibp2} 
    \end{align*} 
  Moreover, $\int_0^{t_1}\ud t_2 (1-\exp(-\frac{\tdet +t_2}{T}))\bar{S}(t_1)\times a\nabla S(t_2)ae$ is approximately local on scale $1 \vee \sqrt{t_1}$ relative to some stationary $\bar{g}$ with \begin{equation}\label{eqn:smalltibp3}
        \|\bar{g}\|_{2-}\lesssim (1 \wedge \dfrac{\tdet}{T}) \frac{1}{\sqrt{t_1}} (1 \wedge \frac{1}{\sqrt{t_1}})^\frac{d}{2}.
    \end{equation} 
\end{lemma}
\FloatBarrier

\section{Proofs}\label{sec:proofs}

\subsection{Proof of Theorem \ref{thm:mainthm}} 

\noindent \emph{Step 1:} Stochastic estimates of $\rss$. 
The idea of the proof is based on \cite[Theorem 1 (ii)]{fischer2017sublinear}. Our goal is to estimate the ``failure probability'' $\langle I(\rss \ge R_0) \rangle$ for an arbitrary $R_0$. We separate the event $\rss \ge R_0$ into three possible scenarios. The first scenario to leave out is $\rr \ge \sqrt{R_0}$ which has probability at most $\exp(-\frac{1}{C}R_0^{\frac{3}{2}})$, as $\rr$ has stochastic integrability $d$ \eqref{eqn:estrstar}. For the rest of the proof we assume $\rr\le \sqrt{R_0}\le R_0$.

Next we look at the failure probability due to first-order correctors $(\phi,\sigma)$. We recall \eqref{intrphigr} which defines the constraint on $\rss$ coming from $(\phi,\sigma)$. By \eqref{eqn:lm4eq1} in Lemma \ref{lem:phisigb1}, which we reformulate in terms of boxes instead of balls, we have for all $R\ge 1$,
\[\lVert (\fint_{Q_R} |(\phi,\sigma)|^2)^\frac{1}{2}\rVert_{2-} \lesssim 1,\] so that by Chebyshev's inequality, we have\footnote{Here and for the rest of the paper $X \lesssim \exp(-R^{s-})$ denotes for any $s_0\in (0,s)$, $X \lesssim_{s_0} \exp(-R^{s_0})$.}
\begin{align*}
    \Big\langle I\Bigl((\fint_{Q_R} |(\phi,\sigma)|^2)^\frac{1}{2}\ge R^{1-\alpha}R_0^\alpha \Bigr) \Big\rangle & = \Big\langle I\Bigl( \exp\bigl((\fint_{Q_R} |(\phi,\sigma)|^2)^\frac{2-}{2}\bigr)\ge \exp\bigl((R^{1-\alpha}R_0^\alpha)^{2-}\bigr) \Bigr) \Big\rangle \\ & \lesssim\exp\Bigl(-\frac{1}{C}(R^{1-\alpha}R_0^\alpha)^{2-}\Bigr). 
\end{align*}Here the second line uses the third expression of the norm $\lVert \cdot \rVert_s$ in \eqref{eqn:equistocnorm}. Therefore if $\rss\ge R_0$ because of the ``failure" of first-order correctors, which means that there must be some $R\ge R_0$ (which we may assume to be dyadic) such that \begin{equation*}
    (\fint_{Q_R} |(\phi,\sigma)|^2)^\frac{1}{2}\ge R^{1-\alpha}R_0^\alpha,
\end{equation*} the probability is dominated by\footnote{here $\sim$ means both $\lesssim$ and $\gtrsim$ hold} \[\sum_{\substack {R\ge R_0 \\ \text{dyadic}}}\exp\Bigl(-\frac{1}{C}(R^{1-\alpha}R_0^\alpha)^{2-}\Bigr)\sim \exp(-\frac{1}{C}R_0^{2-}), \]where we used $\alpha<1$.

We now look at the failure probability due to second-order correctors. We first argue that it is enough to prove for any $R\ge 1$, \begin{equation}\label{eqn:thm1fr}
    \Bigl\lVert I(\sqrt{R}\ge \rr)\Bigl(\fint_{B_R}\bigl\lvert(\psi,\Psi)-\fint_{B_R}(\psi,\Psi)\bigr\rvert^2\Bigr)^\frac{1}{2}\Bigr\rVert_{1-} \lesssim R^\frac{1}{2}.
\end{equation}
Indeed, if \eqref{eqn:thm1fr} holds true, then we can bound the failure probability using Chebyshev inequality: \begin{align*}
    \sum_{\substack {R\ge R_0 \\ \text{dyadic}}} &\Bigl\langle  I\Bigl(\bigl(\fint_{B_R}\bigl\lvert(\psi,\Psi)-\fint_{B_R}(\psi,\Psi)\bigr\rvert^2\bigr)^\frac{1}{2}\ge R_0^\beta R^{2-\beta}   \Bigr) I(\sqrt{R_0}\ge \rr)\Bigr\rangle  \\ & \lesssim \sum_{\substack {R\ge R_0 \\ \text{dyadic}}} \exp\Bigl(-\bigl(\dfrac{R_0^{\beta } R^{2-\beta}}{cR^\frac{1}{2}}\bigr)^{1-}\Bigr) \sim \exp(-\frac{1}{C}R_0^{\frac{3}{2}-}),
\end{align*}
where we used $\beta<\frac{3}{2}$. To prove \eqref{eqn:thm1fr}, we again divide into dyadic series and use \eqref{eqn:pp6eq1} and \eqref{eqn:pp6eq3} in Proposition \ref{prop:psibounds}. We will abuse notation and use $f_r$ to denote $\fint_{Q_r} f$, as \cite[Lemma 13, Step 4]{GO2015most} shows its equivalence to $f*G_r$ under stochastic norm $\lVert \cdot \rVert$. \begin{align*}
    \Bigl\lVert  I&( \sqrt{R}\ge \rr)\Bigl(\fint_{B_R}\bigl\lvert(\psi,\Psi)-\fint_{B_R}(\psi,\Psi)\bigr\rvert^2\Bigr)^\frac{1}{2}\Bigr\rVert_{1-}   \\ & \le \biggl\lVert I(\sqrt{R}\ge \rr)\biggl(\Bigl(\fint_{B_R}\bigl\lvert(\psi,\Psi)-(\psi,\Psi)_{\rr} \bigr\rvert^2\Bigr)^\frac{1}{2}  +  \sum_{\substack {R/2 \ge r\ge \rr \\ \text{dyadic}}} \Bigl(\fint_{B_R}\bigl\lvert(\psi,\Psi)_r-(\psi,\Psi)_{2r}\bigr\rvert^2\Bigr)^\frac{1}{2}\biggr)\biggr\rVert_{1-} \\ & \le  \Bigl\lVert I(\sqrt{R}\ge \rr)\Bigl(\fint_{B_R}\bigl\lvert(\psi,\Psi)-(\psi,\Psi)_{\rr} \bigr\rvert^2\Bigr)^\frac{1}{2}\Bigr\rVert_{1-} \\ & \qquad + \sum_{\substack {R/2 \ge r\ge 1 \\ \text{dyadic}}} \Bigl\lVert I(r\ge \rr)\Bigl(\fint_{B_R}\bigl\lvert(\psi,\Psi)_r-(\psi,\Psi)_{2r}\bigr\rvert^2\Bigr)^\frac{1}{2}\Bigr\rVert_{1-}\\ & \lesssim \bigl\lVert I(\sqrt{R}\ge \rr)\rr \bigl(\int \eta_{R}\lvert \nabla(\psi,\Psi) \bigr\rvert^2\bigr)^\frac{1}{2}\bigr\rVert_{1-}+ \sum_{\substack {R/2 \ge r\ge 1 \\ \text{dyadic}}} \Bigl\lVert I(r\ge \rr) r \bigl(\int \eta_{R}\bigl\lvert \nabla (\psi,\Psi)_r \bigr\rvert^2\Bigr)^\frac{1}{2}\Bigr\rVert_{1-} \\ & \le \sqrt{R}\bigl\lVert I(\sqrt{R}\ge \rr) \bigl(\int \eta_{R}\lvert \nabla(\psi,\Psi) \bigr\rvert^2\bigr)^\frac{1}{2}\bigr\rVert_{1-}+ \sum_{\substack {R/2 \ge r\ge 1 \\ \text{dyadic}}} r\Bigl\lVert I(r\ge \rr)  \nabla (\psi,\Psi)_r \Bigr\rVert_{1-} \\ & \leftstackrel{\eqref{eqn:pp6eq1},\eqref{eqn:pp6eq3}}{\lesssim}\sqrt{R} + \sum_{\substack {R/2 \ge r\ge 1 \\ \text{dyadic}}} \sqrt{r} \lesssim \sqrt{R}.
\end{align*} 
Here in the third inequality we used the Poincar\'e inequality in convolution (see \cite[(201)]{GO2015most}): for any $R\gtrsim r_0$,
\begin{equation}\label{eqn:convpoincare}
    \int \eta_R (f-f_r)^2 \lesssim r^2 \int \eta_R|\nabla f|^2.
\end{equation}This establishes $\langle I(\rss\ge R_0)\rangle \lesssim \exp(-\frac{1}{C}R_0^{\frac{3}{2}-})$ for any $R_0\ge 1$, which is equivalent to \eqref{eqn:thm1rss}.

\medskip 

\noindent \emph{Step 2:} Estimation of failure probabilities. The plan is to pass from the deterministic estimate Proposition \ref{prop:apriori} to the probabilistic statement Theorem \ref{thm:mainthm}, and estimate the probability for the assumptions \eqref{eqn:pp1assump}-\eqref{eqn:pp1approxah} in the Proposition to hold. We would like to comment that the rest of the assumptions in Proposition \ref{prop:apriori} are standard, and, thanks to our finite range assumption, hold with probability 1 (see Subsection \ref{subsec:multipole} for the standard properties of correctors and $\rss$, and \cite[Corollary 1]{bella2017effective} for an argument of \eqref{eqn:qualtlgscconv}).

\smallskip

The starting point is the stochastic bounds that are established on stationary approximations of these quantities (\textit{i.e.}, the quantities without $L$ in superscripts), namely Proposition \ref{prop:phismallscale}, Lemma \ref{lem:ahahL}, and Proposition \ref{prop:psibounds}. Therefore, to estimate the ``failure probability'' of $|a_he_i-\int \omega q_{i,T}|$, we use a Chebyshev inequality as well as \eqref{eqn:ppahahL} (notice that $\rss \ge 1$): \begin{equation}\label{eqn:thmchebyshev}\begin{aligned}
   \Bigl\langle I\bigl(|a_he_i-\int \omega q_{i,T}| \ge (\frac{\rss}{L})^{\beta}\bigr) \Bigr\rangle & \le \langle I(|a_he_i-\int \omega q_{i,T}| \ge L^{-\beta}) \rangle \\ & \lesssim \exp\bigl(-\frac{1}{C}(L^{\beta}T^{-\frac{3}{4}})^{2-}\bigr)=\exp(-\frac{1}{C}L^{\varepsilon-}).
\end{aligned}\end{equation}   
We may also replace $q_{i,T}$ with $q_{i,T}^{(L)}$ in the above \eqref{eqn:thmchebyshev}, since by the deterministic Proposition \ref{prop:truncation}, we can estimate \begin{align*}
    |\int \omega (q_{i,T}-q_{i,T}^{(L)})| \hspace{1em} \leftstackrel{\eqref{eqn:defqiTL},\eqref{intrhomcoeff}}{=}|\int \omega a\nabla (\phi_{i,T}-\phi_{i,T}^{(L)})|  \lesssim (\int \omega^2)^\frac{1}{2}\bigl(\int_{Q_L} |\nabla (\phi_{i,T}-\phi_{i,T}^{(L)})|^2\bigr)^\frac{1}{2} \ \leftstackrel{\eqref{eqn:phitruncation}}{\lesssim} (\frac{\sqrt{T}}{L})^p,
\end{align*} which is much smaller than $L^{-\beta}$ for $\sqrt{T}=L^{1-\varepsilon}$ and $p$ sufficiently large. This shows \begin{equation*}
    \Bigl\langle I\bigl(|a_he_i-\int \omega q_{i,T}^{(L)}| \ge (\frac{\rss}{L})^{\beta}\bigr) \Bigr\rangle \lesssim\exp(-\frac{1}{C}L^{\varepsilon-}).
\end{equation*}

The failure probability for the terms in \eqref{eqn:pp1assump} will be estimated slightly differently as they involve $\rr$. We take $\bigl(\fint_{B_\ell} |\nabla (\psi_T^{(L)}-\psi)|^2\bigr)^\frac{1}{2}$ as an example since it is among the terms that have the worst stochastic integrability. Again by Chebyshev inequality, we obtain \begin{equation*}
    \Bigl\langle I\Bigl(\bigl(\fint_{B_\ell} |\nabla(\psi_{T}-\psi)|^2\bigr)^\frac{1}{2} \ge \sqrt{T}(\dfrac{1}{L})^\beta  \Bigr)I(\ell \ge \rr)\Bigr\rangle \leftstackrel{\eqref{eqn:pp6eq4}}{\lesssim} \exp\Bigl(-\frac{1}{C}(L^{\beta }T^{-\frac{3}{4}})^{1-}\Bigr)=\exp(-\frac{1}{C}L^{\frac{\varepsilon}{2}-}),
\end{equation*}which, in view of \eqref{eqn:psiapprox} and $\rss \ge \rr$, can be changed to \begin{equation*}
    \Bigl\langle I\Bigl(\bigl(\fint_{B_\ell} |\nabla(\psi_{T}^{(L)}-\psi)|^2\bigr)^\frac{1}{2} \ge \sqrt{T}(\dfrac{1}{L})^\beta  \Bigr)I(\ell \ge \rr)\Bigr\rangle \lesssim \exp(-\frac{1}{C}L^{\frac{\varepsilon}{2}-}).
\end{equation*}We now remove the constraint $\ell \ge \rr$ using Bayes' formula 
\begin{align*}
    \Bigl\langle I \Bigl(&\bigl(\fint_{B_\ell} |\nabla(\psi_{T}^{(L)}-\psi)|^2\bigr)^\frac{1}{2} \ge \sqrt{T}(\dfrac{\rss}{L})^\beta  \Bigr)\Big\lvert \ell \ge \rss \Bigr\rangle \\ & \le \Bigl\langle I \Bigl(\bigl(\fint_{B_\ell} |\nabla(\psi_{T}^{(L)}-\psi)|^2\bigr)^\frac{1}{2} \ge \sqrt{T}(\dfrac{1}{L})^\beta  \Bigr)\Big\lvert \ell \ge \rss \Bigr\rangle \\ & =  \dfrac{\Bigl\langle I\Bigl(\bigl(\fint_{B_\ell} |\nabla(\psi_{T}^{(L)}-\psi)|^2\bigr)^\frac{1}{2} \ge \sqrt{T}(\frac{1}{L})^\beta  \Bigr)I(\ell \ge \rss)\Bigr\rangle}{\langle I(\ell\ge \rss) \rangle} \\ & \le   \dfrac{\Bigl\langle I\Bigl(\bigl(\fint_{B_\ell} |\nabla(\psi_{T}^{(L)}-\psi)|^2\bigr)^\frac{1}{2} \ge \sqrt{T}(\frac{1}{L})^\beta  \Bigr)I(\ell \ge \rr)\Bigr\rangle}{\langle I(\ell\ge \rss) \rangle} \\ & \leftstackrel{\eqref{eqn:thm1rss}}{\lesssim}
    \dfrac{\exp(-\frac{1}{C}L^{\frac{\varepsilon}{2}-})}{1-\exp(-\frac{1}{C}\ell^{\frac{3}{2}-})}.
\end{align*} 
We take $\ell$ to be reasonably large and  independent of $L,T,\beta$ so that $1-\exp(-\frac{1}{C}\ell^{\frac{3}{2}-}) \ge \frac{1}{2}$. Since the probability of all failure events can be estimated as such, we derive the total failure probability estimate as desired. 

Finally for the failure probability for \eqref{eqn:pp1psiTgrowth}, we first obtain by Chebyshev inequality \begin{equation*}
    \Bigl\langle I\Bigl(\bigl(\fint_{B_L} \psi_{T}^2\bigr)^\frac{1}{2} \ge T(\dfrac{1}{L})^\beta  \Bigr)I(L \ge \rr)\Bigr\rangle \leftstackrel{\eqref{eqn:pp6eq2}}{\lesssim} \exp(-\frac{1}{C}L^{\frac{\varepsilon}{2}-}),\end{equation*} which, in view of \eqref{eqn:psinoshiftgrowth} and \eqref{eqn:psiapprox}, as well as $\rss \ge \rr$, can be changed to \begin{equation*}
    \Bigl\langle I\Bigl(\bigl(\fint_{B_L} |\psi_{T}^{(L)}-\psi)|^2\bigr)^\frac{1}{2} \ge L^2 (\frac{\rss}{L})^\beta \Bigr)I(L \ge \rss)\Bigr\rangle \lesssim \exp(-\frac{1}{C}L^{\frac{\varepsilon}{2}-}),
\end{equation*}and we again finish the proof by dropping the constraint $L\ge \rss$ using Bayes' formula.
\qed

\subsection{Proof of Proposition \ref{prop:effectivemultipole}}
By \cite[Theorem 2]{bella2017effective}, for $u_h$ that satisfies \eqref{eqn:lm1homphsm}, we have
\begin{equation*}
    \sup_{R\ge \rss} (\dfrac{R}{\rss})^{d+\beta} \Bigl(\frac{1}{R^d}\int_{B_R^c} |\nabla \bigl(u-(1+\phi_i\partial_i+\psi_{ij}\partial_{ij})u_h\bigr)|^2\Bigr)^\frac{1}{2} \lesssim \sup_{R\ge \rss} (\dfrac{R}{\rss})^d (\frac{1}{R^d}\int_{B_R^c} |\nabla u|^2)^\frac{1}{2}.
\end{equation*} Thus, it suffices to prove that, for every $R\ge \rss$, \begin{equation*}
    R^\frac{d}{2} (\int_{B_R^c} |\nabla u|^2)^\frac{1}{2} \lesssim \ell^d.
\end{equation*}
For $\rss \le R\le 2\ell$, the proof is a standard energy estimate: \begin{align*}
    & R^\frac{d}{2} (\int_{B_R^c} |\nabla u|^2)^\frac{1}{2}  \le R^\frac{d}{2}(\int |\nabla u|^2)^\frac{1}{2}   \leftstackrel{\eqref{eqn:intrbaseq}}{\lesssim} R^\frac{d}{2}(\int_{B_{\ell}} |g|^2)^\frac{1}{2}  \le R^\frac{d}{2}\ell^\frac{d}{2} \lesssim \ell^d.
\end{align*}
The proof for $R\ge 2\ell$ uses \cite[Lemma 4 (c)]{bella2017effective}: since $u\in Y_1(l)$, \begin{equation*}
    (\frac{1}{R^d}\int_{B_R^c} |\nabla u|^2)^\frac{1}{2}  \lesssim (\dfrac{\ell}{R})^d (\frac{1}{\ell^d}\int_{B_\ell^c} |\nabla u|^2)^\frac{1}{2} \lesssim  (\dfrac{\ell}{R})^d\ell^{-\frac{d}{2}}(\int |\nabla u|^2)^\frac{1}{2} \leftstackrel{\eqref{eqn:intrbaseq}}{\lesssim} \dfrac{\ell^\frac{d}{2}}{R^d}(\int |g|^2)^\frac{1}{2} \lesssim (\dfrac{\ell}{R})^d.  
\end{equation*}
\qed

\subsection{Proof of Corollary \ref{cor:multipole}}
Define $w=\hat{u}-u$ and $w_D=(1+\phi_i\partial_i+\psi_{ij}\partial_{ij})u_h-u$, then \begin{equation}\label{co1form}
-\nabla\cdot a\nabla w=0\hspace{0.1in}\text{in }Q_L,\hspace{0.5in} w=w_D\hspace{0.1in}\text{ on }\partial Q_L.
\end{equation}
Since $w$ is $a$-harmonic, by \eqref{eqn:ellipticmvp} (note that $\rss\ge \rr$),\begin{equation}\label{co1minrad}
\fint_{Q_R}|\nabla w|^2 \lesssim \fint_{Q_L}|\nabla w|^2\hspace{0.1in}\text{ for }\rss\le  R \le L.
\end{equation}
Thus, it suffices to prove \begin{equation}\label{co1estw}
(\int_{Q_L} |\nabla w|^2)^{\frac{1}{2}} \lesssim (\int_{Q_{2L}-Q_L} |\nabla w_D|^2)^{\frac{1}{2}},
\end{equation}
as \eqref{co1minrad} and \eqref{co1estw} together with Proposition \ref{prop:effectivemultipole} yield Corollary \ref{cor:multipole}.

By rescaling, we may without loss of generality assume $L=1$, and we can further assume $\int_{Q_2-Q_1} w_D=0$ as the expressions in both sides of \eqref{co1estw} are invariant through the subtraction of constants. This allows us to define an extension $\bar{w}_D$ of $w_D$ on $Q_2$ such that $\bar{w}_D=w_D$ on $Q_1$ and (using Poincar\'{e} inequality)\begin{equation}\label{co1bar}
(\int_{Q_2} |\nabla \bar{w}_D|^2)^{\frac{1}{2}}\lesssim (\int_{Q_2-Q_1} |\nabla w_D|^2)^{\frac{1}{2}}.
\end{equation}
Hence \eqref{co1form} can be reformulated as \begin{equation*}
-\nabla\cdot a\nabla (w-\bar{w}_D)=\nabla \cdot a\nabla\bar{w}_D \text{ in }Q_1,\ \  w-\bar{w}_D=0\hspace{0.1in}\text{ on }\partial Q_1.
\end{equation*}
Now \eqref{co1estw} follows from the standard energy estimate $(\int_{Q_1}|\nabla(w-\bar{w}_D)|^2)^{\frac{1}{2}}\lesssim (\int_{Q_1} |\nabla \bar{w}_D|^2)^{\frac{1}{2}}$, triangle inequality and \eqref{co1bar}.
\qed

\subsection{Proof of Proposition \ref{prop:apriori}}
The proof largely resembles that of \cite[Proposition 1]{lu2018optimal}. We divide the proof into four steps.\\
\noindent \emph{Step 1:} We upgrade \eqref{intrpsiorig} and show a seemingly stronger condition \begin{equation}\label{eqn:upgradegrowth}
\dfrac{1}{r^2}(\fint_{B_r}|(\psi,\Psi)-\fint_{B_{\rss}} (\psi,\Psi)|^2)^{\frac{1}{2}} \le (\frac{\rss}{r})^\beta\hspace{0.1in}\text{ for all }r\ge \rss. 
\end{equation}
What separates \eqref{eqn:upgradegrowth} from \eqref{intrpsiorig} is \begin{equation*}
    |\fint_{B_r} (\psi,\Psi)-\fint_{B_{\rss}} (\psi,\Psi)|\lesssim r^2(\dfrac{\rss}{r})^\beta \ \mbox{ for } \ r\ge \rss.
\end{equation*}
To prove this, we apply dyadic decomposition, and (since $\beta<2$) reduce this to \begin{equation*}
    |\fint_{B_r} (\psi,\Psi)-\fint_{B_{r'}} (\psi,\Psi)|\lesssim r^2(\dfrac{\rss}{r})^\beta \ \mbox{ for } \ 2r'\ge r \ge r' \ge \rss,
\end{equation*}
which, by triangle inequality, is a consequence of \eqref{intrpsiorig}. Hence, as $(\psi,\Psi)$ are defined up to a constant, we may assume
\begin{equation}\label{eqn:psinoshiftgrowth}
    \dfrac{1}{r^2}(\fint_{B_r}|(\psi,\Psi)|^2)^{\frac{1}{2}} \lesssim (\frac{\rss}{r})^\beta\hspace{0.1in}\text{ for all }r\ge \rss.
\end{equation}
\smallskip

\noindent For the rest of this proof, by a scaling argument we may without loss of generality assume $\ell=1$.

\noindent \emph{Step 2:} We compare $\tilde{u}_h^{(L)}$, defined in \eqref{eqn:alguhtildeL} which satisfies $-\nabla\cdot a_h^{(L)} \nabla \tilde{u}_h^{(L)} =\nabla \cdot g$, with $\tilde{u}_h$, the solution of \eqref{intruhtilde}, and claim that \begin{equation}\label{eqn:pp1tildeu}
|\nabla(\tilde{u}_h^{(L)}-\tilde{u}_h)|+L|\nabla^2(\tilde{u}_h^{(L)}-\tilde{u}_h)|+L^2|\nabla^3(\tilde{u}_h^{(L)}-\tilde{u}_h)|\lesssim (\dfrac{1}{L})^d(\dfrac{\rss}{L})^\beta\hspace{0.1in}\text{ on }Q_L^c.
\end{equation}
To prove \eqref{eqn:pp1tildeu}, we use that the support of $g$ is contained in $B_1$. It is well known that $|\nabla^n G_h(x)| \lesssim |x|^{2-d-n}$ for any multi-index $n$. Using the representation formula, \begin{equation*}
    \tilde{u}_h = G_h * \nabla\cdot g =  \nabla G_h * g
\end{equation*} for $|x| \gg 1$, and standard Schauder theory for $|x| \lesssim 1$ (here we use the regularity of $g\in C^{2,\gamma}$), we obtain for all $x$, \begin{equation}\label{eqn:pp1tildeunab} \begin{aligned}
   & |\tilde{u}_h(x)|\lesssim \dfrac{1}{(1+|x|)^{d-1}}, & |\nabla \tilde{u}_h(x)|\lesssim \dfrac{1}{(1+|x|)^d}, \\ & |\nabla^2\tilde{u}_h(x)|\lesssim \dfrac{1}{(1+|x|)^{d+1}}, & |\nabla^3 \tilde{u}_h(x)|\lesssim \dfrac{1}{(1+|x|)^{d+2}},
\end{aligned}\end{equation} and \begin{equation}\label{eqn:improveschauder}
     \sup_{0<|y|\le 1}\frac{|\nabla^3 \tilde{u}_h (x+y)-\nabla^3\tilde{u}_h(x)|}{|y|^\gamma} \lesssim \frac{1}{(1+|x|)^{d+2}}.
\end{equation}
Next we consider $w=\tilde{u}_h^{(L)}-\tilde{u}_h$, which satisfies the equation \begin{equation*}
    -\nabla\cdot a_h^{(L)}\nabla w = \nabla \cdot (a_h^{(L)}-a_h)\nabla\tilde{u}_h.
\end{equation*} Similarly we have the representation \begin{equation*}
    w = G_h^{(L)} *\nabla \cdot (a_h^{(L)}-a_h)\nabla \tilde{u}_h = \nabla  G_h^{(L)} *  (a_h^{(L)}-a_h)\nabla\tilde{u}_h.
\end{equation*} Therefore, using  \eqref{eqn:pp1tildeunab} and \eqref{eqn:improveschauder}, again we appeal to representation formula and standard Schauder theory and derive \begin{equation*}\begin{aligned}
    |\nabla w(x)| \lesssim |a_h^{(L)}-a_h| \dfrac{1}{(1+|x|)^d} \leftstackrel{\eqref{eqn:pp1approxah}}{\lesssim} (\dfrac{\rss}{L})^\beta \dfrac{1}{(1+|x|)^d}, \\ |\nabla^2 w(x)| \lesssim |a_h^{(L)}-a_h| \dfrac{1}{(1+|x|)^{d+1}} \leftstackrel{\eqref{eqn:pp1approxah}}{\lesssim} (\dfrac{\rss}{L})^\beta \dfrac{1}{(1+|x|)^{d+1}},\\ |\nabla^3 w(x)| \lesssim |a_h^{(L)}-a_h| \dfrac{1}{(1+|x|)^{d+2}} \leftstackrel{\eqref{eqn:pp1approxah}}{\lesssim} (\dfrac{\rss}{L})^\beta \dfrac{1}{(1+|x|)^{d+2}}, \end{aligned}
\end{equation*}
which is exactly \eqref{eqn:pp1tildeu}.
\smallskip

\noindent \emph{Step 3:} We compare $u_h^{(L)}$ and $u_h$ defined in \eqref{eqn:algapproxbdry} and \eqref{eqn:effectivequadp} and claim \begin{equation}\label{pp1uh}
|\nabla(u_h^{(L)}-u_h)|+L|\nabla^2(u_h^{(L)}-u_h)|+L^2|\nabla^3(u_h^{(L)}-u_h)|\lesssim (\dfrac{1}{L})^d(\dfrac{\rss}{L})^\beta \ \text{ on }Q_L^c
\end{equation}
and \begin{equation}\label{pp1uhmag}
|\nabla u_h^{(L)}|+L|\nabla^2 u_h^{(L)}|+L^2|\nabla^3 u_h^{(L)}|\lesssim (\dfrac{1}{L})^d \ \text{ on }Q_L^c.
\end{equation}
To this purpose, we define 
\begin{equation}\label{eqn:xixil}
    \xi_i:=\int \nabla  g \cdot\phi_i \ \mbox{ and } \ \xi_{i,T}^{(L)}:=\int  g \cdot\nabla \phi_{i,T}^{(L)}.
\end{equation}
In view of \eqref{eqn:pp1tildeu}, it suffices to control the dipole and quadrupole terms: \begin{equation}\begin{aligned}\label{pp1dpdiff}
|\nabla(\xi_{i,T}^{(L)}\partial_iG_h^{(L)}-\xi_i\partial_iG_h)|+L|\nabla^2(\xi_{i,T}^{(L)}\partial_iG_h^{(L)}-\xi_i\partial_iG_h)|\\+L^2|\nabla^3(\xi_{i,T}^{(L)}\partial_iG_h^{(L)}-\xi_i\partial_iG_h)| \lesssim (\dfrac{1}{L})^d(\dfrac{\rss}{L})^\beta
\end{aligned} \end{equation} and
\begin{equation}
\label{pp1qpdiff}\begin{aligned} |\nabla(c_{ij,T}^{(L)}\partial_{ij}G_h^{(L)}-c_{ij}\partial_{ij}G_h)|+L|\nabla^2(c_{ij,T}^{(L)}\partial_{ij}G_h^{(L)}-c_{ij}\partial_{ij}G_h)|\\ +L^2|\nabla^3(c_{ij,T}^{(L)}\partial_{ij}G_h^{(L)}-c_{ij}\partial_{ij}G_h)|\lesssim (\dfrac{1}{L})^d(\dfrac{\rss}{L})^\beta, \ \mbox{ both on }Q_L^c. \end{aligned}
\end{equation} 
We have some obvious estimates for constant coefficient Green's functions:
\begin{equation*}|\nabla^2 G_h^{(L)}|+L|\nabla^3 G_h^{(L)}|+L^2|\nabla^4 G_h^{(L)}|+L^3|\nabla^5 G_h^{(L)}| \lesssim (\dfrac{1}{L})^d\hspace{0.1in}\text{ on }Q_L^c,\end{equation*} and \begin{equation*}\begin{aligned}|\nabla^2 (G_h^{(L)}-G_h)|+L|\nabla^3 (G_h^{(L)}-G_h)|+L^2|\nabla^4 (G_h^{(L)}-G_h)| +L^3|\nabla^5 (G_h^{(L)}-G_h)| \\ \lesssim |a_h^{(L)}-a_h|  (\dfrac{1}{L})^d \leftstackrel{\eqref{eqn:pp1approxah}}{\lesssim}(\dfrac{1}{L})^d(\dfrac{\rss}{L})^\beta  \ \text{ on }Q_L^c.\end{aligned}\end{equation*}
Hence, recalling the definitions of $c_{ij}$ and $c_{ij,T}^{(L)}$, c.f. \eqref{eqn:cij} and \eqref{eqn:cijlt}, it suffices to show \begin{equation}\label{eqn:pp1xicbd}
    |\xi| \lesssim 1 \ \mbox{ and } \ |c_{ij}|\lesssim 1
\end{equation} and \begin{equation}\label{eqn:pp1xicdiff}
    |\xi_T^{(L)}-\xi|\lesssim(\dfrac{\rss}{L})^\beta \hspace{0.1in}\text{     and     }\hspace{0.1in}\dfrac{1}{L}|c_{ij,T}^{(L)}-c_{ij}|\lesssim (\dfrac{\rss}{L})^\beta .
\end{equation} 
The arguments of \eqref{eqn:pp1xicbd} and \eqref{eqn:pp1xicdiff} for $\xi$ and $\xi_T^{(L)}$ are straightforward: \[|\xi|=|\int \phi \nabla \cdot g|\lesssim (\fint_{B_1} \phi^2)^\frac{1}{2} \leftstackrel{\eqref{intrphigr}}{\lesssim} 1,\] and \begin{equation*}
|\xi_T^{(L)}-\xi|= |\int   g \cdot\nabla(\phi_T^{(L)}-\phi)|  \lesssim (\fint_{B_1}|\nabla(\phi_{T}^{(L)}-\phi)|^2)^{\frac{1}{2}} \leftstackrel{\eqref{eqn:pp1assump}}{\lesssim} (\dfrac{\rss}{L})^\beta.
\end{equation*} 
We now prove \eqref{eqn:pp1xicbd} for $c_{ij}$. By the growth condition \eqref{intrpsiorig} of $\psi$ for $r=1$, we have \begin{equation}\label{eqn:pp1psitlpsi}
    |\int \psi \nabla \cdot g |= |\int (\psi-\fint_{B_1} \psi)\nabla \cdot g  | \lesssim (\fint_{B_1}| \psi-\fint_{B_1}\psi|^2)^{\frac{1}{2}} \lesssim 1.
\end{equation}
Moreover, since $\partial_k v_{h,ij}$, defined in \eqref{eqn:harmpol}, is a polynomial homogeneous of order $1$, we have \begin{equation}\label{eqn:pp1vhilbd}
    |\partial_k v_{h,ij}(x)| \lesssim 1 \ \text{ and } \ |\nabla \partial_k v_{h,ij}(x)| \lesssim 1 \ \text{ for } x\in B_1,
\end{equation} thus \begin{equation}\label{eqn:pp1phivh}
    |\int \phi_k\partial_k v_{h,ij} \nabla\cdot g |  \lesssim (\fint_{B_1}(\phi_k\partial_k v_{h,ij})^2)^{\frac{1}{2}}  \lesssim  (\fint_{B_1}\phi_k^2)^{\frac{1}{2}} \leftstackrel{\eqref{intrphigr}}{\lesssim} 1.
\end{equation}Combining \eqref{eqn:pp1psitlpsi} and \eqref{eqn:pp1phivh} we established \eqref{eqn:pp1xicbd} for $c_{ij}$.

The proof of \eqref{eqn:pp1xicdiff} for $|c_T^{(L)}-c|$ is more complicated. Recalling the definitions of $c_{ij}$ and $c_{ij,T}^{(L)}$, c.f. \eqref{eqn:cij} and \eqref{eqn:cijlt}, we may write $|c_T^{(L)}-c|:=|\int \Xi \cdot g|$ and decompose $\Xi$ into different terms:  \begin{align*}
  \Xi_{ij} & =(\partial_k v_{h,ij}-\partial_k v_{h,ij}^{(L)}) \nabla \phi_k  + \partial_k v_{h,ij}^{(L)}\nabla (\phi_k-\phi_{k,T}^{(L)}) + \phi_k\nabla(\partial_k v_{h,ij}-\partial_k v_{h,ij}^{(L)})    \\ & \qquad +(\phi_k-\phi_{k,T}^{(L)})\nabla  \partial_k v_{h,ij}^{(L)} + (2-\delta_{ij})\nabla (\psi_{ij}-\psi_{ij,T}^{(L)})+ (2-\delta_{ij})(\dfrac{a_{hij}}{a_{h11}}-\dfrac{a_{hij}^{(L)}}{a_{h11}^{(L)}})\nabla \psi_{11}.
\end{align*} 
Since $x_i+\phi_i$ is $a$-harmonic, we can use Caccioppoli's estimate and the growth condition \eqref{intrphigr} on $\phi$  to derive  \begin{equation}\label{eqn:pp1gradphiconst}
    (\fint_{B_R} |\nabla \phi|^2)^\frac{1}{2} \lesssim 1+ \dfrac{1}{R}(\fint_{B_{2R}} \phi^2)^\frac{1}{2} \lesssim 1+(\dfrac{\rss}{R})^\alpha \lesssim 1, \ \mbox{ for }R\ge \rss.
\end{equation}
Using the same argument as \eqref{eqn:pp1gradphiconst} for $\psi$, we obtain \begin{equation}\begin{aligned}\label{eqn:pp1gradpsi}
    (\fint_{B_R} |\nabla \psi|^2)^\frac{1}{2} & \leftstackrel{\eqref{eqn:2ndcordef}}{\lesssim} \dfrac{1}{R}(\fint_{B_{2R}} |\psi-\fint_{B_{2R}} \psi|^2)^\frac{1}{2} + (\fint_{B_{2R}} |(\phi,\sigma)|^2)^\frac{1}{2}\\ & \leftstackrel{\eqref{intrpsiorig},\eqref{intrphigr}}{\lesssim} R(\dfrac{\rss}{R})^\beta + R(\dfrac{\rss}{R})^\alpha \lesssim R \ \mbox{ for } R\ge \rss.
\end{aligned}\end{equation}  
Since $\partial_k v_{h,ij}$ and $\partial_k v_{h,ij}^{(L)}$ (defined in \eqref{eqn:harmpolL}) are two polynomials that are homogeneous of order $1$, and their coefficients differ by $|a_h^{(L)}-a_h|$, which by \eqref{eqn:pp1approxah} is bounded by $(\frac{\rss}{L})^\beta$, we have \begin{equation}\label{eqn:pp1vhildiff}
    |\partial_k v_{h,ij}-\partial_k v_{h,ij}^{(L)}| \lesssim (\dfrac{\rss}{L})^\beta \ \text{ and } \ |\nabla (\partial_k v_{h,ij}-\partial_k v_{h,ij}^{(L)})| \lesssim (\dfrac{\rss}{L})^\beta \ \text{ on }  B_1.
\end{equation} 
Therefore combine \eqref{eqn:pp1vhildiff} and \eqref{eqn:pp1vhilbd} we obtain \begin{equation}\label{eqn:pp1vhilbdL}
    |\partial_k v_{h,ij}^{(L)}(x)| \lesssim 1 \ \text{ and } \ |\nabla \partial_k v_{h,ij}^{(L)}(x)| \lesssim 1 \ \text{ on }  B_1.
\end{equation}We are now ready to estimate \begin{equation*}
    |c_T^{(L)}-c|= |\int \Xi\cdot g|\lesssim (\fint_{B_1} |\Xi|^2)^\frac{1}{2},
\end{equation*} For the first term of $\Xi$ we use \eqref{eqn:pp1vhildiff} and \eqref{eqn:pp1gradphiconst}; the second term uses \eqref{eqn:pp1vhilbdL} and \eqref{eqn:pp1assump}; the third term is controlled by \eqref{eqn:pp1vhildiff} and \eqref{intrphigr}; the fourth term can be estimated using \eqref{eqn:pp1vhilbdL} and \eqref{eqn:pp1assump}; the fifth term can be bounded above by \eqref{eqn:pp1assump}; finally the six term is controlled by \eqref{eqn:pp1approxah} as well as \eqref{eqn:pp1gradpsi}. This finishes the proof of \eqref{eqn:pp1xicdiff}.
\smallskip

\noindent \emph{Step 4:} Conclusion. We finally compare $u^{(L)}$ defined through \eqref{eqn:finalapprox} with $\hat{u}$ defined through Corollary \ref{cor:multipole} and claim that \begin{equation*}
(\fint_{B_R}|\nabla (u^{(L)}-\hat{u})|^2)^{\frac{1}{2}}\lesssim (\dfrac{1}{L})^d(\dfrac{\rss}{L})^\beta \hspace{0.1in}\text{ for }L\ge R\ge \rss .
\end{equation*}
By Corollary \ref{cor:multipole} this leads to the conclusion of the proposition. The difference $w:=u^{(L)}-\hat{u}$ satisfies \begin{equation*}
-\nabla\cdot a\nabla w=0\hspace{0.1in}\text{ in }Q_L,\hspace{0.1in}w=w_D\text{ on }\partial Q_L
\end{equation*}
where $w_D:=(1+\phi_{i,T}^{(L)}\partial_i+\psi_{ij,T}^{(L)}\partial_{ij})u_h^{(L)}-(1+\phi_i\partial_i+\psi_{ij}\partial_{ij})u_h$. By (a slight adaptation of) the argument \eqref{co1estw} we have \begin{equation*}
(\fint_{Q_L}|\nabla w|^2)^{\frac{1}{2}}\lesssim (\fint_{Q_{\frac{5}{4}L}-Q_L} |\nabla w_D|^2)^{\frac{1}{2}}.
\end{equation*}
Hence by \eqref{co1minrad} it suffices to show \begin{equation*}
(\fint_{Q_{\frac{5}{4}L}-Q_L} |\nabla w_D|^2)^{\frac{1}{2}} \lesssim (\dfrac{1}{L})^d(\dfrac{\rss}{L})^\beta.
\end{equation*} We break down $w_D$ into eight parts, prove an estimate for each of these, and use a triangle inequality to get our desired result: \begin{align*}
\nabla w_D & =\partial_i(u_h^{(L)}-u_h)(e_i+\nabla\phi_i)+\phi_i\nabla\partial_i(u_h^{(L)}-u_h) +\partial_iu_h^{(L)}\nabla(\phi_{i,T}^{(L)}-\phi_i) \\ & \qquad +(\phi_{i,T}^{(L)}-\phi_i)\nabla\partial_i u_h^{(L)}  +\partial_{ij}(u_h^{(L)}-u_h)\nabla\psi_{ij}+\psi_{ij}\nabla\partial_{ij}(u_h^{(L)}-u_h) \\ & \qquad +\partial_{ij}u_h^{(L)}\nabla(\psi_{ij,T}^{(L)}-\psi_{ij})+(\psi_{ij,T}^{(L)}-\psi_{ij})\nabla\partial_{ij} u_h^{(L)}.
\end{align*}
The first estimate follows from \eqref{eqn:pp1gradphiconst} and \eqref{pp1uh}. The second estimate follows from \eqref{pp1uh} and \eqref{intrphigr}. The third and the fourth come from \eqref{eqn:pp1assump} and \eqref{pp1uhmag}. The fifth term is bounded due to \eqref{pp1uh} and \eqref{eqn:pp1gradpsi}. The sixth term is controlled by \eqref{pp1uh} and \eqref{eqn:psinoshiftgrowth}. The seventh term is good thanks to \eqref{pp1uhmag} and \eqref{eqn:pp1assump}. Finally the eighth term is controlled by \eqref{pp1uhmag} and \eqref{eqn:pp1psiTgrowth}.  
\qed

\subsection{Proof of Proposition \ref{prop:phismallscale}}

\noindent \emph{Step 1:} proof of \eqref{eqn:pp4eq1}. We will omit the proof for the second term since their proofs are identical. Thanks to the decomposition of $\phi$ \eqref{eqn:expphi} and $\phi_T$ \eqref{eqn:expphiT}, it suffices to prove
\begin{equation*}
\Bigl\lVert I(R\ge \rr)\Bigl(\int \eta_{R}\big|\int_0^\infty\ud t \bigl(1-\exp(-\dfrac{t}{T})\bigr)\nabla S(t)ae\big|^2\Bigr)^\frac{1}{2} \Bigr\rVert_{2-} \lesssim T^{-\frac{3}{4}} .
\end{equation*}
The proof is direct using \eqref{eqn:ppsmvpgrad} and Lemma \ref{cor:GOCor4}:  \begin{align*}
    \int_0^\infty & \ud t  (1-\exp(-\dfrac{t}{T}))\Bigl\lVert  I(R\ge \rr) \Bigl(\int \eta_{R}\big|  \nabla S(t)ae\big|^2\Bigr)^\frac{1}{2} \Bigr\rVert_{2-}  \\ & \qquad \leftstackrel{\eqref{eqn:ppsmvpgrad}}{\lesssim} \int_0^\infty \ud t(1\wedge\dfrac{t}{T})\Bigl\lVert   \Bigl(\int \eta_{\sqrt{t}}\big|  \nabla S(\dfrac{t}{2})ae\big|^2\Bigr)^\frac{1}{2} \Bigr\rVert_{2-} \\ & \qquad \leftstackrel{\eqref{eqn:semigroup}}{\lesssim} \int_0^\infty \ud t (1\wedge\dfrac{t}{T}) \frac{1}{t}(1\wedge \frac{1}{\sqrt{t}})^\frac{d}{2} \ \leftstackrel{d\ge 3}{\sim} T^{-\frac{3}{4}}.
\end{align*}
    
\smallskip  
    
\noindent \emph{Step 2:} proof of \eqref{eqn:pp4eq2}. Similar to Step 1, we need to show
\begin{equation*}
    \Bigl\lVert\Bigl(\int \eta_{\sqrt{T}}\big|\int_0^\infty\ud t \Bigl(1-\exp(-\dfrac{t}{T})\Bigr) \nabla S(t)ae\big|^2\Bigr)^\frac{1}{2} \Bigr\rVert_{2-} \lesssim T^{-\frac{3}{4}}.
\end{equation*}
For the range $t\in(0,T)$ we have as above
\begin{align*}
	\Bigl\lVert \Bigl(\int \eta_{\sqrt{T}} & \,\bigl\lvert \int_0^T\ud t \Bigl(1-\exp(-\dfrac{t}{T})\Bigr)\nabla S(t)ae  \bigr\rvert^2\Bigr)^\frac{1}{2} \Bigr\rVert_{2-}  \\ 
	 &\lesssim \int_0^T \ud t\bigl(1-\exp(-\dfrac{t}{T})\bigr)\bigl\lVert  \bigl(\int \eta_{\sqrt{T}} |\nabla S(t)ae|^2 \bigr)^\frac{1}{2} \bigr\rVert_{2-}   \leftstackrel{\eqref{eqn:semigroup}, d\ge 3}{\lesssim} \; T^{-\frac{3}{4}}.
\end{align*} For the range $t\in (T,\infty)$, we divide the integral into dyadic intervals $(2^kT,2^{k+1}T)$ and use the Meyers' inequality \cite[Lemma 7]{GO2015most} in form of \begin{equation}\label{eqn:meyers}
    \Bigl(\int \eta_{\sqrt{T}} \int_{2^kT}^{2^{k+1}T} \ud t|\nabla S(t)ae|^2 \Bigr)^\frac{1}{2} \lesssim (\dfrac{\sqrt{2^{k+1}T}}{\sqrt{T}})^{\frac{d}{2}-\gamma} \Bigl(\int \eta_{\sqrt{2^{k+1}T}} \int_{2^{k-1}T}^{2^{k+1}T} \ud t |\nabla S(t)ae|^2\Bigr)^\frac{1}{2},
\end{equation} where $\gamma=\gamma(d,\lambda)>0$. We then use the third characterization of the norm $\|\cdot \|_s$ in \eqref{eqn:equistocnorm} and the stochastic estimate \eqref{eqn:semigroup} to obtain \begin{align*}
    \Bigl\lVert  \Bigl(\int \eta_{\sqrt{2^{k+1}T}} \int_{2^{k-1}T}^{2^{k+1}T} \ud t |\nabla S(t)ae|^2\Bigr)^\frac{1}{2} \Bigr\rVert_{2-} & \leftstackrel{\eqref{eqn:equistocnorm}}{\sim} \Bigl\lVert  \int \eta_{\sqrt{2^{k+1}T}} \int_{2^{k-1}T}^{2^{k+1}T} \ud t |\nabla S(t)ae|^2 \Bigr\rVert_{1-}^\frac{1}{2} \\ & \le \Bigl(\int_{2^{k-1}T}^{2^{k+1}T} \ud t\Bigl\lVert  \int \eta_{\sqrt{2^{k+1}T}}  |\nabla S(t)ae|^2 \Bigr\rVert_{1-}\Bigr)^\frac{1}{2} \\ & \leftstackrel{\eqref{eqn:equistocnorm}}{\sim} \Bigl(\int_{2^{k-1}T}^{2^{k+1}T} \ud t\Bigl\lVert \bigl( \int \eta_{\sqrt{2^{k+1}T}}  |\nabla S(t)ae|^2 \bigr)^\frac{1}{2}\Bigr\rVert_{2-}^2\Bigr)^\frac{1}{2} \\ & \leftstackrel{\eqref{eqn:semigroup}}{\lesssim} (\int_{2^{k-1}T}^{2^{k+1}T} \ud tt^{-2-\frac{d}{2}})^\frac{1}{2} \sim (2^{k-1}T)^{-\frac{1}{2}-\frac{d}{4}}.
\end{align*}With the above two estimates we are ready to finish the proof:
\begin{align*}
    \Bigl\lVert \Bigl(\int \eta_{\sqrt{T}} & \, \lvert \int_T^\infty\ud t \bigl(1-\exp(-\dfrac{t}{T})\bigr)\nabla S(t)ae \rvert^2\Bigr)^\frac{1}{2} \Bigr\rVert_{2-}  \\ & \lesssim \sum_{k=0}^\infty \Bigl\lVert  \Bigl(\int \eta_{\sqrt{T}} |\int_{2^kT}^{2^{k+1}T}\ud t \bigl(1-\exp(-\dfrac{t}{T})\bigr)\nabla S(t)ae |^2\Bigr)^\frac{1}{2} \Bigr\rVert_{2-}  \\  & \lesssim \sum_{k=0}^\infty \Bigl(\int_{2^kT}^{2^{k+1}T} \ud t\bigl(1-\exp(-\dfrac{t}{T})\bigr)^2\Bigr)^\frac{1}{2} \Bigl\lVert  \Bigl(\int \eta_{\sqrt{T}} \int_{2^kT}^{2^{k+1}T} \ud t|\nabla S(t)ae|^2 \Bigr)^\frac{1}{2} \Bigr\rVert_{2-}  \\  & \leftstackrel{\eqref{eqn:meyers}}{\lesssim} \sum_{k=0}^\infty (2^kT)^\frac{1}{2} 2^{\frac{k+1}{2}(\frac{d}{2}-\gamma)} \Bigl\lVert  \Bigl(\int \eta_{\sqrt{2^{k+1}T}} \int_{2^{k-1}T}^{2^{k+1}T} \ud t |\nabla S(t)ae|^2\Bigr)^\frac{1}{2} \Bigr\rVert_{2-}  \\ & \lesssim \sum_{k=0}^\infty (2^kT)^\frac{1}{2}  2^{\frac{k+1}{2}(\frac{d}{2}-\gamma)} (2^{k-1}T)^{-\frac{1}{2}-\frac{d}{4}} \lesssim T^{-\frac{d}{4}}.
\end{align*} 
\qed

\subsection{Proof of Lemma \ref{lem:ahahL}}
Substituting the definition \eqref{intrhomcoeff} into \eqref{eqn:psqeq1}, where we may change the convolution kernel from Gaussian to any Schwartz function with the same scale and preserve the CLT-scaling (see \cite[Lemma 13, Step 4]{GO2015most}), and choosing $R=L$, we obtain
\begin{equation*}
    \lVert a_he_i-\int \omega_L q_i\rVert_{2-} \lesssim L^{-\frac{3}{2}}.
\end{equation*} The rest of the proof follows from a Cauchy-Schwarz, $q-q_T=a\nabla (\phi-\phi_T)$, and Proposition \ref{prop:phismallscale}, \eqref{eqn:pp4eq2}, as well as \eqref{eqn:mixingscale}: \begin{equation*}
    \lVert \int \omega_L(q-q_T)\rVert_{2-}   \le (\int \omega_L)^\frac{1}{2} \lVert (\int \omega_L|q-q_T|^2)^\frac{1}{2}\rVert_{2-}  \lesssim \lVert (\int \eta_L|\nabla (\phi-\phi_T)|^2)^\frac{1}{2}\rVert_{2-} \lesssim   \sqrt{T}^{-\frac{3}{2}}.
\end{equation*} With a triangle inequality we arrive at $\lVert a_he_i-\int \omega_L q_{i,T}\rVert_{2-}  \lesssim \sqrt{T}^{-\frac{3}{2}}$. 
\qed

\subsection{Proof of Propositions~\ref{prop:psibounds} and \ref{prop:psirss}}

We will focus on the proof of \eqref{eqn:pp6eq4} with full details and discuss the other estimates afterwards. We recall the decomposition formulas \eqref{eqn:psidecom}, \eqref{eqn:psiTdecom} for $\psi$ and $\psi_T$ and write at least formally \begin{multline}\label{eqn:psiminpsit}
    \psi-\psi_T = \underbrace{\int_0^\infty\ud t_0 \int_0^\infty \ud t_1(1-\exp(-\dfrac{t_0+t_1}{T})) S(t_0)(aS(t_1)ae-\bar{S}(t_1)\times ae) }_{\text{double integral term}} \\ - \underbrace{\int_0^\infty\ud t_0 \int_0^\infty \ud t_1  \int_0^\infty \ud t_2 (1-\exp(-\dfrac{t_0+t_1+t_2}{T})) S(t_0)\bar{S}(t_1)\times a\nabla S(t_2)ae }_{\text{triple integral term}}.
\end{multline}
The strategy is to divide these integrals into several regimes depending on the ordering of $t_0,t_1,t_2$, then estimate the layers one-by-one, using either the deterministic Lemma \ref{lem:GOlm1} when the time variable in question is not the largest, or the stochastic Propositions \ref{prop:applocqt} and \ref{prop:Stglocality} and Lemma \ref{cor:GOCor4} when the time variable being estimated is the largest of all remaining variable, with Corollary \ref{cor:upgrademvp} being used to adjust the averaging scales whenever necessary.  As we mentioned in Remark \ref{lem:psidecom}, our proof actually first shows that the integrals of \eqref{eqn:psidecom} -- \eqref{eqn:convPsi} satisfy the desired bounds \eqref{eqn:pp6eq1}, \eqref{eqn:pp6eq2}, which guarantee that the integrals are stationary functions, and hence must correspond to correctors thanks to their algebraic equivalence and uniqueness of stationary correctors.

\medskip

\noindent \emph{Proof of \eqref{eqn:pp6eq4}.} We divide the estimates into two parts, one for the double integral term, the other for the triple integral term.\\
\noindent \emph{Step 1:} estimates for the double integral term. We only estimate the term with $S(t_0)aS(t_1)ae$ since the estimates for the other term are identical. We divide the double integral into two regimes: $t_0>t_1$ and $t_1> t_0$. 

\smallskip 
\noindent\emph{Case 1.1:} $t_0>t_1$. The main ingredient we use in this case is \eqref{eqn:pp5eq2} with $g=\sqrt{t_1}aS(t_1)ae$ which by Lemma \ref{cor:GOCor4} is approximately local on scale $1 \vee \sqrt{t_1}$ with stochastic integrability $2-$, so applying the variable-coefficient semigroup $S(\frac{t_0}{2})$ will lead to the CLT factor $(\frac{1 \vee \sqrt{t_1}}{1 \vee \sqrt{t_0}})^\frac{d}{2}$ with a loss of stochastic integrability to $1-$: \begin{align*}
    \|I( &R\ge  \rr)(\int \eta_R| \int_0^\infty  \ud t_0  \int_0^{t_0} \ud t_1  (1-\exp(-\dfrac{t_0+t_1}{T}))\nabla S(t_0)aS(t_1)ae|^2)^\frac{1}{2}\|_{1-} \\ & \le \int_0^\infty  \ud t_0  \int_0^{t_0} \ud t_1  (1-\exp(-\dfrac{t_0+t_1}{T})) \|I( R\ge  \rr)(\int \eta_R|\nabla S(t_0)aS(t_1)ae|^2)^\frac{1}{2}\|_{1-} \\ & \leftstackrel{\eqref{eqn:ppsmvpgrad}}{\lesssim} \int_0^\infty  \ud t_0 (1\wedge \frac{t_0}{T})  \int_0^{t_0} \ud t_1  \|(\int\eta_{R\vee\sqrt{t_0}}|\nabla S(\dfrac{t_0}{2})aS(t_1)ae|^2)^\frac{1}{2}\|_{1-} \\ & \leftstackrel{\eqref{eqn:pp5eq2}}{\lesssim} \int_0^\infty  \ud t_0 \frac{1}{t_0} (1\wedge \frac{t_0}{T}) \int_0^{t_0} \ud t_1 (\dfrac{1 \vee \sqrt{t_1}}{1 \vee \sqrt{t_0}})^\frac{d}{2}  \frac{1}{\sqrt{t_1}}\Bigl(\|(\int\eta_{\sqrt{t_1}}|\sqrt{t_1}S(t_1)ae|^2)^\frac{1}{2}\|_{2-}+\|\bar{g}\|_{2-}\Bigr) \\ & \leftstackrel{\eqref{eqn:semigroup},\eqref{eqn:expressionbarg}}{\lesssim} \int_0^\infty  \ud t_0 \frac{1}{t_0}(1\wedge \frac{1}{\sqrt{t_0}})^\frac{d}{2} (1\wedge \frac{t_0}{T}) \int_0^{t_0} \ud t_1 (1 \vee \sqrt{t_1})^{\frac{d}{2}} \frac{1}{\sqrt{t_1}}(1 \wedge \frac{1}{\sqrt{t_1}})^\frac{d}{2} \leftstackrel{d\ge 3}{\lesssim} T^{-\frac{1}{4}}.
\end{align*}

\smallskip
\noindent\emph{Case 1.2}: $t_1>t_0$. Before we start the proof, we present here a general integration by parts formula in time which will be used repeatedly hereafter: for any $t,\tdet,T\ge 0$, \begin{equation}\label{eqn:generalibp} \begin{aligned}
    \int_0^t \ud \tau (1-\exp(-\frac{\tau+\tdet}{T}))S(\tau) g & = (1-\exp(-\frac{t+\tdet}{T}))\int_0^t \ud \tau S(\tau) g \\ & \qquad - \dfrac{1}{T}\exp(-\dfrac{\tdet}{T}) \int_0^t \ud t' \exp(-\dfrac{t'}{T}) \int_0^{t'} \ud \tau S(\tau)g.
\end{aligned}
\end{equation}To estimate this case, we first switch the order of integration of $t_0$ and $t_1$ and use \eqref{eqn:generalibp}  with $\tau=t_0, \, t=\tdet=t_1$ and $g=aS(t_1)ae$, then appeal to Lemma \ref{lem:cheating} with $ae$ playing the role of $g$:  \begin{align*}
     \Bigl\|I( &R\ge  \rr)\Bigl(\int \eta_R\bigl| \int_0^\infty  \ud t_1  \int_0^{t_1} \ud t_0  (1-\exp(-\dfrac{t_0+t_1}{T}))\nabla S(t_0)aS(t_1)ae\bigr|^2\Bigr)^\frac{1}{2}\Bigr\|_{2-} \\ & \leftstackrel{\eqref{eqn:generalibp}}{\le}  \Bigl\|I( R\ge  \rr)\Bigl(\int \eta_R\bigl|\int_0^\infty  \ud t_1    (1-\exp(-\dfrac{2t_1}{T})) \int_0^{t_1} \ud t_0\nabla S(t_0)aS(t_1)ae\bigr|^2\Bigr)^\frac{1}{2}\Bigr\|_{2-}\\ & \qquad +\Bigl\|I( R\ge  \rr)\Bigl(\int \eta_R\bigl|  \dfrac{1}{T}\int_0^\infty  \ud t_1 \exp(-\dfrac{t_1}{T}) \int_0^{t_1}\ud t_0 \exp(-\dfrac{t_0}{T})\int_0^{t_0} \ud \tau  \nabla S(\tau) a S(t_1)ae \bigr|^2\Bigr)^\frac{1}{2}\Bigr\|_{2-} \\ & \le \int_0^\infty  \ud t_1    (1-\exp(-\dfrac{2t_1}{T})) \Bigl\|I( R\ge  \rr)\Bigr(\int \eta_R| \int_0^{t_1} \ud t_0\nabla S(t_0)aS(t_1)ae|^2\Bigr)^\frac{1}{2}\Bigr\|_{2-}\\ & \qquad + \dfrac{1}{T}\int_0^\infty  \ud t_1 \exp(-\dfrac{t_1}{T})\int_0^{t_1}\ud t_0 \exp(-\dfrac{t_0}{T})\Bigl\|I( R\ge  \rr)\Bigl(\int \eta_R| \int_0^{t_0} \ud \tau  \nabla S(\tau) a S(t_1)ae |^2\Bigr)^\frac{1}{2}\Bigr\|_{2-} \\ & \leftstackrel{\eqref{eqn:nobarcheating}}{\lesssim} \int_0^\infty  \ud t_1    (1 \wedge \dfrac{t_1}{T} )\Bigl\|\Bigl(\int \eta_{\sqrt{t_1}}\bigl| \bigl(\sqrt{t_1}\nabla S(\frac{t_1}{2})ae,S(\frac{t_1}{2})ae\bigr)\bigr|^2\Bigr)^\frac{1}{2}\Bigr\|_{2-} \\ & \leftstackrel{\eqref{eqn:semigroup}}{\lesssim} \int_0^\infty  \ud t_1  (1\wedge \frac{t_1}{T})\frac{1}{\sqrt{t_1}}(1\wedge \frac{1}{\sqrt{t_1}})^{\frac{d}{2}} \ \leftstackrel{d\ge 3}{\lesssim} T^{-\frac{1}{4}}.
\end{align*} 
\smallskip

\noindent \emph{Step 2:} estimates for the triple integral term. We divide the triple integral into five regimes. In particular, we use Lemma \ref{lem:cheating} whenever $t_0<t_1$, and Lemma \ref{lem:smalltibp} whenever $t_1>t_2$ (in which case the $t_2$ integral always stays inside the stochastic norm $\lVert\cdot \rVert$). 

\smallskip

\noindent
\emph{Case 2.1}: $t_0>t_1>t_2$. In this case, notice that $\int_0^{t_1} \ud t_2 (1-\exp(-\frac{t_0+t_1+t_2}{T})) \bar{S}(t_1)\times a\nabla S(t_2)ae$ has locality scale $1 \vee \sqrt{t_1}$ by Lemma \ref{lem:smalltibp} and stochastic integrability $2-$, therefore, by Proposition \ref{prop:applocqt}, applying the semigroup $S(t_0)$ will induce stochastic cancellations by a CLT-factor $(\frac{1 \vee \sqrt{t_1}}{1 \vee \sqrt{t_0}})^d$ with a loss of stochastic integrability to $1-$:
\begin{align*}
     \int_0^\infty & \ud t_0  \int_0^{t_0}\ud t_1\Bigl\lVert I (R \ge  \rr) \Bigl(\int \eta_R\bigl| \nabla S(t_0)\int_0^{t_1}\ud t_2 (1-\exp(-\dfrac{t_0+t_1+t_2}{T}))\bar{S}(t_1)\times a\nabla S(t_2)ae\bigr|^2\Bigr)^\frac{1}{2}\Bigr\rVert_{1-}\\ & \leftstackrel{\eqref{eqn:ppsmvpgrad}}{\lesssim } \int_0^\infty \ud t_0  \int_0^{t_0}\ud t_1\Bigl\lVert   \Bigl(\int \eta_{R \vee\sqrt{t_0}}\bigl|\nabla S(\dfrac{t_0}{2}) \int_0^{t_1}\ud t_2 (1-\exp(-\dfrac{t_0+t_1+t_2}{T}))\bar{S}(t_1)\times a\nabla S(t_2)ae\bigr|^2\Bigr)^\frac{1}{2}\Bigr\rVert_{1-} \\ & \leftstackrel{\eqref{eqn:pp5eq2}}{\lesssim} \int_0^\infty \ud t_0 \frac{1}{t_0} \int_0^{t_0}\ud t_1 (\dfrac{1 \vee \sqrt{t_1}}{1 \vee \sqrt{t_0}})^\frac{d}{2}\Bigl\lVert  \Bigl(\int \eta_{\sqrt{t_1}}\bigl| \int_0^{t_1}\ud t_2 (1-\exp(-\dfrac{t_0+t_1+t_2}{T}))\bar{S}(t_1)\times a\nabla S(t_2)ae\bigr|^2\Bigr)^\frac{1}{2}\Bigr\rVert_{2-} \\ & \leftstackrel{\eqref{eqn:smalltibp}}{\lesssim}  \int_0^\infty \ud t_0 \frac{1}{t_0}(1 \wedge \frac{1}{\sqrt{t_0}})^{\frac{d}{2}} \int_0^{t_0}\ud t_1 (1 \vee \sqrt{t_1})^\frac{d}{2} (1 \wedge \dfrac{t_0}{T})\frac{1}{\sqrt{t_1}}(1 \wedge \frac{1}{\sqrt{t_1}})^{\frac{d}{2}}  \leftstackrel{d\ge 3}{\lesssim} T^{-\frac{1}{4}}.
\end{align*}Here the $\bar{g}$ term is not necessary when using \eqref{eqn:pp5eq2}, thanks to Lemma \ref{lem:smalltibp}, which we used for $t_3=t_0+t_1\sim t_0$.

\smallskip \noindent
\emph{Case 2.2}: $t_0>t_2>t_1$. Now $\sqrt{t_1}\bar{S}(t_1)\times a\nabla S(t_2)ae$ has locality scale $1 \vee \sqrt{t_2}$ (using Proposition \ref{prop:Stglocality} for $t_1$ and Lemma \ref{cor:GOCor4} for $t_2$) and stochastic integrability $2-$ (using Lemma \ref{cor:GOCor4} and Remark \ref{rmk:barSintg}), so Proposition \ref{prop:applocqt} tells us that applying $S(t_0)$ on it will result in stochastic cancellations by a CLT-factor $(\frac{1 \vee \sqrt{t_2}}{1 \vee \sqrt{t_0}})^d$ and a loss of stochastic integrability to $1-$: \begin{align*}
    \int_0^\infty & \ud t_0  \int_0^{t_0}\ud t_1 \int_{t_1}^{t_0}\ud t_2(1-\exp(-\dfrac{t_0+t_1+t_2}{T}))\Bigl\lVert I (R \ge  \rr) \Bigl(\int \eta_R\bigl| \nabla S(t_0)\bar{S}(t_1)\times a\nabla S(t_2)ae\bigr|^2\Bigr)^\frac{1}{2}\Bigr\rVert_{1-}\\ & \leftstackrel{\eqref{eqn:ppsmvpgrad}}{\lesssim} \int_0^\infty  \ud t_0 (1 \wedge \dfrac{t_0}{T}) \int_0^{t_0}\ud t_1 \int_{t_1}^{t_0}\ud t_2\Bigl\lVert  \Bigl(\int\eta_{\sqrt{t_0}}\bigl| \nabla S(\dfrac{t_0}{2})\bar{S}(t_1)\times a\nabla S(t_2)ae\bigr|^2\Bigr)^\frac{1}{2}\Bigr\rVert_{1-} \\ & \leftstackrel{\eqref{eqn:pp5eq2}}{\lesssim} \int_0^\infty  \ud t_0 (1 \wedge \dfrac{t_0}{T})\frac{1}{t_0} \int_0^{t_0}\ud t_1 \int_{t_1}^{t_0}\ud t_2 (\dfrac{1 \vee \sqrt{t_2}}{1 \vee \sqrt{t_0}})^\frac{d}{2}\frac{1}{\sqrt{t_1}} \\ & \qquad \times \Bigl(\Bigl\lVert  \Bigl(\int\eta_{\sqrt{t_2}}\bigl| \sqrt{t_1}\bar{S}(t_1)\times a\nabla S(t_2)ae\bigr|^2\Bigr)^\frac{1}{2}\Bigr\rVert_{2-}+\|\bar{G}\|_{2-}\Bigr) \\ & \leftstackrel{\eqref{eqn:GOlm1},\eqref{eqn:bargbound},\eqref{eqn:BARG}}{\lesssim} \int_0^\infty  \ud t_0 (1 \wedge \dfrac{t_0}{T})\frac{1}{t_0}\int_0^{t_0}\ud t_1 \frac{1}{\sqrt{t_1}} \int_{t_1}^{t_0}\ud t_2 (\frac{1 \vee\sqrt{t_2}}{1 \vee\sqrt{t_0}})^\frac{d}{2}\frac{1}{t_2}\Bigl(\lVert  (\int\eta_{\sqrt{t_2}}|t_2\nabla S(t_2)ae|^2)^\frac{1}{2}\rVert_{2-}+\|\bar{g}\|_{2-}\Bigr) \\ & \leftstackrel{\eqref{eqn:semigroup},\eqref{eqn:expressionbarg}}{\lesssim} \int_0^\infty  \ud t_0 (1 \wedge \dfrac{t_0}{T})\frac{1}{t_0}(1\wedge \frac{1}{\sqrt{t_0}})^\frac{d}{2} \int_0^{t_0}\ud t_1 \frac{1}{\sqrt{t_1}} \int_{t_1}^{t_0}\ud t_2 (1 \vee\sqrt{t_2})^\frac{d}{2}\frac{1}{t_2}(1 \wedge \frac{1}{\sqrt{t_2}})^{\frac{d}{2}} \ \leftstackrel{d\ge 3}{\lesssim}T^{-\frac{1}{4}}.
\end{align*}

\smallskip\noindent
\emph{Case 2.3}: $t_2>t_0>t_1$. In this case, since both $t_0$ and $t_1$ are small compared to the innermost $t_2$, we use the deterministic Lemma \ref{lem:GOlm1} on both $t_0$ and $t_1$, and the Lipschitz estimate \eqref{eqn:ppsmvpgrad} and finally Lemma \ref{cor:GOCor4} on $t_2$: \begin{align*}
     \int_0^\infty & \ud t_0  \int_0^{t_0}\ud t_1 \int_{t_0}^\infty\ud t_2(1-\exp(-\dfrac{t_0+t_1+t_2}{T}))\Bigl\lVert I (R \ge  \rr) \Bigl(\int \eta_R\bigl| \nabla S(t_0)\bar{S}(t_1)\times a\nabla S(t_2)ae\bigr|^2\Bigr)^\frac{1}{2}\Bigr\rVert_{2-} \\ & \leftstackrel{\eqref{eqn:ppsmvpgrad}}{\lesssim} \int_0^\infty  \ud t_0  \int_0^{t_0}\ud t_1 \int_{t_0}^\infty\ud t_2(1 \wedge \dfrac{t_2}{T})\Bigl\lVert I (R \vee \sqrt{t_0} \ge  \rr) \Bigl(\int \eta_{R\vee\sqrt{t_0}}\bigl| \nabla S(\dfrac{t_0}{2})\bar{S}(t_1)\times a\nabla S(t_2)ae\bigr|^2\Bigr)^\frac{1}{2}\Bigr\rVert_{2-} \\ & \leftstackrel{\eqref{eqn:GOlm1}}{\lesssim}\int_0^\infty  \ud t_0 \frac{1}{t_0} \int_0^{t_0}\ud t_1 \int_{t_0}^\infty\ud t_2(1 \wedge \dfrac{t_2}{T})\Bigl\lVert I (R \vee\sqrt{t_0} \ge  \rr) \Bigl(\int \eta_{R\vee\sqrt{ t_0}}\bigl|\bar{S}(t_1)\times a\nabla S(t_2)ae\bigr|^2\Bigr)^\frac{1}{2}\Bigr\rVert_{2-} \\ & \leftstackrel{\eqref{eqn:GOlm1}}{\lesssim}\int_0^\infty  \ud t_0 \frac{1}{t_0} \int_0^{t_0}\ud t_1 \frac{1}{\sqrt{t_1}} \int_{t_0}^\infty\ud t_2(1 \wedge \dfrac{t_2}{T})\Bigl\lVert I (R \vee \sqrt{t_0} \ge  \rr) \Bigl(\int \eta_{R\vee\sqrt{t_0}}\bigl|\nabla S(t_2)ae\bigr|^2\Bigr)^\frac{1}{2}\Bigr\rVert_{2-} \\ & \leftstackrel{\eqref{eqn:ppsmvpgrad}}{\lesssim}\int_0^\infty  \ud t_0 \frac{1}{\sqrt{t_0}}  \int_{t_0}^\infty\ud t_2(1 \wedge \dfrac{t_2}{T})\Bigl\lVert  \Bigl(\int \eta_{\sqrt{t_2}}\bigl|\nabla S(\dfrac{t_2}{2})ae\bigr|^2\Bigr)^\frac{1}{2}\Bigr\rVert_{2-} \\ & \leftstackrel{\eqref{eqn:semigroup}}{\lesssim}  \int_0^\infty  \ud t_0 \frac{1}{\sqrt{t_0}}  \int_{t_0}^\infty\ud t_2(1 \wedge \dfrac{t_2}{T}) \frac{1}{t_2}(1 \wedge \frac{1}{\sqrt{t_2}})^{\frac{d}{2}} \ \leftstackrel{d\ge 3}{\lesssim} T^{-\frac{1}{4}}.
\end{align*}

\smallskip \noindent 
\emph{Case 2.4}: $t_2>t_1> t_0$. We start the proof by using \eqref{eqn:generalibp} on $t_0$ with $t=t_1$, $t_3=t_1+t_2$, and $g=\bar{S}(t_1)\times a \nabla S(t_2)ae$,
\begin{align*}
     \int_0^\infty & \ud t_2 \int_0^{t_2} \ud t_1 \Bigl\lVert I (R \ge  \rr) \Bigl(\int \eta_R\bigl|  \int_0^{t_1}\ud t_0 (1-\exp(-\dfrac{t_0+t_1+t_2}{T}))\nabla S(t_0)\bar{S}(t_1)\times a\nabla S(t_2)ae\bigr|^2\Bigr)^\frac{1}{2}\Bigr\rVert_{2-} \\ & \leftstackrel{\eqref{eqn:generalibp}}{\le}  \int_0^\infty \ud t_2\int_0^{t_2} \ud t_1 (1-\exp(-\frac{2t_1+t_2}{T}))\Bigl\lVert I (R \ge  \rr) \Bigl(\int \eta_R\bigl| \int_0^{t_1}\ud t_0 \nabla S(t_0)\bar{S}(t_1)\times a\nabla S(t_2)ae\bigr|^2\Bigr)^\frac{1}{2}\Bigr\rVert_{2-} \\ & \qquad + \dfrac{1}{T}\int_0^\infty \ud t_2\exp(-\frac{t_2}{T}) \int_0^{t_2} \ud t_1 \exp(-\frac{t_1}{T})\int_0^{t_1}\ud t_0 \exp(-\frac{t_0}{T})  \\ & \qquad \qquad \times \Bigl\lVert I (R \ge  \rr) \Bigl(\int \eta_R\bigl|\int_0^{t_0} \ud \tau \nabla S(\tau)\bar{S}(t_1)\times a\nabla S(t_2)ae\bigr|^2\Bigr)^\frac{1}{2}\Bigr\rVert_{2-}.
\end{align*}
We now state a slight modification of Lemma \ref{lem:cheating} in the sense of \begin{multline}\label{eqn:pp7case4}\Bigl\|I(R \ge \rr) \Bigl(\int \eta_{R} \bigl|\int_0^{t_0} \ud \tau \nabla S(\tau)  \bar{S}(t_1)g\bigr|^2\Bigr)^\frac{1}{2} \Bigr\|_{2-} \\ \lesssim \Bigl\|I(R \vee \sqrt{t_1}\ge \rr)\Bigl(\int \eta_{R \vee \sqrt{t_1}} \bigl(|\sqrt{t_1}\nabla \bar{S}(\frac{t_1}{2})g|^2 + |\bar{S}(\frac{t_1}{2})g|^2\bigr)\Bigr)^\frac{1}{2}\Bigr\|_{2-}.\end{multline}In this regime we apply \eqref{eqn:pp7case4} to $g=a\nabla S(t_2)ae$, so applying Lemma \ref{lem:GOlm1} on $t_1$ and finally the Lipschitz estimate \eqref{eqn:ppsmvpgrad} and Lemma \ref{cor:GOCor4} on $t_2$ yields
\begin{align*}
    \Bigl\|I(R &\ge \rr) \Bigl(\int \eta_{R} \bigl|\int_0^{t_0} \ud \tau \nabla S(\tau)  \bar{S}(t_1)\times a \nabla S(t_2)ae\bigr|^2\Bigr)^\frac{1}{2} \Bigr\|_{2-} \\ & \leftstackrel{\eqref{eqn:pp7case4}}{\lesssim} \Bigl\|I(R \vee \sqrt{t_1}\ge \rr)\Bigl(\int \eta_{R \vee \sqrt{t_1}} \bigl(|\sqrt{t_1}\nabla \bar{S}(\frac{t_1}{2})\times a \nabla S(t_2)ae|^2 + |\bar{S}(\frac{t_1}{2})\times a \nabla S(t_2)ae|^2\bigr)\Bigr)^\frac{1}{2}\Bigr\|_{2-} \\ & \leftstackrel{\eqref{eqn:GOlm1}}{\lesssim}\frac{1}{\sqrt{t_1}}  \Bigl\|\Bigl(\int \eta_{R \vee \sqrt{t_1}} | \nabla S(t_2)ae|^2 \Bigr)^\frac{1}{2}\Bigr\|_{2-} \leftstackrel{\eqref{eqn:ppsmvpgrad}}{\lesssim}\frac{1}{\sqrt{t_1}}\Bigl\lVert  \Bigl(\int \eta_{ \sqrt{t_2}}\bigl|\nabla S(\frac{t_2}{2})ae\bigr|^2\Bigr)^\frac{1}{2}\Bigr\rVert_{2-}  \leftstackrel{\eqref{eqn:semigroup}}{\lesssim}\frac{1}{\sqrt{t_1}}\frac{1}{t_2}(1 \wedge \frac{1}{\sqrt{t_2}})^\frac{d}{2}.
\end{align*}
Combining all the arguments, we are able to conclude
\begin{align*}
     \int_0^\infty & \ud t_2 \int_0^{t_2} \ud t_1 \Bigl\lVert I (R \ge  \rr) \Bigl(\int \eta_R\bigl|  \int_0^{t_1}\ud t_0 (1-\exp(-\dfrac{t_0+t_1+t_2}{T}))\nabla S(t_0)\bar{S}(t_1)\times a\nabla S(t_2)ae\bigr|^2\Bigr)^\frac{1}{2}\Bigr\rVert_{2-} \\ & \lesssim \int_0^\infty \ud t_2 (1\wedge \frac{t_2}{T}) \frac{1}{t_2} (1\wedge \frac{1}{\sqrt{t_2}})^\frac{d}{2} \int_0^{t_2} \ud t_1 \frac{1}{\sqrt{t_1}} + \frac{1}{T} \int_0^\infty \ud t_2  \frac{1}{t_2} (1\wedge \frac{1}{\sqrt{t_2}})^\frac{d}{2} \int_0^{t_2} \ud t_1 \frac{1}{\sqrt{t_1}} \int_0^{t_1} \ud t_0 \exp(-\frac{t_0}{T}) \\ & \leftstackrel{d\ge 3}{\lesssim}T^{-\frac{1}{4}}.
\end{align*}

\smallskip \noindent
\emph{Case 2.5}: $t_1>t_0\vee t_2$. We start by using \eqref{eqn:generalibp} on $t_0$ as in Case 2.4:
\begin{align*}
     \int_0^\infty & \ud t_1 \Bigl\lVert I (R \ge  \rr) \Bigl(\int \eta_R\bigl|  \int_0^{t_1}\ud t_0\int_0^{t_1} \ud t_2 (1-\exp(-\dfrac{t_0+t_1+t_2}{T}))\nabla S(t_0)\bar{S}(t_1)\times a\nabla S(t_2)ae\bigr|^2\Bigr)^\frac{1}{2}\Bigr\rVert_{2-} \\ & \leftstackrel{\eqref{eqn:generalibp}}{\le} \int_0^\infty  \ud t_1 \Bigl\lVert I (R \ge  \rr) \Bigl(\int \eta_R\bigl|  \int_0^{t_1}\ud t_0\int_0^{t_1} \ud t_2 (1-\exp(-\dfrac{2t_1+t_2}{T}))\nabla S(t_0)\bar{S}(t_1)\times a\nabla S(t_2)ae\bigr|^2\Bigr)^\frac{1}{2}\Bigr\rVert_{2-} \\ & \qquad + \dfrac{1}{T}\int_0^\infty  \ud t_1  \int_0^{t_1}\ud t_0\exp(-\dfrac{t_1+t_0}{T})\\ & \qquad \qquad \times \Bigl\lVert I (R \ge  \rr) \Bigl(\int \eta_R\bigl|  \int_0^{t_1} \ud t_2 \exp(-\dfrac{t_2}{T})\int_0^{t_0}\nabla S(\tau)\bar{S}(t_1)\times a\nabla S(t_2)ae\bigr|^2\Bigr)^\frac{1}{2}\Bigr\rVert_{2-}.
\end{align*}
We again continue by using \eqref{eqn:pp7case4} with $g=a\nabla S(t_2)ae$, though this time we use Lemma \ref{lem:smalltibp} with $\tdet=2t_1$ on $t_1$, $t_2$: \begin{align*}
    \int_0^\infty & \ud t_1 \Bigl\lVert I (R \ge  \rr) \Bigl(\int \eta_R\bigl|  \int_0^{t_1}\ud t_0\int_0^{t_1} \ud t_2 (1-\exp(-\dfrac{t_0+t_1+t_2}{T}))\nabla S(t_0)\bar{S}(t_1)\times a\nabla S(t_2)ae\bigr|^2\Bigr)^\frac{1}{2}\Bigr\rVert_{2-}  \\ & \leftstackrel{\eqref{eqn:pp7case4}}{\lesssim} \int_0^\infty  \ud t_1 \Bigl\lVert  \Bigl(\int \eta_{\sqrt{t_1}}\bigl|  \int_0^{t_1} \ud t_2 (1-\exp(-\dfrac{2t_1+t_2}{T}))\bigl(|\sqrt{t_1}\nabla \bar{S}(\frac{t_1}{2})\times a\nabla S(t_2)ae|^2+|\bar{S}(\frac{t_1}{2})\times a\nabla S(t_2)ae|^2\bigr)\Bigr)^\frac{1}{2}\Bigr\rVert_{2-}  \\ & \qquad + \dfrac{1}{T}\int_0^\infty  \ud t_1  \int_0^{t_1}\ud t_0\exp(-\dfrac{t_0}{T}) \\ & \qquad \qquad \times \Bigl\lVert \Bigl(\int \eta_{\sqrt{t_1}}\bigl|  \int_0^{t_1} \ud t_2 \exp(-\dfrac{t_2}{T})\bigl(|\sqrt{t_1}\nabla \bar{S}(\frac{t_1}{2})\times a\nabla S(t_2)ae|^2 +|\bar{S}(\frac{t_1}{2})\times a\nabla S(t_2)ae|^2\bigr)\Bigr)^\frac{1}{2}\Bigr\rVert_{2-}  \\ & \leftstackrel{\eqref{eqn:smalltibp},\eqref{eqn:smalltibp2}}{\lesssim} \int_0^\infty (1\wedge \frac{t_1}{T}) \frac{1}{\sqrt{t_1}}(1\wedge \frac{1}{\sqrt{t_1}})^\frac{d}{2}+ \dfrac{1}{T}\int_0^\infty  \ud t_1 \int_0^{t_1}\ud t_0\exp(-\dfrac{t_0}{T}) \frac{1}{\sqrt{t_1}}(1\wedge\frac{1}{\sqrt{t_1}})^\frac{d}{2}  \ \leftstackrel{d\ge 3}{\lesssim} T^{-\frac{1}{4}}.
\end{align*}

This finishes the proof of \eqref{eqn:pp6eq4}. The estimates \eqref{eqn:pp6eq1} and \eqref{eqn:pp6eq2} follow by simple modifications of the above argument: Indeed, notice that when any of the $t_i$'s is larger than $T$, the weight $1-\exp(-\frac{t_0+t_1+t_2}{T})$ is $\sim 1$. If we replace the weights $1-\exp(-\frac{t_0+t_1+t_2}{T})$ by $1$ and repeat the integral estimates we get \eqref{eqn:pp6eq1} for $\psi$ (the calculations are actually easier since integration by parts on the weights are no longer needed). If we replace the weights by $\exp(-\frac{t_0+t_1+t_2}{T})$ and remove the gradient in front of $S(t_0)$ (which will make the bound worse by $\sqrt{T}$) we obtain \eqref{eqn:pp6eq2}.
\medskip

\noindent \emph{Proof of \eqref{eqn:pp6eq3} for $\nabla \psi$.} When any one of the $t_i$'s in representation \eqref{eqn:psidecom} is larger than $R^2$, we can use $|f_R|\lesssim (\int \eta_{R} |f|^2)^\frac{1}{2}$ to derive the bound $R^{1-\frac{d}{2}}$. For example, in the case $t_0>t_1>t_2$ of the triple integral in \eqref{eqn:psidecom} (which means $t_0>R^2$), we can estimate as \begin{align*}
    \Bigl\|I( & R\ge \rr)\Bigl(\int_{R^2}^\infty \ud t_0 \int_0^{t_0}\ud t_1 \int_0^{t_1}\ud t_2 \nabla S(t_0)\bar{S}(t_1)\times a\nabla S(t_2)ae\Bigr)_R\Bigr\|_{1-}  \\ & \lesssim   \int_{R^2}^\infty \ud t_0\int_0^{t_0}\ud t_1\Bigl\|I( R\ge \rr)\Bigl(\int \eta_R  \bigl|\int_0^{t_1}\ud t_2 \nabla S(t_0)\bar{S}(t_1)\times a\nabla S(t_2)ae\bigr|^2\Bigr)^\frac{1}{2}\Bigr\|_{1-} \\ & \leftstackrel{\eqref{eqn:ppsmvpgrad}}{\lesssim} \int_{R^2}^\infty \ud t_0\int_0^{t_0}\ud t_1\Bigl\|\Bigl(\int \eta_{\sqrt{t_0}}  \bigl|\int_0^{t_1}\ud t_2 \nabla S(\frac{t_0}{2})\bar{S}(t_1)\times a\nabla S(t_2)ae\bigr|^2\Bigr)^\frac{1}{2}\Bigr\|_{1-} \\ & \leftstackrel{\eqref{eqn:pp5eq2}}{\lesssim} \int_{R^2}^\infty \ud t_0 \frac{1}{t_0}\int_0^{t_0}\ud t_1 (\frac{1 \vee \sqrt{t_1}}{\sqrt{t_0}})^\frac{d}{2} \frac{1}{\sqrt{t_1}}\Bigl(\Bigl\|\Bigl(\int \eta_{\sqrt{t_1}}  \bigl|\int_0^{t_1}\ud t_2 \sqrt{t_1}\bar{S}(t_1)\times a\nabla S(t_2)ae\bigr|^2\Bigr)^\frac{1}{2}\Bigr\|_{2-} +\|\bar{g}\|_{2-}\Bigr) \\ & \lesssim \int_{R^2}^\infty \ud t_0 \frac{1}{t_0}\int_0^{t_0}\ud t_1 (\frac{1 \vee \sqrt{t_1}}{\sqrt{t_0}})^\frac{d}{2} \frac{1}{\sqrt{t_1}}\bigl(\bigl\|( \sqrt{t_1}\nabla \times (\opS(t_1)ae-ae))_{\sqrt{t_1}}\bigr\|_{2-} +\|\bar{g}\|_{2-}\bigr) \\ & \leftstackrel{\eqref{eqn:psqeq1},\eqref{eqn:GOlm1},\eqref{eqn:expressionbarg}}{\lesssim} \int_{R^2}^\infty \ud t_0 \frac{1}{t_0}\int_0^{t_0}\ud t_1 (\frac{1 \vee \sqrt{t_1}}{\sqrt{t_0}})^\frac{d}{2} \frac{1}{\sqrt{t_1}}(1 \wedge \frac{1}{\sqrt{t_1}})^\frac{d}{2} \leftstackrel{d\ge 3}{\sim} \, R^{1-\frac{d}{2}}.
\end{align*} We may apply the same modifications for all other cases and obtain the same bound. Hence we only need to consider the case when all $t_i$'s are in $(0,R^2)$, which will be done now.

For the double integral term of $\psi$ in \eqref{eqn:psidecom}, notice that since $aS(t_1)ae$ has locality scale $1 \vee \sqrt{t_1}$ by Lemma \ref{cor:GOCor4}, we can apply \eqref{eqn:pp5eq1} of Proposition \ref{prop:applocqt} with $aS(t_1)ae$ playing the role of $g$ and estimate
\begin{align*}
    \Big\lVert  \Bigl(  \int_0^{R^2}& \ud t_0\int_0^{R^2}\ud t_1 \nabla S(t_0)aS(t_1)ae \Bigr)_R\Big\rVert_{1-}   \\ & \le \int_0^{R^2}\ud t_1\Big\lVert  \Bigl(\int_0^{R^2}\ud t_0 \nabla S(t_0)aS(t_1)ae \Bigr)_R\Big\rVert_{1-} \\ & \leftstackrel{\eqref{eqn:pp5eq1}}{\lesssim} \int_0^{R^2} \ud t_1 (\frac{1\vee \sqrt{t_1}}{R})^\frac{d}{2} \frac{1}{\sqrt{t_1}} (\lVert (\int \eta_{\sqrt{t_1}}|\sqrt{t_1}aS(t_1)ae|^2)^\frac{1}{2}\rVert_{2-} +\|\bar{g}\|_{2-})  \\ & \leftstackrel{\eqref{eqn:semigroup},\eqref{eqn:expressionbarg}}{\lesssim} R^{-\frac{d}{2}} \int_0^{R^2}\ud t_1 (1\vee \sqrt{t_1})^\frac{d}{2}\frac{1}{\sqrt{t_1}}(1\wedge \frac{1}{\sqrt{t_1}})^\frac{d}{2} \lesssim R^{1-\frac{d}{2}}.
\end{align*}
For triple integral in \eqref{eqn:psidecom} we divide the integral into the two ranges $t_1>t_2$ and $t_1<t_2$. In the regime $t_1>t_2$, we appeal to \eqref{eqn:pp5eq1} of Proposition \ref{prop:applocqt} with $g=( \nabla \times (\opS(t_1)ae-ae))_{\sqrt{t_1}}$ which has locality scale $1 \vee \sqrt{t_1}$ (which follows from a combination of Proposition \ref{prop:Stglocality} and Lemma \ref{cor:GOCor4}). In the regime $t_1<t_2$, we also apply \eqref{eqn:pp5eq1}, this time with $g=\bar{S}(t_1)\times a \nabla S(t_2)ae$, which has locality scale $1 \vee \sqrt{t_2}$ (again using Proposition \ref{prop:Stglocality} and Lemma \ref{cor:GOCor4}). \begin{align*}
     \Big\lVert \Bigl( \int_0^{R^2}& \ud t_0\int_0^{R^2}\ud t_1\int_0^{R^2}\ud t_2 \nabla S(t_0)\bar{S}(t_1)\times a\nabla S(t_2)ae \Bigr)_R\Big\rVert_{1-}  \\ &  \leftstackrel{\eqref{eqn:cals}}{\le}\int_0^{R^2}\ud t_1 \Big\lVert  \Bigl(\int_0^{R^2} \ud t_0 \nabla S(t_0)( \nabla \times (\opS(t_1)ae-ae))_{\sqrt{t_1}}\Bigr)_R\Big\rVert_{1-}  \\ & \qquad + \int_0^{R^2}\ud t_1\int_{t_1}^{R^2} \ud t_2 \Big\lVert  \Bigl(\int_0^{R^2}\ud t_0\nabla S(t_0)\bar{S}(t_1)\times a \nabla S(t_2)ae\Bigr)_R\Big\rVert_{1-}  \\ &  \leftstackrel{\eqref{eqn:pp5eq1}}{\lesssim}  \int_0^{R^2}  \ud t_1 (\frac{1 \vee \sqrt{t_1}}{R})^\frac{d}{2} \frac{1}{\sqrt{t_1}}\Bigl(\Bigl\lVert \sqrt{t_1} ( \nabla \times (\opS(t_1)ae-ae))_{\sqrt{t_1}}\Bigr\rVert_{2-} +\|\bar{g}\|_{2-}\Bigr) \\ & \qquad + \int_0^{R^2}\ud t_1\int_{t_1}^{R^2}  \ud t_2 (\frac{1\vee \sqrt{t_2}}{R})^\frac{d}{2} \frac{1}{\sqrt{t_1}}\Bigl( \Bigl\lVert \bigl(\int \eta_{\sqrt{ t_2}}|\sqrt{t_1}\bar{S}(t_1)\times a \nabla S(t_2)ae|^2\bigr)^\frac{1}{2}\Bigr\rVert_{2-}+ \|\bar{G}\|_{2-}\Bigr)   \\ &  \leftstackrel{\eqref{eqn:psqeq1},\eqref{eqn:GOlm1},\eqref{eqn:expressionbarg}}{\lesssim} R^{-\frac{d}{2}}\int_0^{R^2} (1 \vee \sqrt{t_1})^\frac{d}{2}\frac{1}{\sqrt{t_1}}(1 \wedge \frac{1}{\sqrt{t_1}})^\frac{d}{2}\ud t_1 \\ & \qquad \qquad  + R^{-\frac{d}{2}}\int_0^{R^2}\ud t_1 \frac{1}{\sqrt{t_1}}\int_{t_1}^{R^2}\ud t_2 (1 \vee \sqrt{t_2})^\frac{d}{2} \frac{1}{t_2}\Bigl(\Bigl\lVert \bigl(\int \eta_{\sqrt{ t_2}}|t_2 \nabla S(t_2)ae|^2\bigr)^\frac{1}{2}\Bigr\rVert_{2-} +\|\bar{g}\|_{2-}\Bigr) \\ & \leftstackrel{\eqref{eqn:semigroup}}{\lesssim} R^{1-\frac{d}{2}}+ R^{-\frac{d}{2}}\int_0^{R^2}\ud t_1 \frac{1}{\sqrt{t_1}}\int_{t_1}^{R^2}\ud t_2 (1 \vee\sqrt{t_2})^\frac{d}{2}\frac{1}{t_2} (1\wedge \frac{1}{\sqrt{t_2}})^\frac{d}{2} \lesssim R^{1-\frac{d}{2}} .
\end{align*}
Using almost identical arguments we can also establish \begin{equation}\label{eqn:anabpsiR}
    \|I(R\ge \rr)(a\nabla \psi-\langle a\nabla \psi \rangle)_R\|_{1-} \lesssim R^{1-\frac{d}{2}},
\end{equation}  the only difference is that whenever we use \eqref{eqn:pp5eq1}, we use the result for $\opS(t)g- \langle \opS(t)g \rangle$. We thus do not write the argument here.

\medskip

\noindent \emph{Proof of \eqref{eqn:pp6eq1} and \eqref{eqn:pp6eq3} for $\nabla \Psi$.} The proof is based on the representation \eqref{eqn:convPsi}. For \eqref{eqn:pp6eq1} we can estimate \begin{align*}
    \|I(R & \ge \rr)(\int \eta_R|\nabla \Psi|^2)^\frac{1}{2}\|_{1-} \\ & \le \|I(R \ge \rr)(\int \eta_R|\int_0^{R^2}\ud t \nabla\bar{S}(t) \times (a\nabla\psi,a\phi,\sigma)|^2)^\frac{1}{2}\|_{1-}   \\ & \qquad + \int_{R^2}^\infty \ud t \|I(R \ge \rr)(\int \eta_R|  \nabla(\nabla\times (a\nabla\psi,a\phi,\sigma))_{\sqrt{t}}|^2)^\frac{1}{2}\|_{1-} \\ & \leftstackrel{\eqref{eqn:GOlm1}}{\lesssim} \|I(R \ge \rr)(\int \eta_R|  (a\nabla \psi,a\phi,\sigma)|^2)^\frac{1}{2}\|_{1-} + \int_{R^2}^\infty \ud t \|I(R\ge \rr) \nabla(\nabla\times (a\nabla \psi,a\phi,\sigma))_{\sqrt{t}}\|_{1-} \\ & \leftstackrel{\eqref{eqn:pp6eq1},\eqref{eqn:anabpsiR},\eqref{eqn:lm4eq1},\eqref{eqn:lm4eq2},\eqref{eqn:aphiR}}{\lesssim} 1+ \int_{R^2}^\infty  \sqrt{t}^{-2+1-\frac{d}{2}} \ud t  \leftstackrel{d\ge 3}{\lesssim} 1.
\end{align*}
The strategy for \eqref{eqn:pp6eq3} is similar, \begin{align*}
     \|I(R  \ge \rr)\nabla \Psi_R\|_{1-} & \lesssim \int_0^\infty \ud t \|I(R\ge \rr)  \nabla(\nabla\times (a\nabla\psi,a\phi,\sigma))_{\sqrt{R^2+t}}\|_{1-} \\ & \leftstackrel{\eqref{eqn:aphiR},\eqref{eqn:lm4eq2},\eqref{eqn:anabpsiR}}{\lesssim} \int_0^\infty \ud t \sqrt{R^2+t}^{-2+1-\frac{d}{2}} \ \leftstackrel{d\ge 3}{\lesssim} R^{1-\frac{d}{2}}.  
\end{align*}
\qed

\subsection{Proof of Proposition \ref{prop:truncation}}
\noindent \emph{Step 1:} proof of \eqref{eqn:phitruncation}. In fact, in this step we prove a more general statement: for any $\frac{1}{2}\le r'<r$ to be specified later, if $F_T$ satisfies the whole-space equation \begin{equation*}
     \dfrac{1}{T}F_T-\nabla\cdot a \nabla F_T= \nabla\cdot g,
\end{equation*}and $F_T^{(L)}$ satisfies the following equation in $Q_{rL}$: \begin{equation}\label{eqn:pp7ftl} \begin{cases}
    \dfrac{1}{T}F_T^{(L)}-\nabla\cdot a \nabla F_T^{(L)}= \nabla\cdot g & \mbox{ in }Q_{rL}, \\ F_T^{(L)}=0 & \mbox{ on }\partial Q_{rL}, \end{cases}
\end{equation} then for any $p<\infty$, \begin{equation}\label{eqn:pp7st1}
    \Bigl(\fint_{Q_{r'L}} \bigl|\bigl(\sqrt{T}\nabla (F_T-F_T^{(L)}),F_T-F_T^{(L)}\bigr) \bigr|^2\Bigr)^\frac{1}{2} \lesssim_{r,r',p} (\frac{\sqrt{T}}{L})^p (\int \eta_{rL} |g|^2)^\frac{1}{2}. 
\end{equation}
Choosing $g=ae$ so that $F_T=\phi_T, F_T^{(L)}=\phi_T^{(L)}$, then pick $r=2$, $r'=\frac{7}{4}$ and \eqref{eqn:pp7st1} turns into \eqref{eqn:phitruncation}. 

\smallskip

To prove \eqref{eqn:pp7st1}, we define a solution operator $S^{(L)}$ for finite domain: given $r>0$, let $v(x,t)$ be the solution for  \begin{equation}\label{eqn:truncvt}
\begin{cases}
 \partial_t v-\nabla \cdot a \nabla v=0 & \text{ in }Q_{rL}\times (0,\infty), \\ 
  v=0 & \text{ on }\partial Q_{rL}\times(0,\infty),\\
 v(t=0)=\nabla \cdot g & \text{ in }Q_{rL},  
\end{cases}
\end{equation} and we define \[S^{(L)}(t)g:=v(t).\]
Following the discussions which leads to \eqref{eqn:masparadecom}, \begin{equation*}
    F_T^{(L)}=\int_0^\infty \exp(-\dfrac{t}{T})S^{(L)}(t) g \ud t.
\end{equation*}
Following the steps of \cite[Lemma 1]{GO2015most}, we can show 
\begin{equation}\label{eqn:engyestvtrunc}
 \Bigl(\fint_{Q_{rL}} \bigl|\bigl(t\nabla S^{(L)}(t)g,\sqrt{t} S^{(L)}(t)g\bigr)\bigr|^2\Bigr)^\frac{1}{2} \lesssim (\fint_{Q_{rL}} |g|^2)^\frac{1}{2}.
\end{equation}
Now we introduce the intermediate length scale $\tilde{T}=\sqrt{T}L$. The idea is to divide the $t$-interval into $[0,\tilde{T}]$ and $[\tilde{T},\infty)$. In $[0,\tilde{T}]$, $\exp(-\frac{t}{T})\lesssim 1$ and we use Lemma \ref{lem:GOlm12} for small $t$, while in the large $t$ regime we directly use the small factor $\exp(-\frac{t}{T})$ to derive subalgebraic bound:
\begin{align*}
     \Bigl(\fint_{Q_{r'L}} &\bigl|\nabla (F_T-F_{T}^{(L)})\bigr|^2\Bigr)^\frac{1}{2} \\ & \leftstackrel{\eqref{eqn:masparadecom},\eqref{eqn:pp7ftl}}{=} \Bigl(\fint_{Q_{r'L}} \Bigl|\int_0^\infty\ud t \exp(-\dfrac{t}{T})\nabla \bigl(S(t)-S^{(L)}(t)\bigr)g\Bigr|^2\Bigr)^\frac{1}{2} \\ &    \le  \Bigl(\fint_{Q_{r'L}} \Bigl|\int_0^{\tilde{T}}\ud t \exp(-\dfrac{t}{T})\nabla \bigl(S(t)-S^{(L)}(t)\bigr)g\Bigr|^2\Bigr)^\frac{1}{2} \\
    & \qquad  +\Bigl(\fint_{Q_{r'L}} \Bigl|\int_{\tilde{T}}^\infty \ud t \exp(-\dfrac{t}{T})\nabla \bigl(S(t)-S^{(L)}(t)\bigr)g\Bigr|^2\Bigr)^\frac{1}{2} \\ &  \le \bigl(\int_0^{\tilde{T}}\ud t \exp(-\dfrac{2t}{T})\bigr)^\frac{1}{2}\Bigl(\fint_{Q_{r'L}} \int_0^{\tilde{T}} \ud t\bigl|\nabla \bigl(S(t)-S^{(L)}(t)\bigr)g\bigr|^2\Bigr)^\frac{1}{2} \\ &  \qquad +\int_{\tilde{T}}^\infty \ud t\exp(-\dfrac{t}{T})\Bigl(\fint_{Q_{r'L}} \bigl|\nabla \bigl(S(t)-S^{(L)}(t)\bigr)g\bigr|^2\Bigr)^\frac{1}{2}\\ &  \lesssim  \sqrt{T}\Bigl(\fint_{Q_{r'L}} \int_0^{\tilde{T}}\ud t \bigl|\nabla \bigl(S(t)-S^{(L)}(t)\bigr)g\bigr|^2\Bigr)^\frac{1}{2} \\ &  \qquad +\int_{\tilde{T}}^\infty \ud t\exp(-\dfrac{t}{T})\Bigl(\fint_{Q_{r'L}} \bigl|\nabla \bigl(S(t)-S^{(L)}(t)\bigr)g\bigr|^2\Bigr)^\frac{1}{2}.
\end{align*}
For the first term, note that $(S(t)-S^{(L)}(t))g$ satisfies  \begin{equation}\label{eqn:localityu-v} 
\begin{cases}
    \partial_t (S(t)-S^{(L)}(t))g-\nabla \cdot a \nabla (S(t)-S^{(L)}(t))g=0 & \text{ in }Q_{rL}\times (0,\infty), \\ 
    (S(0)-S^{(L)}(0))g=0 & \text{ in }Q_{rL} .
\end{cases}
\end{equation}
Thus, by Lemma \ref{lem:GOlm1} and a slight modification of Lemma \ref{lem:GOlm12} with $\tilde{T}$ playing the role of $T$ and $rL,r'L$ replacing $2L,L$, for any $p<\infty$, 
\begin{align*}
 \sqrt{T}\Bigl(\fint_{Q_{r'L}} & \int_0^{\tilde{T}}\ud t \bigl|\nabla \bigl(S(t)-S^{(L)}(t)\bigr)g\bigr|^2\Bigr)^\frac{1}{2} \\ & \leftstackrel{\eqref{eqn:locality}}{ \lesssim} \dfrac{\sqrt{T}}{\tilde{T}}(\dfrac{\sqrt{\tilde{T}}}{L})^{2p-1} \int_0^{\tilde{T}}\ud t(\dfrac{\sqrt{t}}{\sqrt{\tilde{T}}})^{p} \Bigl(\fint_{Q_{rL}\backslash Q_{r'L}}\bigl|(S(t)-S^{(L)}(t))g\bigr|^2\Bigr)^\frac{1}{2} \\ 
&  \lesssim  \frac{1}{\sqrt{\tilde{T}}}(\frac{\sqrt{T}}{L})^p\int_0^{\tilde{T}} \ud t(\dfrac{\sqrt{t}}{\sqrt{\tilde{T}}})^p \Bigl((\fint_{Q_{rL}} |S(t)g|^2)^\frac{1}{2}+(\fint_{Q_{rL}} |S^{(L)}(t)g|^2)^\frac{1}{2}\Bigr) \\
& \leftstackrel{\eqref{eqn:GOlm1},\eqref{eqn:engyestvtrunc}}{\lesssim}   \frac{1}{\sqrt{\tilde{T}}}(\frac{\sqrt{T}}{L})^p(\int \eta_{rL} |g|^2)^\frac{1}{2}\int_0^{\tilde{T}}\ud t \dfrac{\sqrt{t}^{p-1}}{\sqrt{\tilde{T}}^p}\lesssim (\frac{\sqrt{T}}{L})^p(\int \eta_{rL} |g|^2)^\frac{1}{2}.
\end{align*}
For the second term $$\int_{\tilde{T}}^\infty \ud t\exp(-\dfrac{t}{T})(\fint_{Q_{r'L}} |\nabla (S(t)-S^{(L)}(t))g|^2)^\frac{1}{2},$$ 
since $\exp(-\frac{t}{T})\le \exp(-\frac{\tilde{T}}{T})\lesssim_p (\frac{\sqrt{T}}{L})^p$ decays faster than any algebraic power of $\frac{\sqrt{T}}{L}$, we only need to show that $(\fint_{Q_{r'L}} |\nabla (S(t)-S^{(L)}(t))g|^2\ud t)^\frac{1}{2}$ grows at most algebraically fast in $\frac{L}{\sqrt{T}}$. This can be easily achieved by energy estimates: \begin{align*}
 \int_{\tilde{T}}^\infty &\ud t\exp(-\dfrac{t}{T}) \Bigl(\fint_{Q_{r'L}} \lvert\nabla (S(t)-S^{(L)}(t))g\rvert^2\Bigr)^\frac{1}{2} \\
 & \lesssim \int_{\tilde{T}}^\infty \ud t\exp(-\dfrac{t}{T})\Bigl(\fint_{Q_{r'L}} |\nabla S(t)g|^2\ud t\Bigr)^\frac{1}{2} +\int_{\tilde{T}}^\infty \ud t\exp(-\dfrac{t}{T})\Bigl(\fint_{Q_{r'L}} |\nabla S^{(L)}(t)g|^2\ud t\Bigr)^\frac{1}{2} \\ &  \leftstackrel{\eqref{eqn:engyestvtrunc}}{\lesssim} L^{-\frac{d}{2}}  \int_{\tilde{T}}^\infty \ud t\exp(-\dfrac{t}{T})\sqrt{t}^{\frac{d}{2}}\bigl(\fint_{Q_{\sqrt{t}}}|\nabla S(t)g|^2\ud t\bigr)^\frac{1}{2}+\int_{\tilde{T}}^\infty \ud t\exp(-\dfrac{t}{T}) t^{-1} (\int \eta_{rL} |g|^2)^\frac{1}{2}   \\ & \leftstackrel{\eqref{eqn:GOlm1}}{\lesssim} \Bigl(L^{-\frac{d}{2}}  \int_{\tilde{T}}^\infty \ud t\exp(-\dfrac{t}{T})t^{\frac{d}{4}-1}+\dfrac{T}{\tilde{T}}\exp(-\dfrac{\tilde{T}}{T})\Bigr)(\int \eta_{rL} |g|^2)^\frac{1}{2}  \lesssim (\frac{\sqrt{T}}{L})^p(\int \eta_{rL} |g|^2)^\frac{1}{2}. 
\end{align*}
The proof for $F_T-F_T^{(L)}$ without gradient is identical, and now we finish the proof of \eqref{eqn:pp7st1}. 

\smallskip

\noindent \emph{Step 2:} proof of \eqref{eqn:psiapprox}. As $\psi$'s are defined through $\sigma$'s, we need to prove the intermediate estimate \begin{equation}\label{eqn:sigtruncation}
    \Bigl(\fint_{Q_{\frac{3}{2}L}} \bigl|\bigl(\sqrt{T}\nabla (\sigma_T-\sigma_{T}^{(L)}),\sigma_T-\sigma_T^{(L)}\bigr)\bigr|^2\Bigr)^\frac{1}{2}\lesssim_p (\frac{\sqrt{T}}{L})^p.\end{equation} Let $\tilde{\sigma}_T$ be the solution of \begin{equation}\label{eqn:intmtildesig} \left\{ \begin{aligned}
    & \dfrac{1}{T}\tilde{\sigma}_T - \Delta \tilde{\sigma}_T = \nabla\cdot (q_{k,T}e_j-q_{j,T}e_k) & & \mbox{ in }Q_{\frac{7}{4}L},\\ & \tilde{\sigma}_T=0 & & \mbox{ on }\partial Q_{\frac{7}{4}L}. \end{aligned} \right.
\end{equation} Now substitute $g$ with $q_{k,T}e_j-q_{j,T}e_k$ and pick $r=\frac{7}{4}$, $r'=\frac{3}{2}$ in \eqref{eqn:pp7st1}, we obtain \begin{equation*}
    \Bigl(\fint_{Q_{\frac{3}{2}L}} \bigl|\bigl(\sqrt{T}\nabla (\sigma_T-\tilde{\sigma}_{T}),\sigma_T-\tilde{\sigma}_T\bigr)\bigr|^2\Bigr)^\frac{1}{2}\ \leftstackrel{\eqref{eqn:pp7st1}}{\lesssim} (\frac{\sqrt{T}}{L})^p (\int \eta_{\frac{7}{4}L} |q_T|^2)^\frac{1}{2}\leftstackrel{\eqref{eqn:GOlm3ell}}{\lesssim} (\frac{\sqrt{T}}{L})^p.
\end{equation*} To finish the proof of \eqref{eqn:sigtruncation} we need to pass from $\tilde{\sigma}_T$ to $\sigma_T^{(L)}$. We subtract \eqref{eqn:intmtildesig} by \eqref{eqn:algsigma}, and appeal to the standard elliptic energy estimate \begin{align*}
    (\fint_{Q_{\frac{3}{2}L}} \bigl|\bigl(\sqrt{T}\nabla (\sigma_T^{(L)}-\tilde{\sigma}_{T}),\sigma_T^{(L)}-\tilde{\sigma}_T\bigr)\bigr|^2)^\frac{1}{2} & \lesssim \sqrt{T}(\fint_{Q_{\frac{7}{4}L}} |q_T-q_{T}^{(L)}|^2)^\frac{1}{2}\\ & \lesssim \sqrt{T}(\fint_{Q_{\frac{7}{4}L}} |\nabla(\phi_T-\phi_{T}^{(L)})|^2)^\frac{1}{2}\leftstackrel{\eqref{eqn:phitruncation}}{\lesssim} (\frac{\sqrt{T}}{L})^p.
\end{align*}  The proof of \eqref{eqn:psiapprox} now follows from an identical argument as that of \eqref{eqn:sigtruncation}: pick $r=\frac{3}{2}$ and $r'=\frac{5}{4}$ in \eqref{eqn:pp7st1} and we have\begin{align*}
     \Bigl(\fint_{Q_{\frac{5}{4}L}}& \bigl|\bigl(\sqrt{T}\nabla (\psi_T-\psi_{T}^{(L)}),\psi_T-\psi_T^{(L)}\bigr)\bigr|^2\Bigr)^\frac{1}{2} \\ & \lesssim (\frac{\sqrt{T}}{L})^p\Big(\int \eta_{\frac{3}{2}L} (\phi_T^2+|\sigma_T|^2)\Big)^\frac{1}{2} +\sqrt{T}\Big(\fint_{Q_{\frac{3}{2}L}} \big((\phi_T-\phi_{T}^{(L)})^2 + |\sigma_T-\sigma_{T}^{(L)}|^2\big)\Big)^\frac{1}{2}\\ & \leftstackrel{\eqref{eqn:GOlm3ell},\eqref{eqn:phitruncation},\eqref{eqn:sigtruncation}}{\lesssim}\sqrt{T}(\frac{\sqrt{T}}{L})^p. 
\end{align*}
\qed

\subsection{Proof of Proposition \ref{prop:applocqt}}

We only include the proof for $\opS(T)g$ since the proof for $\int_0^T \ud t \nabla S(t)g$ is identical. The proof of  \eqref{eqn:pp5eq1} is divided into four steps:
\smallskip

\noindent\emph{Step 1}: Proof of \eqref{eqn:pp5eq1} when $\sqrt{T}\le r_0$.  When $r\le r_0$, using the fact that Gaussian is dominated by exponential averaging functions, we have \begin{equation}\label{eqn:pp5st1rsm} \begin{aligned}
    r^\frac{d}{2}\lVert (\opS(T)g-\langle \opS(T)g \rangle)_r \rVert_s & \lesssim  r^\frac{d}{2}\lVert (\int G_{r} |\opS(T)g|^2)^\frac{1}{2} \rVert_s \\ & \lesssim r_0^\frac{d}{2}\lVert(\int \eta_{r_0} |\opS(T)g|^2)^\frac{1}{2} \rVert_s \leftstackrel{\eqref{eqn:GOlm1}}{\lesssim} r_0^\frac{d}{2}\lVert(\int \eta_{r_0} |g|^2)^\frac{1}{2} \rVert_s.
\end{aligned}\end{equation}
When $r\ge r_0$, in view of Proposition \ref{prop:Stglocality} (in this regime the bounds on $\bar{G}$ only uses Lemma \ref{lem:GOlm1}), we may apply Lemma~\ref{lem:GOlm13} to $F=\opS(T)g-\langle \opS(T)g \rangle$ and derive with the help of \eqref{eqn:pp5st1rsm} \begin{equation}\label{eqn:pp5st1rlg}\sup_{r\ge r_0} r^\frac{d}{2}\lVert (\opS(T)g-\langle \opS(T)g\rangle)_r \rVert_s \lesssim  r_0^{\frac{d}{2}}\Bigl(\lVert  (\int \eta_{r_0} |g|^2)^\frac{1}{2}\rVert_s  + \lVert \bar{g}\rVert_s\Bigr).\end{equation}

\smallskip 
\noindent\emph{Step 2}: We claim that there exists a constant $C=C(\lambda, d)$ such that for any $\delta \ll 1$,
\begin{equation}\label{eqn:pp5st3}
    F_0:= \dfrac{1}{\sqrt{T}}\bigl(\int \eta_{\sqrt{T}}|(\phi,\sigma)_{\frac{\delta\sqrt{T}}{C}}|  ^2\bigr)^\frac{1}{2} \bigr)\le \delta
\end{equation}
implies 
\begin{equation*} 
\dfrac{1}{\sqrt{T}} \Bigl(\int \eta_{\sqrt{T}}|(\phi,\sigma)|^2\Bigr)^\frac{1}{2}\le 2\delta.
\end{equation*}Here $\phi_{\frac{\delta\sqrt{T}}{C}}$ denotes the convolution of $\phi$ and Gaussian $G_{\frac{\delta\sqrt{T}}{C}}$ (instead of the massive corrector).

We first establish the simple Caccioppoli-type estimate \begin{equation}\label{eqn:pp5engyphi}
    (\int \eta_{\sqrt{T}}|\nabla(\phi,\sigma)|^2)^\frac{1}{2}\lesssim 1+\dfrac{1}{\sqrt{T}}(\int \eta_{\sqrt{T}}|(\phi,\sigma)|^2)^\frac{1}{2}.
\end{equation}We test \eqref{intrphi} with $\eta_{\sqrt{T}}\phi$ and obtain \begin{equation*}
    \int \nabla(\eta_{\sqrt{T}}\phi)\cdot a \nabla \phi = -\int \nabla(\eta_{\sqrt{T}}\phi)\cdot ae=-\int \phi\nabla\eta_{\sqrt{T}}\cdot ae-\int \eta_{\sqrt{T}}\nabla\phi\cdot ae.
\end{equation*}We then use $|\nabla \eta_{\sqrt{T}}|\le \frac{1}{\sqrt{T}}\eta_{\sqrt{T}}$, ellipticity \eqref{eqn:intrunifell} and Young's inequality to obtain for any $\varepsilon \in (0,1)$, \begin{align*}
    \int \eta_{\sqrt{T}}|\nabla \phi|^2 & \lesssim \int \eta_{\sqrt{T}}\nabla \phi\cdot a \nabla \phi = \int \nabla(\eta_{\sqrt{T}}\phi)\cdot a \nabla \phi - \int \phi\nabla \eta_{\sqrt{T}}\cdot a \nabla\phi\\ & \lesssim \int |\phi||\nabla\eta_{\sqrt{T}}|+\int \eta_{\sqrt{T}}|\nabla\phi|+\int |\phi||\nabla \eta_{\sqrt{T}}||\nabla\phi| \\ & \lesssim \dfrac{1}{\sqrt{T}}\int \eta_{\sqrt{T}}|\phi| +\int \eta_{\sqrt{T}}|\nabla\phi| + \dfrac{1}{\sqrt{T}}\int \eta_{\sqrt{T}}|\phi||\nabla \phi| \\ & \lesssim  \varepsilon\int \eta_{\sqrt{T}}|\nabla \phi|^2 + \dfrac{1}{\varepsilon}(1+\dfrac{1}{T}\int \eta_{\sqrt{T}}\phi^2),
\end{align*}
and we finish the proof of \eqref{eqn:pp5engyphi} for $\phi$ by choosing a small $\varepsilon>0$ to absorb the first r.h.s.\ term into l.h.s. Similarly for $\sigma$ we obtain from \eqref{eqn:intrsig} \begin{align*}
    \int \eta_{\sqrt{T}}|\nabla \sigma|^2 & \lesssim \varepsilon \int \eta_{\sqrt{T}}|\nabla \sigma|^2+ \dfrac{1}{\varepsilon}(1+ \int \eta_{\sqrt{T}}|\nabla \phi|^2 + \dfrac{1}{T}\int \eta_{\sqrt{T}}\sigma^2)\\ & \lesssim \varepsilon \int \eta_{\sqrt{T}}|\nabla \sigma|^2+ \dfrac{1}{\varepsilon}(1+ \dfrac{1}{T}\int \eta_{\sqrt{T}} \phi^2 + \dfrac{1}{T}\int \eta_{\sqrt{T}}\sigma^2),
\end{align*}and we therefore finish the proof of \eqref{eqn:pp5engyphi} after absorption.

By Poincar\'{e} inequality in convolution \eqref{eqn:convpoincare}, we have for some constant $C$ which may change from line to line, \begin{align*}
     \dfrac{1}{\sqrt{T}}\bigl(\int \eta_{\sqrt{T}}\phi ^2\bigr)^\frac{1}{2} & \le \dfrac{1}{\sqrt{T}}\bigl(\int \eta_{\sqrt{T}}|\phi-\phi_{\frac{\delta\sqrt{T}}{C}}| ^2\bigr)^\frac{1}{2} +  \dfrac{1}{\sqrt{T}}\bigl(\int \eta_{\sqrt{T}}|\phi_{\frac{\delta\sqrt{T}}{C}}|  ^2\bigr)^\frac{1}{2} \\ & \leftstackrel{\eqref{eqn:pp5st3}}{\le} \dfrac{\delta}{C}\bigl(\int \eta_{\sqrt{T}}|\nabla \phi| ^2\bigr)^\frac{1}{2}+\delta \\ & \leftstackrel{\eqref{eqn:pp5engyphi}}{\le} \dfrac{\delta}{C}+\dfrac{\delta}{C\sqrt{T}}\bigl(\int \eta_{\sqrt{T}} \phi ^2\bigr)^\frac{1}{2} + \delta.
\end{align*}The term $\frac{\delta }{C\sqrt{T}}\bigl(\int \eta_{\sqrt{T}}| \phi| ^2\bigr)^\frac{1}{2}$ can be absorbed into the l.h.s., and for $C$ sufficiently large, we obtain 
\begin{equation*}
    \dfrac{1}{\sqrt{T}}(\int \eta_{\sqrt{T}}\phi^2)^\frac{1}{2}  \le 2 \delta.
\end{equation*}
 Similarly, \begin{align*}
    \dfrac{1}{\sqrt{T}}(\int \eta_{\sqrt{T}}\sigma^2)^\frac{1}{2} & \le\dfrac{1}{\sqrt{T}}\bigl(\int \eta_{\sqrt{T}}|\sigma-\sigma_{\frac{\delta\sqrt{T}}{C}}| ^2\bigr)^\frac{1}{2} +  \dfrac{1}{\sqrt{T}}\bigl(\int \eta_{\sqrt{T}}|\sigma_{\frac{\delta\sqrt{T}}{C}}|  ^2\bigr)^\frac{1}{2} \\ & \leftstackrel{\eqref{eqn:pp5st3}}{\le} \dfrac{\delta}{C}\bigl(\int \eta_{\sqrt{T}}|\nabla \sigma| ^2\bigr)^\frac{1}{2}+\delta \\ & \leftstackrel{\eqref{eqn:pp5engyphi}}{\le} \dfrac{\delta}{C}\bigl(1+ \dfrac{1}{\sqrt{T}}(\int \eta_{\sqrt{T}} \phi^2)^\frac{1}{2} + \dfrac{1}{\sqrt{T}}(\int \eta_{\sqrt{T}}\sigma^2)^\frac{1}{2}\bigr)+\delta\\ & \leftstackrel{\eqref{eqn:pp5st3}}{\le}\dfrac{\delta}{C}(1+2\delta)+\dfrac{\delta}{C\sqrt{T}}\bigl(\int \eta_{\sqrt{T}} \sigma ^2\bigr)^\frac{1}{2} + \delta ,
\end{align*}
and we finish the proof for Step 2 after absorbing $\frac{\delta}{C\sqrt{T}} (\int \eta_{\sqrt{T}}\sigma^2)^\frac{1}{2}$ into l.h.s.\ and using $\delta\ll 1$.

\smallskip 
\noindent\emph{Step 3}: Estimation of $\vertiii{\bigl(\opS(T)-\opS^h(\frac{T}{2})\opS(\frac{T}{2})\bigr)g}_{\frac{2s}{s+2}-}$. Here and for the rest of the paper we use $\vertiii{\cdot}$ to denote a norm on stationary random fields with a CLT-scaling built in: \begin{equation*}\vertiii{f}_s:=\sup_{R\ge 1}R^\frac{d}{2}\lVert (f-\langle f\rangle)_R\rVert _{s},\end{equation*} The strategy is to first estimate $ R^\frac{d}{2}\lVert ( \bigl(\opS(T)-\opS^h(\frac{T}{2})\opS(\frac{T}{2})\bigr)g)_R\rVert_{\frac{2s}{s+2}-}$ for $R\le \sqrt{T}$ and then use Lemma \ref{lem:GOlm13} for larger $R$.

Following the proof of Lemmas 4, 5, \& 6 and Proposition 1 in \cite{GO2015most}, with $(\phi,\sigma, a_h)$ playing the role of $(\phi_T,\sigma_T,a_{hT})$ (and thus the $g$ there vanishes), we obtain the following deterministic estimate of the homogenization error on the level of the flux semigroup $\opS(t)$ in a weak topology, i.e., there exists some possibly large $p=p(\lambda,d)$ such that for all $R\le \sqrt{T}$, \begin{equation}\label{eqn:pp5ff1}\begin{aligned}
    |F|:=(\dfrac{R}{\sqrt{T}})^{\frac{d}{2}}& \bigl\lvert \bigl(  \opS(T)g-\opS^h(\frac{T}{2}) \opS(\frac{T}{2})g \bigr)_R\bigr \rvert  \\ & \lesssim \delta^\frac{1}{p}\fint_{\frac{T}{4}}^{\frac{T}{2}} \ud t\fint_0^{\sqrt{t}}\ud r (\dfrac{r}{\sqrt{t}})^\frac{d}{2} \int \eta_{\sqrt{T}}\bigl\lvert\bigl(\opS(t)g-\langle \opS(t)g\rangle\bigr)_r\bigr\rvert=:\delta^\frac{1}{p}F_1,
\end{aligned}\end{equation} provided \begin{equation}\label{eqn:pp5cond}
    \dfrac{1}{\sqrt{T}} \Bigl(\int \eta_{\sqrt{T}}|(\phi,\sigma)|^2\Bigr)^\frac{1}{2}\le 2\delta.
\end{equation}
The key feature of \eqref{eqn:pp5ff1} is that the r.h.s.\ $F_1$ is of the same nature as $F$: it is a weak norm of the flux.

Using Step 2, we may replace \eqref{eqn:pp5cond} by $F_0\le 2\delta$ with the random variable $F_0$ defined in \eqref{eqn:pp5st3} that has the desired cancellation bound $\sqrt{T}^{-\frac{d}{2}}$, as we may derive from \eqref{eqn:lm4eq2},
\begin{equation}\label{eqn:F0stcbound}
    \|F_0\|_{2-} \lesssim  \delta(\delta \sqrt{T})^{-\frac{d}{2}}.
\end{equation} 
We also need the following pointwise bound on $F$:
	\begin{equation}\label{eqn:ptwsbd}
	(\dfrac{R}{\sqrt{T}})^\frac{d}{2}\Bigl|\Bigl(  \bigl(\opS(T)-\opS^h(\frac{T}{2})\opS(\frac{T}{2})\bigr)g\Bigr)_R\Bigr|\lesssim (\int \eta_{\sqrt{T}} |g|^2)^\frac{1}{2} .
	\end{equation} The proof is identical to \cite[Lemma ~ 15,~(3)]{GO2015most}, with the only difference that $ae$ is replaced by $(\int \eta_{\sqrt{T}} |g|^2)^\frac{1}{2}$, and is thus omitted.
	
For any $s_0<\frac{2s}{s+2}$, let $s_1<2$ be such that $\frac{1}{s_0}=\frac{1}{s_1}+\frac{1}{s}$, by Lemma \ref{lem:stochasticintegrability} we have for $R \le \sqrt{T}$ and $\tilde{s}_1\in (s_1,2)$, with $2\delta$ playing the role of $\delta$, \begin{equation}\begin{aligned}\label{eqn:pp5pf1}
  \lVert F\rVert_{s_0}  & \lesssim \delta^\frac{1}{p}\lVert F_1\rVert_{s_0}+ (\dfrac{\lVert F_0\rVert_{\tilde{s}_1}}{\delta})^{\frac{\tilde{s}_1}{s_1}}\lVert F \rVert_s
  \\ & \leftstackrel{\eqref{eqn:F0stcbound},\eqref{eqn:ptwsbd}}{\lesssim} \delta^\frac{1}{p}\sqrt{T}^{-\frac{d}{2}}\sup_{ t\le T} \vertiii{\opS(t)g}_{s_0}+(\delta\sqrt{T})^{-\frac{d}{2}\frac{\tilde{s}_1}{s_1}}\lVert (\int \eta_{\sqrt{T}} |g|^2)^\frac{1}{2}\rVert_s.
\end{aligned}\end{equation}
Multiply both sides by $\sqrt{T}^{\frac{d}{2}}$ and we obtain for $R\le \sqrt{T}$, \begin{equation*}
    R^\frac{d}{2}\bigl\lVert \bigl( \opS(T)g-\opS^h(\frac{T}{2})\opS(\dfrac{T}{2})g\bigr)_R\bigr\rVert_{s_0} \lesssim \delta^\frac{1}{p}\sup_{0\le t\le T} \vertiii{\opS(t)g}_{s_0}+(\delta\sqrt{T})^{-\frac{d}{2}(\frac{\tilde{s}_1}{s_1}-1)}\lVert (\int \eta_{\sqrt{T}} |g|^2)^\frac{1}{2}\rVert_s.
\end{equation*}
Since $\tilde{s}_1>s_1$, by optimizing in $\delta$, it is clear that there exists an exponent $\gamma>0$ such that for all $R\le \sqrt{T}$,\begin{equation*}
R^\frac{d}{2}\bigl\lVert \bigl(  \opS(T)g-\opS^h(\frac{T}{2})\opS(\dfrac{T}{2})g\bigr)_R\bigr\rVert _{s_0} 
\lesssim (\dfrac{1}{\sqrt{T}})^\gamma \bigl(\sup_{t\le T}\vertiii{\opS(t)g}_{s_0}+ \lVert (\int \eta_{\sqrt{T}} |g|^2)^\frac{1}{2}\rVert _s\bigr).
\end{equation*}
In order to deal with the range $R\ge \sqrt{T}$, we appeal to Lemma \ref{lem:GOlm13} with $F=\bigl(\opS(T)-\opS^h(\frac{T}{2})\opS(\frac{T}{2})\bigr)g$ minus its expectation, and $\bar{F}$ defined as in \eqref{eqn:barF} of Lemma \ref{lem:buckling}, such that $\sqrt{T}^\frac{d}{2}\|\bar{F}\|_{s_0} \lesssim \displaystyle{\sup_{t\le T} \vertiii{\opS(t)g}_{s_0}}$. Using Young's inequality for the last inequality, we get
\begin{align*}
& \vertiii{(\opS(T) -\opS^h(\frac{T}{2})  \opS(\dfrac{T}{2}))g}_{s_0}  \\ & \hspace{3em} \lesssim (\dfrac{1}{\sqrt{T}})^\gamma \Bigl(\sup_{t\le T}\vertiii{\opS(t)g}_{s_0}+ \lVert (\int \eta_{\sqrt{T}} |g|^2)^\frac{1}{2}\rVert _s\Bigr) \\ & \hspace{3em} \qquad +\Bigl((\dfrac{1}{\sqrt{T}})^\gamma \bigl(\sup_{t\le T}\vertiii{\opS(t)g}_{s_0}+ \lVert (\int \eta_{\sqrt{T}} |g|^2)^\frac{1}{2}\rVert _s\bigr)\Bigr)^{\frac{1}{2}}  \Bigl(\sup_{t\le T} \vertiii{\opS(t)g}_{s_0}+r_0^\frac{d}{2}\lVert \bar{g}\rVert _{s_0}\Bigr)^\frac{1}{2} \\ & \hspace{3em} \leftstackrel{T\ge 1}{\lesssim} (\dfrac{1}{\sqrt{T}})^{\frac{\gamma}{2}} \Bigl(\sup_{t\le T}\vertiii{\opS(t)g}_{s_0} + \lVert (\int \eta_{\sqrt{T}} |g|^2)^\frac{1}{2}\rVert _s +r_0^{\frac{d}{2}}\lVert \bar{g}\rVert_s\Bigr).\stepcounter{equation}\tag{\theequation}\label{eqn:pp5timecanc}
\end{align*}

\smallskip 
\noindent\emph{Step 4}: Conclusion of the proof by decomposing the time interval into dyadic pieces and using the semigroup property \eqref{eqn:GOlm2} of $\opS(t)$. The proof follows from \cite[Theorem 1]{GO2015most}. For any $t_1\ge r_0^2$, using \cite[Lemma 16]{GO2015most}, $\vertiii{\opS^h(t) g}_s \lesssim \vertiii{g}_s$ for all $g$ and $s$,
\begin{align*}
   \vertiii{\opS(T)g}_{s_0} &= \vertiii{\opS^h(T-t_1) \opS(t_1)g +\sum_{t_1<t\le T} \opS^h(T-t)\bigl(\opS(t)-\opS^h(\frac{t}{2})\opS(\dfrac{t}{2})\bigr)g}_{s_0} \\ & \le \vertiii{\opS^h(T-t_1) \opS(t_1)g}_{s_0} +\sum_{t_1<t\le T} \vertiii{ \opS^h(T-t)\bigl(\opS(t)-\opS^h(\frac{t}{2})\opS(\dfrac{t}{2})\bigr)g}_{s_0} \\ &  \lesssim \vertiii{\opS(t_1)g}_{s_0}+\sum_{t_1<t\le T} \vertiii{ \opS(t)g-\opS^h(\frac{t}{2})\opS(\dfrac{t}{2})g}_{s_0} \\  & \leftstackrel{\eqref{eqn:pp5timecanc}}{\lesssim } \vertiii{\opS(t_1)g}_{s_0}+\sum_{t_1<t\le T} (\dfrac{1}{\sqrt{t}})^{\frac{\gamma}{2}} \Bigl(\sup_{\tau\le t}\vertiii{\opS(\tau)g}_{s_0} + \lVert (\int \eta_{\sqrt{t}} |g|^2)^\frac{1}{2}\rVert _s +r_0^{\frac{d}{2}}\lVert \bar{g}\rVert_s\Bigr)  \\ &  \leftstackrel{\eqref{eqn:mixingscale}}{\lesssim}  \vertiii{\opS(t_1)g}_{s_0}+(\frac{1}{\sqrt{t_1}})^{\frac{\gamma}{2}} \Bigl(\sup_{0\le t\le T}\vertiii{\opS(t)g}_{s_0} + \lVert (\int \eta_{\sqrt{t_1}} |g|^2)^\frac{1}{2}\rVert _s +r_0^{\frac{d}{2}}\lVert \bar{g}\rVert_s\Bigr).
\end{align*}
Obviously, the estimate still holds if the l.h.s.\ is replaced by $\displaystyle{\sup_{t_1\le t\le T}\vertiii{\opS(t)g}_{s_0}}$. Provided $t_1\gg 1$, we may absorb
 $\displaystyle{\sup_{t_1\le t\le T}}\vertiii{\opS(t)g}_{s_0}$ into the l.h.s. Fixing $t_1$ to be a large multiple of $r_0^2\ge 1$ we may appeal \eqref{eqn:pp5eq1} for $\displaystyle{\sup_{0\le t\le t_1}}\vertiii{\opS(t)g}_{s_0}$ (which we proved in Step 1, here we used $s_0\le s$) and end up with \begin{align*}
    \sup_{0\le t\le T}\vertiii{\opS(t)g}_{s_0} & \lesssim \sup_{0\le t\le t_1}\vertiii{\opS(t)g}_{s_0} +\lVert (\int \eta_{\sqrt{t_1}} |g|^2)^\frac{1}{2}\rVert _s +r_0^{\frac{d}{2}}\lVert \bar{g}\rVert_s \\ & \leftstackrel{\eqref{eqn:pp5eq1},\eqref{eqn:mixingscale}, r_0\ge 1}{\lesssim}  r_0^{\frac{d}{2}}(\lVert (\int \eta_{r_0} |g|^2)^\frac{1}{2}\rVert_s+\lVert \bar{g}\rVert_s).
\end{align*}
This finishes the proof of \eqref{eqn:pp5eq1}. The proof of \eqref{eqn:pp5eq2} is now immediate: we use \cite[Lemma 6]{GO2015most} to obtain
\begin{equation*}
     \Bigl(\int \eta_{\sqrt{T}}|(\sqrt{T}\nabla S(T)g,S(T)g)|^2\Bigr)^\frac{1}{2} \lesssim \dfrac{1}{\sqrt{T}}\fint_{\frac{T}{2}}^T \ud t \fint_0^{\sqrt{t}}\ud r (\dfrac{r}{\sqrt{t}})^\frac{d}{2} \int \eta_{\sqrt{T}}|(\opS(t)g-\langle \opS(t)g \rangle)_r|,
\end{equation*} 
and then apply the stochastic bound \eqref{eqn:pp5eq1}. \qed

\subsection{Proof of Proposition \ref{prop:Stglocality}}

Our goal is to prove, for any $a=\tilde{a}$ in $B_{2R}$ for $R\ge r_0\vee \sqrt{T}$,\begin{equation*}
    \Bigl(\fint_{B_R} |\bigl(T\nabla \bigl(S(T)g-\tilde{S}(T)\tilde{g}\bigr), \sqrt{T}\bigl(S(T)g-\tilde{S}(T)\tilde{g}\bigr),\opS(T)g-\tilde{\opS}(T)\tilde{g}\bigr)|^2\Bigr)^\frac{1}{2}  \lesssim  \bigl(\dfrac{r_0\vee \sqrt{T}}{R}\bigr)^p \int \eta_{R} (\bar{G}+\bar{\tilde{G}}),
\end{equation*}
where $\bar{G}$ is defined through \begin{equation}\label{eqn:BARG}
        \bar{G}(0):=\dfrac{1}{\sqrt{T}}\int_0^T \ud t(\dfrac{\sqrt{t}}{\sqrt{T}})^\frac{d}{2} (\int \eta_{(r_0\vee\sqrt{t})} |S(t)g|^2)^\frac{1}{2}+(1\wedge \dfrac{r_0}{\sqrt{T}})^\frac{d}{2} \bar{g}.
    \end{equation} Here $\tilde{g}$, $\tilde{S}$, $\tilde{\opS}$ and $\bar{\tilde{G}}$ are the quantities corresponding to $a$ replaced by $\tilde{a}$. For simplicity we only present the estimates for $S(T)g-\tilde{S}(T)\tilde{g}$ since the proofs for the other two quantities are identical.
    
We decompose $S(T)g-\tilde{S}(T)\tilde{g}$ as $S(T)g-\tilde{S}(T)\tilde{g}= w_1(T)+w_2(T)$, where $w_1,w_2$ satisfy the following equations: 
	\begin{equation*}
	\partial_t w_1-\nabla \cdot a \nabla w_1 = \nabla \cdot \bigl((\tilde{a}-a)\nabla \tilde{S}(T)\tilde{g}\bigr),\hspace{0.1in}w_1(t=0)=\nabla \cdot\Bigl(I(B_{2R}^c) (g-\tilde{g})\Bigr)
	\end{equation*}
	and \begin{equation*}
	\partial_t w_2-\nabla \cdot a \nabla w_2 =0 ,\hspace{0.1in}w_2(t=0)=\nabla \cdot \Bigl(I(B_{2R})(g-\tilde{g})\Bigr).
	\end{equation*}
	 Since $a=\tilde{a}$ in $B_{2R}$, $w_1$ satisfies the assumptions of Lemma \ref{lem:GOlm12}, and hence \begin{align*}
	\Bigl(\fint_{B_R} \bigl\lvert \sqrt{T}w_1(T)\bigr\rvert^2\Bigr)^\frac{1}{2}  \lesssim  \frac{1}{\sqrt{T}}(\dfrac{\sqrt{T}}{R})^p \int_0^T \ud t(\dfrac{\sqrt{t}}{\sqrt{T}})^p (\fint_{B_{2R}\backslash B_R} w_1^2(t))^\frac{1}{2}.
	\end{align*}
	In order to estimate the r.h.s.\ we write $w_1(T)=S(T)g-\tilde{S}(T)\tilde{g}-w_2(T)$ and apply the triangle inequality, the first contribution \begin{equation*}
      \dfrac{1}{\sqrt{T}}(\dfrac{\sqrt{T}}{R})^p \int_0^T \ud t (\dfrac{\sqrt{t}}{\sqrt{T}})^p (\fint_{B_{2R}\backslash B_R} |S(t)g|^2)^\frac{1}{2}
	\end{equation*} is controlled by the first term of $\bar{G}$ using \cite[(204)]{GO2015most}. The second contribution is controlled the same way. For the third term, $\frac{1}{\sqrt{T}}(\frac{\sqrt{T}}{R})^p \int_0^T (\frac{\sqrt{t}}{\sqrt{T}})^p (\fint_{B_{2R}\backslash B_R} w_2^2(t))^\frac{1}{2}\ud t$, we use Lemma~\ref{lem:GOlm1}: 
\begin{equation*}
	\frac{1}{\sqrt{T}}(\dfrac{\sqrt{T}}{R})^p \int_0^T \ud t (\dfrac{\sqrt{t}}{\sqrt{T}})^p (\fint_{B_{2R}\backslash B_R} w_2^2(t))^\frac{1}{2} \lesssim \frac{1}{\sqrt{T}}(\dfrac{\sqrt{T}}{R})^p \int_0^T\ud t (\dfrac{\sqrt{t}}{\sqrt{T}})^p\dfrac{1}{\sqrt{t}}(\fint_{B_{2R}} |g-\tilde{g}|^2)^\frac{1}{2}.
\end{equation*}
By assumption \eqref{eqn:apploc}, the r.h.s.\ $\lesssim (\frac{\sqrt{T}}{R})^p(\frac{r_0}{R})^p \int \eta_{R}(\bar{g}+\bar{\tilde{g}})$, which is controlled by the second term of $\bar{G}$.
	
We are left with $(\fint_{B_R} |\sqrt{T} w_2(T)|^2)^\frac{1}{2} $. By Lemma \ref{lem:GOlm1} we have 
\begin{equation*}
	(\fint_{B_R} |\sqrt{T} w_2(T)|^2)^\frac{1}{2}  \lesssim (\fint_{B_{2R}} |g-\tilde{g}|^2)^\frac{1}{2}.
\end{equation*}
Again by locality assumption \eqref{eqn:apploc}, the above is controlled by the second term of $\bar{G}$. The bound \eqref{eqn:bargbound} is a direct consequence of Proposition \ref{prop:applocqt} for $\sqrt{T}\ge r_0$ and \eqref{eqn:sc5arg} (which only uses Lemma \ref{lem:GOlm1}) for $\sqrt{T}\le r_0$.
\qed

\subsection{Proof of Lemma \ref{lem:buckling}}

We first reduce the problem to estimating $(\fint_{B_R} |\opS(T)g-\tilde{\opS}(T)\tilde{g}|^2)^\frac{1}{2}$ as the difference of the other term $\opS^h(\frac{T}{2}) (\opS(\frac{T}{2})g-\tilde{\opS}(\frac{T}{2})\tilde{g})$ can be estimated similarly. Here we recall that $\tilde{g}$, $\tilde{S}$ and $\tilde{\opS}$ are the quantities correspond to $\tilde{a}$ in place of $a$. Following the proof of Proposition \ref{prop:Stglocality}, we can show for $a=\tilde{a}$ in $B_{2R}$ with $R\ge r_0\vee \sqrt{T}$,\begin{align*}
	    \bigl( \fint_{B_R} & \lvert \opS(T)g-\tilde{\opS}(T)\tilde{g} \rvert^2\bigr)^\frac{1}{2} \lesssim \bigl( \fint_{B_R}  \lvert g-\tilde{g}\rvert^2\bigr)^\frac{1}{2} + \bigl( \fint_{B_R} \lvert \int_0^T \ud t\bigl(a\nabla S(t)g-\tilde{a}\nabla \tilde{S}(t)\tilde{g}\bigr)\rvert^2\bigr)^\frac{1}{2}\\ & \leftstackrel{\eqref{eqn:apploc}}{\lesssim} (\dfrac{r_0}{R})^p\int \eta_{R} (\bar{g}+\bar{\tilde{g}})+ \Bigl(\fint_{B_R} \bigl\lvert\int_0^T \ud t \nabla (S(t)g- \tilde{S}(t)\tilde{g})\bigr\rvert^2\Bigr)^\frac{1}{2} \\ & \lesssim \dfrac{1}{\sqrt{T}}(\dfrac{\sqrt{T}}{R})^p \int_0^T\ud t (\dfrac{t}{T})^p (\fint_{B_{2R}\backslash B_R} |S(t)g|^2)^\frac{1}{2} + (\dfrac{r_0}{R})^p\int \eta_{R} \bar{g} \\ & \qquad + \dfrac{1}{\sqrt{T}}(\dfrac{\sqrt{T}}{R})^p \int_0^T \ud t(\dfrac{t}{T})^p (\fint_{B_{2R}\backslash B_R} |\tilde{S}(t)\tilde{g}|^2)^\frac{1}{2} + (\dfrac{r_0}{R})^p\int \eta_{R} \bar{\tilde{g}}.
	\end{align*} To finish the proof we have one more step\begin{equation*}
	    \dfrac{1}{\sqrt{T}}(\dfrac{\sqrt{T}}{R})^p \int_0^T \ud t(\dfrac{t}{T})^p (\fint_{B_{2R}\backslash B_R} |S(t)g|^2)^\frac{1}{2}\lesssim (\dfrac{\sqrt{T}}{R})^p\int \eta_{R} \bar{F}(a),
	\end{equation*}
	which can be shown as in the proof of \cite[Lemma 12]{GO2015most}.
\qed

\subsection{Proof of Lemma~\ref{lem:phisigb1}}

We only need to establish \eqref{eqn:lm4eq2}, \eqref{eqn:aphiR} and \eqref{eqn:lm4eq1}. Throughout the proof of this lemma we use $\phi_R:=\phi*G_R$ to denote the convolution of $\phi$ and the Gaussian $G_R$ (not to be confused with $\phi_T$), and similarly $\sigma_R$. We first address \eqref{eqn:lm4eq2}. For the $\phi$ part, we use \eqref{eqn:decomlap} in form of \begin{equation*}
    \phi=\int_0^\infty \ud t (-\Delta \phi)_{\sqrt{t}}.
\end{equation*} Therefore \begin{align*}
    \|\phi_R\|_{2-}  =\|\bigl(\int_0^\infty \ud t (-\Delta \phi)_{\sqrt{t}}\bigr)_R\|_{2-}    \le \int_0^\infty \ud t\| (-\Delta \phi)_{\sqrt{t+R^2}}\|_{2-}  &\le \int_0^\infty \ud t\frac{1}{\sqrt{t+R^2}}\| (\nabla \phi)_{\frac{\sqrt{t+R^2}}{2}}\|_{2-}  \\ &\leftstackrel{\eqref{eqn:psqeq1}}{\lesssim} \int_0^\infty \ud t \sqrt{t+R^2}^{-1-\frac{d}{2}} \; \leftstackrel{d\ge 3}{\lesssim} R^{1-\frac{d}{2}}.
\end{align*}
For $\sigma$ part we use decomposition \eqref{eqn:easydecomsig} and appeal to the bound \eqref{eqn:psqeq1} on $q$: \begin{align*}
    \|\sigma_R\| & =\|(\int_0^\infty \ud t(\nabla \times q)_{\sqrt{t}})_R\|_{2-} \le \int_0^\infty  \ud t  \|(\nabla \times q)_{\sqrt{t+R^2}}\|_{2-}\\ &  \leftstackrel{\eqref{eqn:psqeq1}}{\lesssim} \int_0^\infty  \ud t\sqrt{t+R^2}^{-1-\frac{d}{2}} \leftstackrel{d\ge 3}{\lesssim} R^{1-\frac{d}{2}}.
\end{align*}

We now turn to \eqref{eqn:aphiR}. It is slightly different from the above proof of \eqref{eqn:lm4eq2} since the starting point is the semigroup decomposition \eqref{eqn:expphi}: \begin{equation*}
    \|I(R\ge \rr) (a\phi-\langle a\phi \rangle)_R \|_{2-} \le \int_0^\infty \ud t\|I(R\ge \rr) (a S(t)ae-\langle aS(t)ae \rangle)_R \|_{2-}.
\end{equation*} For the range $\sqrt{t}\le R$, in view of Lemma \ref{cor:GOCor4} which ensures that $aS(t)ae$ is local on scale $1 \vee \sqrt{t}$, we use Lemma \ref{lem:GOlm13} on $aS(t)ae-\langle aS(t)ae \rangle$, then a Cauchy-Schwarz and the fact that Gaussian kernel is dominated by exponential kernel, and finally Lemma \ref{cor:GOCor4} to obtain  \begin{align*}
    \int_0^{R^2} \ud t\| (a S(t)ae-\langle aS(t)ae \rangle )_R \|_{2-} & \leftstackrel{\eqref{eqn:clt}}{\lesssim } \int_0^{R^2} \ud t (\dfrac{1 \vee \sqrt{t}}{R})^\frac{d}{2}\frac{1}{\sqrt{t}}\Bigl(\sup_{r\le 1 \vee \sqrt{t}}(\frac{r}{1 \vee \sqrt{t}})^\frac{d}{2}\|\sqrt{t}(a S(t)ae)_r \|_{2-}+\|\bar{g}\|_{2-}\Bigr) \\ & \lesssim \int_0^{R^2} \ud t (\dfrac{1 \vee \sqrt{t}}{R})^\frac{d}{2}\frac{1}{\sqrt{t}}\Bigl(\|(\int \eta_{(1\vee \sqrt{t})}|\sqrt{t}a S(t)ae|^2)^\frac{1}{2}\|_{2-}+\|\bar{g}\|_{2-}\Bigr) \\ & \leftstackrel{\eqref{eqn:semigroup},\eqref{eqn:expressionbarg}}{\lesssim} R^{-\frac{d}{2}}\int_0^{R^2} \ud t (1 \vee\sqrt{t})^\frac{d}{2}\frac{1}{\sqrt{t}}(1\wedge \frac{1}{\sqrt{t}})^\frac{d}{2} \sim R^{1-\frac{d}{2}}.
\end{align*}
For the range $\sqrt{t}\ge R$, we directly use the exponential kernel to control the Gaussian kernel, then appeal to \eqref{eqn:ppsmvpnograd} and Lemma \ref{cor:GOCor4}: \begin{align*}
    \int_{R^2}^\infty \ud t\|I(R\ge \rr) (a S(t)ae-\langle aS(t)ae\rangle)_R \|_{2-} & \lesssim  \int_{R^2}^\infty \ud t\|I(R\ge \rr) (\int \eta_{R} | S(t)ae|^2)^\frac{1}{2}\|_{2-} \\ & \leftstackrel{\eqref{eqn:ppsmvpnograd}}{\lesssim}  \int_{R^2}^\infty \ud t\| (\int \eta_{\sqrt{t}} | S(\frac{t}{2})ae|^2)^\frac{1}{2}\|_{2-} \\ & \leftstackrel{\eqref{eqn:semigroup}}{\lesssim} \int_{R^2}^\infty \ud t \sqrt{t}^{-1-\frac{d}{2}} \; \leftstackrel{d\ge 3}{\sim} R^{1-\frac{d}{2}}.
\end{align*}

We finally prove \eqref{eqn:lm4eq1}. The $\phi$ part uses Poincar\'{e} inequality in convolution \eqref{eqn:convpoincare}, and \eqref{eqn:lm4eq2} for $r=1$, \begin{align*}
     \lVert (\int \eta_R |\phi |^2)^\frac{1}{2}\rVert_{2-} & \lesssim \lVert (\int \eta_R |\phi-\phi_1 |^2)^\frac{1}{2}\rVert_{2-}   + \lVert (\int\eta_R |\phi_1 |^2)^\frac{1}{2}\rVert_{2-}  \\ & \lesssim \lVert (\int \eta_R |\nabla \phi|^2)^\frac{1}{2}\rVert_{2-}   + \lVert \phi_1 \rVert_{2-}  \leftstackrel{\eqref{eqn:psqeq2},\eqref{eqn:lm4eq2}}{\lesssim}  1.
\end{align*}
For $\sigma$ part, we again use decomposition \eqref{eqn:easydecomsig}
\begin{align*}
     \lVert (\int \eta_R |\sigma |^2)^\frac{1}{2}\rVert_{2-} & =  \lVert (\int \eta_R |\int_0^\infty\ud t (\nabla \times q)_{\sqrt{t}}|^2)^\frac{1}{2}\rVert_{2-}  \\ &  \lesssim \lVert (\int \eta_R |\int_0^1\ud t \bar{S}(t)\times q |^2)^\frac{1}{2}\rVert_{2-} + \int_1^\infty\ud t\lVert (\fint_{B_R} |(\nabla \times q)_{\sqrt{t}}|^2)^\frac{1}{2}\rVert_{2-}  \\ &  \leftstackrel{\eqref{eqn:GOlm1}}{\lesssim} \lVert (\int \eta_R | q |^2)^\frac{1}{2}\rVert_{2-}+ \int_1^\infty\ud t\lVert (\nabla \times q)_{\sqrt{t}}\rVert_{2-}  \\ & \leftstackrel{\eqref{eqn:psqeq1},\eqref{eqn:psqeq2}}{\lesssim} 1+ \int_1^\infty\ud t \sqrt{t}^{-1-\frac{d}{2}} \leftstackrel{d\ge 3}{\lesssim} 1.
\end{align*}
\qed

\subsection{Proof of Lemma~\ref{lem:cheating}}

If $R \ge \sqrt{t_0}$, the argument is direct: \begin{align*}
    \|I(R \ge \rr) (\int \eta_{R} |\int_0^{t_0} \ud \tau \nabla S(\tau)  aS(t_1)g|^2)^\frac{1}{2} \|_{2-} & \leftstackrel{\eqref{eqn:GOlm1}}{\lesssim}  \|I(R  \ge \rr) (\int \eta_{R} |S(t_1)g|^2)^\frac{1}{2} \|_{2-} \\ & \leftstackrel{\eqref{eqn:ppsmvpnograd}}{\lesssim}\|I(R \vee \sqrt{t_1}\ge \rr)(\int \eta_{R\vee\sqrt{t_1}} |S(\frac{t_1}{2})g|^2)^\frac{1}{2} \|_{2-}.
\end{align*} We now assume $\sqrt{t_0}\ge R$. The starting point is the following: let $u$ satisfy the equation\begin{equation*}
    -\nabla \cdot a (\nabla u+g)=f,
\end{equation*} then we have \begin{equation}\label{eqn:schauder}
   \sup_{r\in [\rr,R]}\int \eta_r |\nabla u|^2 \\ \lesssim \int \eta_R |\nabla u|^2 +\sup_{r\in [\rr,R]} \int \eta_r \Bigl(|g|^2+ R^2 |\nabla g|^2+ R^2f^2\Bigr).
\end{equation} By a slight modification of \cite[Corollary 3]{gloria2014regularity} as well as Poincar\'e inequality we have
 \begin{equation}\label{eqn:noppsschauder}\begin{aligned}
   \sup_{r\in[\rr,R]} \fint_{B_r} |\nabla u+g|^2 & \lesssim \fint_{B_R} |\nabla u+g|^2 + \sup_{r\in [\rr,R]}\dfrac{R^2}{r^2} \fint_{B_r} |g-\fint_{B_r} g|^2 +R^2\sup_{r\in [\rr,R]} \fint_{B_r} f^2 \\ & \lesssim \fint_{B_R} |\nabla u+g|^2 + R^2\sup_{r\in [\rr,R]} \fint_{B_r} (|\nabla g|^2+f^2).
\end{aligned}\end{equation}
The proof of \eqref{eqn:schauder} then follows from post-processing from $\fint_{B_r}$ to $\int \eta_r$ (see the proof of Lemma \ref{lem:parabolicmvp} for the argument) and a simple triangle inequality. Now, since the function $w:=\int_0^{t_0} \ud \tau  S(\tau) a S(t_1)g$ satisfies the elliptic equation \begin{equation*}
    -\nabla \cdot a (\nabla w +S(t_1)g)=- S(t_0)a S(t_1)g,
\end{equation*} using \eqref{eqn:schauder} with $\sqrt{t_0},\ S(t_1)g, \ -S(t_0)a S(t_1)g$ playing the roles of $R, g, f$ respectively, as well as Corollary \ref{cor:upgrademvp}, we obtain \begin{align*}
      \|I(R &\ge \rr) (\int \eta_{R} |\int_0^{t_0} \ud \tau\nabla S(\tau) a S(t_1)g|^2)^\frac{1}{2} \|_{2-} \\ & \leftstackrel{\eqref{eqn:schauder}}{\lesssim}  \|I(\sqrt{t_0} \ge \rr) (\int \eta_{\sqrt{t_0}} |\int_0^{t_0} \ud \tau \nabla S(\tau) a S(t_1)g|^2)^\frac{1}{2} \|_{2-}  +  \|\sup_{r\in[\rr,\sqrt{t_0}]} (\int \eta_r |S(t_1)g|^2)^\frac{1}{2}\|_{2-} \\ & \qquad+ \sqrt{t_0}\|\sup_{r\in[\rr,\sqrt{t_0}]} (\int \eta_r |\nabla S(t_1)g|^2)^\frac{1}{2}\|_{2-} \\ & \qquad + \sqrt{t_0}\|I(\sqrt{t_0}\ge \rr)\sup_{r\in[\rr,\sqrt{t_0}]}  (\int \eta_r|S(t_0)a S(t_1)g|^2)^\frac{1}{2}\|_{2-} \\ & \leftstackrel{\eqref{eqn:detppsmvpgrad},\eqref{eqn:detppsmvpnograd}}{\lesssim}\|I(\sqrt{t_0} \ge \rr) (\int \eta_{\sqrt{t_0}} |\int_0^{t_0} \ud \tau \nabla S(\tau) a S(t_1)g|^2)^\frac{1}{2} \|_{2-} + \|I(\sqrt{t_1}\ge \rr)(\int \eta_{\sqrt{t_1}} |S(\frac{t_1}{2})g|^2)^\frac{1}{2}\|_{2-}\\ & \qquad+  \sqrt{t_0}\| I(\sqrt{t_1}\ge \rr)(\int \eta_{\sqrt{t_1}} |\nabla S(\frac{t_1}{2})g|^2)^\frac{1}{2}\|_{2-}  + \sqrt{t_0}\|I(\sqrt{t_0}\ge \rr)  (\int \eta_{\sqrt{t_0}}|S(\frac{t_0}{2})a S(t_1)g|^2)^\frac{1}{2}\|_{2-} \\ & \leftstackrel{\eqref{eqn:GOlm1}}{\lesssim} \|I(\sqrt{t_0} \ge \rr) (\int \eta_{\sqrt{t_0}} | S(t_1)g|^2)^\frac{1}{2} \|_{2-} + \|I(\sqrt{t_1}\ge \rr)(\int \eta_{\sqrt{t_1}} |S(\frac{t_1}{2})g|^2)^\frac{1}{2}\|_{2-}\\ & \qquad +  \sqrt{t_0}\| I(\sqrt{t_1}\ge \rr)(\int \eta_{\sqrt{t_1}} |\nabla S(\frac{t_1}{2})g|^2)^\frac{1}{2}\|_{2-} \\ & \leftstackrel{\eqref{eqn:ppsmvpnograd}}{\lesssim} \|I(\sqrt{t_1}\ge \rr)(\int \eta_{\sqrt{t_1}} |S(\frac{t_1}{2})g|^2)^\frac{1}{2}\|_{2-}+  \sqrt{t_0}\| I(\sqrt{t_1}\ge \rr)(\int \eta_{\sqrt{t_1}} |\nabla S(\frac{t_1}{2})g|^2)^\frac{1}{2}\|_{2-}.
\end{align*}
\qed

\subsection{Proof of Lemma~\ref{lem:smalltibp}}
We first prove \eqref{eqn:smalltibp}. The starting point is to show for $R\ge 1$\begin{equation}\label{eqn:lm65ae}
     \|(ae-\langle ae \rangle)_R\|_{2-} \lesssim R^{-\frac{d}{2}}.
\end{equation}The proof uses \eqref{eqn:lm510} for $ae-\langle ae \rangle$ that is mean-zero and exactly local on scale $r=1$:\begin{equation*}
    \lVert (ae-\langle ae \rangle)_R \rVert_{2-} \leftstackrel{\eqref{eqn:lm510}}{\lesssim} (\int G_R^2)^\frac{1}{2}\lVert ae-\langle ae\rangle \rVert_{2-}\lesssim R^{-\frac{d}{2}}.
\end{equation*}
The next step is to prove for $t_2\le t_1$,
\begin{equation}\label{eqn:lm65calst2}
    \lVert \nabla \times (\opS(t_2)ae-ae)_{\sqrt{t_1}}\rVert_{2-} \lesssim \frac{1}{\sqrt{t_1}}(1 \wedge \frac{1}{\sqrt{t_1}})^{\frac{d}{2}}.
\end{equation}
If $t_1 \le 1$ the proof is completely deterministic. Using a Cauchy-Schwarz and that the Gaussian kernel can be dominated by exponential kernel and Lemma \ref{lem:GOlm1} we obtain \begin{align*}
    \lVert \nabla \times (\opS(t_2)ae-ae)_{\sqrt{t_1}}\rVert_{2-} & \lesssim   \frac{1}{\sqrt{t_1}}\lVert (a\int_0^{t_2} \ud \tau \nabla S(\tau)ae)_{\sqrt{t_1}/2}\rVert_{2-} \\ & \lesssim \frac{1}{\sqrt{t_1}}\Bigl\lVert \Bigl(\int \eta_{\frac{\sqrt{t_1}}{2}}\bigl|a\int_0^{t_2} \ud \tau \nabla S(\tau)ae\bigr|^2\Bigr)^\frac{1}{2}\Bigr\rVert_{2-} \lesssim \frac{1}{\sqrt{t_1}}.
\end{align*}
In the case $t_1\ge 1$ the proof follows from \eqref{eqn:psqeq1}: \begin{equation*}
     \lVert \nabla \times (\opS(t_2)ae)_{\sqrt{t_1}}\rVert_{2-} =  \lVert \nabla \times (\opS(t_2)ae-\langle \opS(t_2)ae \rangle)_{\sqrt{t_1}}\rVert_{2-}  \leftstackrel{\eqref{eqn:psqeq1}}{\lesssim} \sqrt{t_1}^{-1-\frac{d}{2}},
\end{equation*}
and we prove \eqref{eqn:lm65calst2} using a triangle inequality and \eqref{eqn:lm65ae}. We are now ready to prove \eqref{eqn:smalltibp}. For simplicity we only prove estimates for the second term as the proof for the other term is identical. Using integration by parts \eqref{eqn:generalibp}, \begin{align*}
     \int_0^{t_1}\ud t_2 & (1-\exp(-\dfrac{\tdet +t_2}{T}))\bar{S}(t_1)\times a\nabla S(t_2)ae  \\ & = (1-\exp(-\dfrac{\tdet+t_1}{T})) \int_0^{t_1}\ud t_2 \bar{S}(t_1)\times a\nabla S(t_2)ae \\ & \qquad - \dfrac{1}{T}  \exp(-\dfrac{\tdet}{T})\int_0^{t_1}\ud t_2 \exp(-\dfrac{t_2}{T})\int_0^{t_2}\ud \tau\bar{S}(t_1)\times a\nabla S(\tau)ae \\ & 
     \leftstackrel{\eqref{eqn:cals},\eqref{eqn:barSgaussian}}{=} (1-\exp(-\dfrac{\tdet+t_1}{T})) \nabla \times (\opS(t_1)ae-ae)_{\sqrt{t_1}} \\ & \qquad - \dfrac{1}{T}  \exp(-\dfrac{\tdet}{T})\int_0^{t_1}\ud t_2 \exp(-\dfrac{t_2}{T})\nabla\times (\opS(t_2)ae-ae)_{\sqrt{t_1}}.
\end{align*}
Therefore using \eqref{eqn:psqeq1} and $\tdet \ge t_1$, 
\begin{align*}
  \bigl\lVert \int_0^{t_1} & \ud t_2 \bigl(1-\exp(-\dfrac{\tdet +t_2}{T})\bigr)\bar{S}(t_1)\times a\nabla S(t_2)ae \bigr\rVert_{2-}  \\ & \le (1-\exp(-\dfrac{\tdet+t_1}{T}))\lVert  (\nabla\times (\opS(t_1)ae-ae))_{\sqrt{t_1}} \rVert_{2-}\\ & \qquad + \dfrac{1}{T} \exp(-\dfrac{\tdet}{T}) \int_0^{t_1}\ud t_2 \exp(-\dfrac{t_2}{T})\lVert (\nabla \times (\opS(t_2)ae-ae))_{\sqrt{t_1}}\rVert_{2-} \\ & \leftstackrel{\eqref{eqn:lm65calst2}}{\lesssim}   (1 \wedge \dfrac{t_2}{T})\frac{1}{\sqrt{t_1}}(1 \wedge \frac{1}{\sqrt{t_1}})^{\frac{d}{2}} +  \exp(-\dfrac{\tdet}{T})(1 \wedge \dfrac{t_2}{T}) \frac{1}{\sqrt{t_1}}(1 \wedge \frac{1}{\sqrt{t_1}})^{\frac{d}{2}} \\ & \lesssim (1 \wedge \dfrac{t_2}{T})\frac{1}{\sqrt{t_1}}(1 \wedge \frac{1}{\sqrt{t_1}})^{\frac{d}{2}}. 
\end{align*}
The proof of \eqref{eqn:smalltibp2} is similar, again we only write down the proof for the second term: \begin{align*}
  \bigl\lVert \int_0^{t_1} & \ud t_2 \exp(-\dfrac{t_2}{T})\bar{S}(t_1)\times a\nabla S(t_2)ae \bigr\rVert_{2-}  \\ & \le \exp(-\dfrac{t_1}{T})\lVert  (\nabla\times (\opS(t_1)ae-ae))_{\sqrt{t_1}} \rVert_{2-}\\ & \qquad + \dfrac{1}{T}  \int_0^{t_1}\ud t_2 \exp(-\dfrac{t_2}{T})\lVert (\nabla \times (\opS(t_2)ae-ae))_{\sqrt{t_1}}\rVert_{2-} \\ & \leftstackrel{\eqref{eqn:lm65calst2}}{\lesssim}  \frac{1}{\sqrt{t_1}}(1 \wedge \frac{1}{\sqrt{t_1}})^{\frac{d}{2}} +  (1 \wedge \dfrac{t_2}{T}) \frac{1}{\sqrt{t_1}}(1 \wedge \frac{1}{\sqrt{t_1}})^{\frac{d}{2}} \\ & \lesssim \frac{1}{\sqrt{t_1}}(1 \wedge \frac{1}{\sqrt{t_1}})^{\frac{d}{2}}. 
\end{align*}

For the proof of approximate locality \eqref{eqn:smalltibp3}, in view of Proposition \ref{prop:Stglocality} applied to $\bar{S}(t_1)$, it suffices to prove that $\int_0^{t_1}  \ud t_2 (1-\exp(-\frac{\tdet +t_2}{T})) a\nabla S(t_2)ae$ is approximately local on scale $1 \vee \sqrt{t_1}$ relative to some stationary $\bar{g}_0$ with $\|\bar{g}_0\|_{2-}\lesssim (1 \wedge \frac{\tdet}{T})(1 \wedge \frac{1}{\sqrt{t_1}})^{\frac{d}{2}}$. Now suppose $a=\tilde{a}$ in $B_{2R}$ for some $R\ge 1 \vee \sqrt{t_1}$, notice that $S(t)ae-\tilde{S}(t)\tilde{a}e$ satisfies the conditions of Lemma \ref{lem:GOlm12}, we can estimate \begin{align*}
    \Bigl(\fint_{B_R}&  \bigl|\int_0^{t_1}\ud t_2 (1-\exp(-\frac{\tdet +t_2}{T})) \nabla (S(t_2)ae-\tilde{S}(t_2)\tilde{a}e) \bigr|^2\Bigr)^\frac{1}{2} \\ & \le \int_0^{t_1}\ud t_2 (1-\exp(-\frac{\tdet +t_2}{T}))\Bigl(\fint_{B_R} \bigl| \nabla (S(t_2)ae-\tilde{S}(t_2)\tilde{a}e) \bigr|^2\Bigr)^\frac{1}{2} \\ & \leftstackrel{\eqref{eqn:locality}}{\lesssim} \dfrac{1}{R}\int_0^{t_1}\ud t_2 (1\wedge \frac{t_3}{T})(\frac{\sqrt{t_2}}{R})^{p+2}\fint_0^{t_2} \ud \tau (\frac{\sqrt{\tau}}{\sqrt{t_2}})^p\Bigl(\fint_{B_{2R}\backslash B_R}(S(\tau)ae-\tilde{S}(\tau)ae )^2  \Bigr)^\frac{1}{2} \\ & \lesssim  \frac{1}{R^{p+3}}(1 \wedge \frac{t_3}{T}) \int_0^{t_1}\ud \tau \tau^\frac{p}{2} \Bigl(\bigl(\fint_{B_{2R}\backslash B_R}|S(\tau)ae|^2\bigr)^\frac{1}{2}+\bigl(\fint_{B_{2R}\backslash B_R}|\tilde{S}(\tau)ae|^2\bigr)^\frac{1}{2}  \Bigr) \int_{\tau}^{t_1} \ud t_2 \\ & \sim \frac{1}{R}(1 \wedge \frac{t_3}{T})(\frac{\sqrt{t_1}}{R})^2 (\frac{1\vee \sqrt{t_1}}{R})^{p} \int_0^{t_1}\ud \tau (\frac{\sqrt{\tau}}{1 \vee \sqrt{t_1}})^{p} \Bigl(\bigl(\fint_{B_{2R}\backslash B_R}|S(\tau)ae|^2\bigr)^\frac{1}{2}+\bigl(\fint_{B_{2R}\backslash B_R}|\tilde{S}(\tau)ae|^2\bigr)^\frac{1}{2}  \Bigr).
\end{align*}
Now using the fact that $R\ge 1\vee \sqrt{t_1}$, as well as \cite[(204)]{GO2015most}, we derive that $\int_0^{t_1}  \ud t_2 (1-\exp(-\frac{\tdet +t_2}{T})) a\nabla S(t_2)ae$ is approximately local on scale $1 \vee \sqrt{t_1}$ relative to some stationary $\bar{g}_0$ with \begin{equation*}
    \bar{g}_0(0):=\frac{1}{1 \vee\sqrt{t_1}}(1 \wedge \frac{t_3}{T}) \int_0^{t_1}\ud \tau (\frac{\sqrt{\tau}}{1 \vee \sqrt{t_1}})^{p}  (\int \eta_{\sqrt{\tau}} |S(\tau)ae|^2)^\frac{1}{2},
\end{equation*} and we finish the proof after applying Lemma \ref{cor:GOCor4} on $\bar{g}_0$. \qed

\section*{Acknowledgement}
The work of JL and LW is supported in part by the National Science Foundation via grants DMS-1454939 and DMS-2012286. We would like to thank Marius Lemm for his suggestion on references. LW would also like to thank the Max Planck Institute for their generous support, and thank Peter Bella, Nicolas Clozeau, Antoine Gloria, Yu Gu and Jim Nolen for helping with the understanding of stochastic homogenization, and Yingzhou Li for assistance on numerical implementation.

\appendix
\section{Properties of stochastic norm $\|\cdot\|_s$}
In this section we collect properties of the stochastic norm $\|\cdot\|_s$ defined in \eqref{eqn:defsnorm}. The first lemma presents two other equivalent formulations that are convenient for this work.
\begin{lemma}\label{lem:strexpnorm} (\cite[Lemma 3.7]{ledoux2013probability})
	For any $s\in (0,2]$ and the corresponding $c$ defined in \eqref{eqn:defsnorm}, \begin{equation}\label{eqn:equistocnorm}
    \lVert F\rVert _s \sim_s \sup_{m\in \N}\dfrac{\langle |F|^m \rangle^\frac{1}{m}}{m^\frac{1}{s}}\sim_s \inf\bigl\{M>0: ~ \log \langle\exp\bigl((\dfrac{|F|}{M})^s\bigr)\rangle\le 1  \bigr\}.
\end{equation}
\end{lemma}
\begin{proof}
  We first state the following auxiliary result that will be used later: there exist constants $c_0,c_1>0$ such that for all positive integers $k$, \begin{equation}\label{eqn:factorial}
    (c_0k)^k \le k!\le (c_1k)^k.
\end{equation}
Our proof consists of three steps, with each one showing one expression of \eqref{eqn:equistocnorm} being dominated by another. \\
\noindent \emph{Step 1:} middle $\lesssim$ left. If $\lVert F \rVert_s \le 1$, or equivalently \[\langle\exp\bigl((|F|+c)^s\bigr)\rangle \le 1+\exp(c^s),\] then we will show there exists a constant $C$, only depending on $s$, such that $\langle |F|^m\rangle^\frac{1}{m} \le Cm^\frac{1}{s}$ for all $m\ge 1$. By the Chebyshev inequality, \begin{equation*}
    \langle I(|F|>\lambda) \rangle \le \exp(-(\lambda+c)^s)\langle \exp((|F|+c)^s) \rangle \le (1+\exp(c^s))\exp(-\lambda^s).
\end{equation*}
Therefore, using the formula for non-negative random variable $F$ \begin{equation}\label{eqn:probdecom}
    \langle \varphi(F) \rangle =\varphi(0)+\int_0^\infty \langle I(F\ge \lambda)\rangle\varphi'(\lambda)\ud \lambda 
\end{equation} with $\varphi(x)=x^m$ and a change of variables $\tilde{\lambda}:=\lambda^s$, \begin{align*}
    \langle |F|^m \rangle \leftstackrel{\eqref{eqn:probdecom}}{\le} m(1+\exp(c^s))\int_0^\infty \lambda^{m-1}\exp(-\lambda^s)\ud\lambda & \le \dfrac{m(1+\exp(c^s))}{s}\int_0^\infty \tilde{\lambda}^{\frac{m}{s}-1}e^{-\tilde{\lambda}}\ud\tilde{\lambda} \\ &= \dfrac{m(1+\exp(c^s))}{s}\Gamma(\dfrac{m}{s}).
\end{align*}
Taking the $m$-th root and using the asymptotics of the Gamma function $\Gamma(\frac{m}{s})\lesssim \sqrt{m}(\frac{m}{se})^{\frac{m}{s}}$, we obtain as desired $\langle |F|^m \rangle^\frac{1}{m}\le Cm^\frac{1}{s}$ using $m^\frac{1}{m}\lesssim 1$. 
\smallskip

\noindent \emph{Step 2:} left $\lesssim$ right. Suppose $\langle \exp(|F|^s)\rangle \le e$, then for any constant $M$, \begin{equation*}
    \langle I(\dfrac{|F|}{M}+c>\lambda ) \rangle  =  \langle I(|F|>M(\lambda-c) ) \rangle \le \left\{\begin{aligned} & 1 & & 0<\lambda<c \\ & \exp\bigl(1-M^s(\lambda-c)^s\bigr) & & \lambda\ge c, \end{aligned} \right.
\end{equation*}
and therefore using again \eqref{eqn:probdecom} with $\varphi(x)=\exp(x^s)$, \begin{align*}
    \langle \exp\bigl((\dfrac{|F|}{M}+c)^s\bigr) \rangle & \leftstackrel{\eqref{eqn:probdecom}}{\le} 1+\int_0^c s\lambda^{s-1}\exp(\lambda^s)\ud \lambda + e\int_c^\infty s\lambda^{s-1}\exp(\lambda^s)\exp\bigl(-M^s(\lambda-c)^s\bigr)\ud\lambda \\ & \le \exp(c^s)+ e\int_c^\infty s\lambda^{s-1}\exp(\lambda^s)\exp\bigl(-M^s(\lambda-c)^s\bigr)\ud\lambda.
\end{align*}The above integrand decreases to zero pointwise as $M$ tends to infinity, and is obviously integrable for $M>1$, thus the integral also converges to zero. Hence for some sufficiently large $M>1$ only depending on $s$, $ \langle \exp\bigl((\frac{|F|}{M}+c)^s\bigr) \rangle \le \exp(c^s)+1$, which is equivalent to $\|F\|_s \le 1$. 
\smallskip

\noindent \emph{Step 3:} right $\lesssim$ middle. This direction of the proof uses the algebraic moments of random variables to construct the stretched exponential norm \eqref{eqn:defsnorm}. If \begin{equation}\label{eqn:2ndgtr1st}\dfrac{\langle |F|^m \rangle^\frac{1}{m}}{m^{1/s}}\le 1, \ \forall m\in \N,\end{equation}  then,  let $k_0$ be the largest integer such that $k_0s<1$. In the range $k\ge k_0+1$, which means $ks\ge 1$, for some large $C$ to be chosen later,\begin{equation*}
    \sum_{k=k_0+1}^\infty \dfrac{\langle |F|^{ks}\rangle}{C^{ks}k!} ~\leftstackrel{\eqref{eqn:2ndgtr1st}}{\le} \sum_{k=k_0+1}^\infty \dfrac{(ks)^k}{C^{ks}k!} \leftstackrel{\eqref{eqn:factorial}}{\le} \sum_{k=k_0+1}^\infty \dfrac{s^k}{C^{ks}c_0^k} =\dfrac{(\frac{s}{C^sc_0})^{k_0+1}}{1-\frac{s}{C^sc_0}} \le \dfrac{s}{C^sc_0-s}.
\end{equation*}
The lower moments are controlled by the higher moments: for $1\le k \le k_0$, using Jensen's inequality on $\langle \cdot \rangle$, we have \[
    \langle |F|^{ks} \rangle \le \langle |F|^{(k+k_0)s} \rangle^{\frac{k}{k+k_0}}\le 1+\langle |F|^{(k+k_0)s} \rangle. 
\] Therefore, using $C\ge 1$,  \begin{equation*}
    \sum_{k=1}^{k_0} \dfrac{\langle |F|^{ks}\rangle}{C^{ks}k!} \le C^{-s}\Bigl(k_0+(2k_0)!\sum_{k=k_0+1}^{2k_0}\dfrac{\langle |F|^{ks}\rangle}{k!}\Bigr) \leftstackrel{\eqref{eqn:2ndgtr1st}}{\le} C^{-s}\Bigl(k_0+(2k_0)!\sum_{k=k_0+1}^{2k_0}\dfrac{s^k}{c_0^k}\Bigr).
\end{equation*} Thus for fixed $s,c_0$ we may choose $C$ sufficiently large such that \begin{equation*}
    \langle \exp(\dfrac{|F|^s}{C^s})\rangle = 1+\sum_{k=1}^\infty  \dfrac{\langle |F|^{ks}\rangle}{C^{ks}k!} \le 1+C^{-s}\Bigl(k_0+(2k_0)!\sum_{k=k_0+1}^{2k_0}\dfrac{s^k}{c_0^k}\Bigr)+\dfrac{s}{C^sc_0-s} \le e.
\end{equation*} This establishes that the middle expression dominates the third one. Hence all three norms in \eqref{eqn:equistocnorm} are indeed equivalent.
\end{proof}

\noindent The second property of $\|\cdot\|_s$ is a ``stochastic H\"older inequality''.
\begin{lemma}\label{lem:holder} For $s,s_1,s_2>0$ with $\frac{1}{s}=\frac{1}{s_1}+\frac{1}{s_2}$, \begin{equation}\label{eqn:holder}\lVert F_1F_2\rVert_s\lesssim_{s,s_1,s_2} \lVert F_1 \rVert_{s_1}\lVert F_2 \rVert_{s_2}.\end{equation} Taking $F_2=1$ we obtain for any $s_1\le s_2$, \[\lVert F\rVert_{s_1}\lesssim \lVert F \rVert_{s_2}.\] \end{lemma}
\begin{proof}
The proof uses the third characterization of $\lVert\cdot\rVert_s$ in \eqref{eqn:equistocnorm}. Indeed, for any $M_1,M_2>0$, by Young's inequality and the convexity of $F\mapsto \log \langle \exp (F) \rangle$, \begin{align*}
	\log \Bigl\langle \exp(|\dfrac{F_1F_2}{M_1M_2}|^s) \Bigr\rangle &  \le \log \Bigl\langle \exp(\dfrac{s}{s_1}|\dfrac{F_1}{M_1}|^{s_1}+\dfrac{s}{s_2}|\dfrac{F_2}{M_2}|^{s_2})\Bigr\rangle \\ & \le \dfrac{s}{s_1}\log \Bigl\langle \exp(|\dfrac{F_1}{M_1}|^{s_1})\Bigr\rangle+\dfrac{s}{s_2}\log \Bigl\langle \exp(|\dfrac{F_2}{M_2}|^{s_2})\Bigr\rangle.
	\end{align*}
	Pick $M_1\ge \lVert F_1\rVert _{s_1}$ and $M_2 \ge \lVert F_2\rVert _{s_2}$, so that in view of the equivalence \eqref{eqn:equistocnorm}, the $\text{r.h.s.}\lesssim 1$. This means $\lVert F_1F_2\rVert _s \lesssim M_1M_2$. Since $M_1,M_2$ are arbitrary, this establishes \eqref{eqn:holder}.
\end{proof}
The inequality \eqref{eqn:holder} allows us to prove the following lemma, which is a generalization of \cite[Lemma~14]{GO2015most}. This is at the origin in the loss of stochastic integrability from $s=\infty$ for $ae$, to $s=2-$ for $\nabla \phi$ and to $s=1-$ for $\nabla \psi$. 
\begin{lemma}\label{lem:stochasticintegrability}
	Suppose for some $\delta\in(0,1]$ and $p<\infty$ and random variables $F,~F_0,~F_1$, we have
\begin{equation}\label{eqn:condholder}
	   F_0 \le \delta ~ \ \implies \ ~ |F|\le \delta^\frac{1}{p}F_1.
	\end{equation}
	Then for $s,~s_0,~s_1>0$ with $\frac{1}{s_0}=\frac{1}{s}+\frac{1}{s_1}$, and $\tilde{s}_1\ge s_1$ we have \begin{equation*}
 \lVert F\rVert_{s_0}\lesssim \delta^\frac{1}{p}\lVert F_1\rVert _{s_0}+ \Bigl(\dfrac{\lVert F_0\rVert _{\tilde{s}_1}}{\delta}\Bigr)^\frac{\tilde{s}_1}{s_1}\lVert F\rVert _s.
	\end{equation*}
\end{lemma}
\begin{proof}
By the triangle inequality, 
\begin{equation*}
\lVert F\rVert_{s_0} \le \lVert I(F_0 \le \delta)F\rVert _{s_0} + \lVert     I(F_0 \ge \delta)F\rVert _{s_0}.
\end{equation*}
By our two assumptions \eqref{eqn:condholder} and $ F_0 \le \delta$, $\lVert I(F_0\le \delta)F\rVert _{s_0} \le \delta^\frac{1}{p}\lVert F_1\rVert _{s_0}$. For the control of $\lVert I(F_0\ge \delta)F\rVert _{s_0}$, we have by Lemma \ref{lem:holder},
	\begin{equation*}
	\lVert I(F_0 \ge \delta)F\rVert _{s_0} \le \lVert I(F_0 \ge \delta)\rVert _{s_1} \lVert F \rVert_s.
	\end{equation*}
	The final step is a Chebyshev inequality: \begin{equation*}
	    \lVert I(F_0 \ge \delta)\rVert _{s_1} = \lVert I(F_0\ge \delta)\rVert_{\tilde{s}_1}^{\frac{\tilde{s}_1}{s_1}} \le \Bigl(\dfrac{\lVert F_0 \rVert _{\tilde{s}_1}}{\delta}\Bigr)^\frac{\tilde{s}_1}{s_1}.
	\end{equation*}
\end{proof}
\noindent The last property of $\lVert \cdot \rVert_s$ is that for a stationary random field $f$, $t\mapsto \lVert(\int \eta_{\sqrt{t}} f^2)^\frac{1}{2}\rVert_s $ is decreasing in $t$:
\begin{lemma}\label{lem:mixingscale}
	There exists a universal constant $C$ such that for any $t \le T$ and stationary random field $f$, \begin{equation}\label{eqn:mixingscale}
	\lVert(\int \eta_{\sqrt{T}} f^2)^\frac{1}{2}\rVert_s  \le C \lVert(\int \eta_{\sqrt{t}} f^2)^\frac{1}{2}\rVert_s .
	\end{equation}
\end{lemma}
\begin{proof}
    The proof uses \cite[(211)]{GO2015most} and the stationarity of $f$: 
\begin{align*}
\lVert (\int \eta_{\sqrt{T}} f^2)^\frac{1}{2}\rVert _{s}=\Big( \big\lVert \int\eta_{\sqrt{T}} f^2\big\rVert _{\frac{s}{2}}\Big)^\frac{1}{2} & \sim \Big(\big\lVert \int \eta_{\sqrt{T}}(y)\int \eta_{\sqrt{t}}(x-y)f^2(x)\ud x\ud y \big\rVert _\frac{s}{2}\Big)^\frac{1}{2}  \\ & \le \Big(\int \eta_{\sqrt{T}}(y)\big\lVert \int \eta_{\sqrt{t}}(x-y)f^2(x)\ud x\big\rVert _\frac{s}{2}\ud y \Big)^\frac{1}{2}= \lVert (\int \eta_{\sqrt{t}} f^2)^\frac{1}{2}\rVert _{s}. 
\end{align*}
\end{proof}

\section{Deterministic Estimates}
The next three lemmas are fundamental deterministic estimates. Lemma \ref{lem:GOlm3} is a localized energy estimate, Lemma \ref{lem:GOlm1} is a semigroup estimate that resolves the behavior in the time variable, and Lemma \ref{lem:GOlm12} is a crucial approximate locality estimate, which is the main tool for finding the approximate locality properties of several objects. These three lemmas are proved in \cite{GO2015most} and we do not present their proofs here.
\begin{lemma}
    (\cite[Lemma 3]{GO2015most}) \label{lem:GOlm3} Let $T>0$, and let $v, f$, and $g$ be related through the elliptic equation \begin{equation*}
        \dfrac{1}{T}v-\nabla\cdot a \nabla v=f+\nabla\cdot g,
    \end{equation*}
    then we have for all $R\ge \sqrt{T}$, \begin{equation}\label{eqn:GOlm3ell}
        \int \eta_R |(\nabla v, \dfrac{v}{\sqrt{T}})|^2 \lesssim \int \eta_R (|\sqrt{T}f|^2+|g|^2).
    \end{equation}
    Let $v, f, g$, and $v_0$ be related through the parabolic equation \begin{equation*}\begin{aligned}[rcl]
        \partial_t v-\nabla\cdot a \nabla v=f+\nabla\cdot g  &  & \mbox{ for } t>0, \\  v=v_0  &  & \mbox{ for }t=0,
    \end{aligned}\end{equation*}
    then we have for all $R\ge \sqrt{T}$, \begin{equation}\label{eqn:GOlm3para}
        \sup_{t\le T}\int \eta_R v^2(t) + \int_0^T\ud t \int \eta_R \bigl|(\nabla v(t), \dfrac{v(t)}{\sqrt{T}})\bigr|^2  \lesssim \int \eta_R v_0^2+ \int_0^T\ud t \int \eta_R (|\sqrt{T}f|^2+|g|^2).
    \end{equation}
\end{lemma}
\begin{lemma}(\cite[Lemma 1]{GO2015most})\label{lem:GOlm1}
    For all $R\ge \sqrt{T}>0$ we have  \begin{equation}\label{eqn:GOlm1}
\int \eta_R |(T\nabla S(T)g,\sqrt{T}S(T)g,\opS(T)g)|^2  + \int \eta_R \big|\int_0^T\ud t(\nabla S(t)g,\dfrac{S(t)g}{\sqrt{T}})\big|^2  \lesssim \int \eta_R |g|^2.
\end{equation}
\end{lemma}
\begin{lemma}\label{lem:GOlm12} (\cite[Lemma 12]{GO2015most})
For some $R\ge \sqrt{T}$ consider the following parabolic equation \begin{equation*} \left. \begin{aligned}
\partial_t w-\nabla\cdot a \nabla w=0 \ & \ & \mbox{ for }t>0 \\ w=0 \ & \ & \mbox{ for }t=0
\end{aligned} \right\} \ \mbox{ in} \  (0,T)\times B_{2R},
\end{equation*}
then we have for any exponent $p<\infty$, \begin{multline}\label{eqn:locality}
    \sup_{t\le T}\bigl(\fint_{B_R} w^2(t)\bigr)^\frac{1}{2} + R(\fint_{B_R} |\nabla w(T)|^2)^\frac{1}{2} +(\fint_{B_R} \int_0^T \ud t|\nabla w(t)|^2)^\frac{1}{2} \\ \lesssim_p (\dfrac{\sqrt{T}}{R})^p \fint_0^T \ud t(\dfrac{\sqrt{t}}{\sqrt{T}})^p \bigl(\fint_{B_{2R}\backslash B_R} w^2(t) \bigr)^\frac{1}{2} .
\end{multline}
\end{lemma}
\noindent Lemma \ref{lem:GOlm12} quantifies that $w\big|_{B_R\times (0,T)}$ depends very little on $w\big|_{B_{2R}^c \times (0,T)}$, provided $R\gg \sqrt{T}$. In addition, the estimate can afford a time average that degenerates as $t\ll T$. The algebraic smallness is sufficient for our purposes.

\section{CLT Estimates}
The following Lemma \ref{lem:GOlm13}, which restates \cite[Lemma 13]{GO2015most}, is the main probabilistic ingredient of this paper. It turns deterministic semigroup estimates on small scales and approximate locality properties into CLT-cancellations on large scales. The CLT-scaling typically arises from a combination of Lemmas \ref{lem:GOlm12} and \ref{lem:GOlm13}. We will give a proof here since the proof in \cite{GO2015most} is only valid when $s\in (1,2]$.
\begin{lemma}(\cite[Lemma 13]{GO2015most})\label{lem:GOlm13}
    Let $s \in(0,2]$. Let $F$ and $\bar{F}$ be stationary random fields such that $\langle F\rangle =0$ and suppose $F$ is approximately local on scale $\sqrt{T}$ relative to $\bar{F}$ in the sense of \eqref{eqn:apploc}. Then,
\begin{equation}\label{eqn:clt}
    \sup_{R\ge\sqrt{T}}(\dfrac{R}{\sqrt{T}})^\frac{d}{2} \lVert F_R\rVert _{s} \lesssim \sup_{r\le\sqrt{T}}(\dfrac{r}{\sqrt{T}})^\frac{d}{2} \lVert F_r\rVert _{s}+\Bigl(\sup_{r\le\sqrt{T}}(\dfrac{r}{\sqrt{T}})^\frac{d}{2} \lVert F_r\rVert _{s}  \Bigr)^{1-\frac{d}{2p}}\lVert \bar{F}\rVert_s^{\frac{d}{2p}}.
\end{equation}
\end{lemma}
\begin{proof}
By other parts of the proof of \cite[Lemma 13]{GO2015most} which are valid in the full range $s\in (0,2]$, we only need to prove the following: \\ Suppose $F$ is a mean-zero stationary random field that is exactly local on scale $r\ge 1$, in the sense that \begin{equation*}
    F(a)=F(\tilde{a}) \ \mbox{ provided } a=\tilde{a} \mbox{ on } B_r,
\end{equation*} then for any convolution kernel $\bar{G}$, \begin{equation}\label{eqn:lm510}
    \lVert \bar{G}*F\rVert _{s}\lesssim r^\frac{d}{2}(\int \bar{G}^2)^\frac{1}{2}\lVert F\rVert _{s}.
\end{equation}
The proof is inspired by \cite[Lemma 2.3]{su1997central}. We divide the proof into three steps:

\smallskip
\noindent \emph{Step 1}: A CLT in iid symmetric random variables.\\
Suppose $X_1,\cdots,X_n$ are iid symmetric random variables (in other words, the law of $X_1$ equals to the law of $-X_1$), then for any real numbers $w_1,\cdots,w_n$, we claim \begin{equation}\label{eqn:rvclt}
    \lVert \sum_{i=1}^n w_iX_i \rVert _s \lesssim (\sum_{i=1}^n w_i^2)^\frac{1}{2}\lVert X_1\rVert _s.
\end{equation}

Without loss of generality we assume $\lVert X_1 \rVert_s=1$. By Lemma \ref{lem:strexpnorm}, it is sufficient to show for any integer $m\ge 1$, \begin{equation*}
    \langle (\sum_{i=1}^n w_iX_i)^{2m} \rangle^\frac{1}{2m} \le C (2m)^\frac{1}{s}(\sum_{i=1}^n w_i^2)^\frac{1}{2}.
\end{equation*}
Let $c_0,c_1$ be constants as defined in \eqref{eqn:factorial}. Since the $X_i$ are symmetric, if we expand $\langle (\sum_{i=1}^n w_iX_i)^{2m}\rangle$, then all odd powers of $X_i$ will vanish, and therefore \begin{align*}
    \langle (\sum_{i=1}^n w_iX_i)^{2m}\rangle & =  \sum_{m_1+\cdots+m_n=m}\binom{2m}{2m_1,\cdots,2m_n}\prod_{i=1}^n \langle (w_iX_i)^{2m_i} \rangle \\ & \leftstackrel{\eqref{eqn:equistocnorm}}{\le} \sum_{m_1+\cdots+m_n=m}\dfrac{(2m)!}{\prod_{i=1}^n (2m_i)!}\prod_{j=1}^n w_j^{2m_j} \bigl(C(2m_j)^\frac{1}{s}\bigr)^{2m_j} \\ & \leftstackrel{\eqref{eqn:factorial}}{\le} C^{2m}\sum_{m_1+\cdots+m_n=m} \dfrac{c_1^{2m}(2m)^{2m}}{\prod_{i=1}^n c_0^{2m_i}(2m_i)^{2m_i}}\prod_{j=1}^n w_j^{2m_j}(2m_j)^{\frac{2m_j}{s}} \\ & \le (\dfrac{Cc_1}{c_0})^{2m}(2m)^{2m}\sum_{m_1+\cdots+m_n=m}\prod_{i=1}^n(2m_i)^{2m_i(\frac{1}{s}-1)}w_i^{2m_i}.
\end{align*}
Let $\varphi(t)=(2m)^{2m}\sum_{m_1+\cdots+m_n=m}\prod_{i=1}^n(2m_i)^{2m_i(t-1)}w_i^{2m_i}$. A direct calculation (only using $\log(2m_i)\le \log(2m)$) yields the differential inequality \begin{align*}
    \varphi'(t)&=(2m)^{2m}\sum_{m_1+\cdots+m_n=m}\Bigl(\prod_{i=1}^n(2m_i)^{2m_i(t-1)}w_i^{2m_i}\Bigr)\Bigl(\sum_{j=1}^n 2m_j\log(2m_j)\Bigr) \\ & \le 2m\log (2m) \varphi(t).
\end{align*} The value $t=\frac{1}{2}$ is special as it corresponds to $s=2$, since we have \begin{align*}
    \varphi(\dfrac{1}{2}) & = (2m)^{2m}\sum_{m_1+\cdots+m_n=m}\prod_{i=1}^n(2m_i)^{-m_i}w_i^{2m_i} \\ & = (2m)^m \sum_{m_1+\cdots+m_n=m} \dfrac{m^m}{\prod_{i=1}^n m_i^{m_i}}\prod_{j=1}^n w_j^{2m_j} \\ & \leftstackrel{\eqref{eqn:factorial}}{\le} (2m)^m \sum_{m_1+\cdots+m_n=m} \dfrac{c_1^m m!}{\prod_{i=1}^n c_0^{m_i}m_i!}\prod_{j=1}^n w_j^{2m_j} \\ &= (2m)^m(\dfrac{c_1}{c_0})^m \sum_{m_1+\cdots+m_n=m} \binom{m}{m_1,\cdots,m_n}\prod_{i=1}^nw_i^{2m_i} = (2m)^m(\dfrac{c_1}{c_0})^m(\sum_{i=1}^n w_i^2)^m.
\end{align*}
Therefore we obtain by integration for $t\ge \frac{1}{2}$ (or $s\le 2$), \begin{align*}
    \varphi(t) & \le (2m)^m(\dfrac{c_1}{c_0})^m(\sum_{i=1}^n w_i^2)^m \exp\bigl(2m\log(2m)(t-\dfrac{1}{2})\bigr) \\ & = (2m)^m(\dfrac{c_1}{c_0})^m(\sum_{i=1}^n w_i^2)^m (2m)^{2m(t-\frac{1}{2})} = (2m)^{2mt}(\dfrac{c_1}{c_0})^m(\sum_{i=1}^n w_i^2)^m.
\end{align*}
Finally we take the $2m$-th root:\begin{equation*}
    \langle (\sum_{i=1}^n w_iX_i)^{2m}\rangle^\frac{1}{2m} \le \dfrac{Cc_1}{c_0}(2m)^\frac{1}{s}(\dfrac{c_1}{c_0})^\frac{1}{2}(\sum_{i=1}^n w_i^2)^\frac{1}{2}.
\end{equation*}This finishes the proof.

\smallskip

\noindent \emph{Step 2}: From symmetric random variables to mean-zero random variables; we claim that \eqref{eqn:rvclt} holds just under the assumption that $\langle X_1\rangle =0$. 

For any $X$ such that $\langle X \rangle=0$, let $Y$ be independent of $X$ and equal to $-X$ in distribution. Let $\mu$ be the distribution of $X$ so that $f\mapsto \int \int f(x-y)\ud \mu(x)\ud\mu(y)$ denotes the distribution of $X+Y$, which is obviously symmetric. It remains to show \begin{equation}\label{eqn:equistep2}
    \lVert X+Y\rVert _s \sim \lVert X\rVert _s.
\end{equation} 
The ``$\lesssim$" direction is a direct consequence of the triangle inequality. For the other direction, notice that for any convex function $f$, \begin{equation}\label{eqn:probsym}
    \int \int f(x-y)\ud\mu(x)\ud \mu(y) \ge \int f(x-\int y \ud \mu(y)) \ud \mu(x) = \int f(x)\ud \mu(x).
\end{equation} We apply this to $f(x)=\exp\bigl((|x|+c)^s\bigr) -\exp(c^s)$, which is convex for $c$ defined in \eqref{eqn:defsnorm}. Hence, \eqref{eqn:probsym} turns into the estimate $ \lVert X+Y\rVert _s \gtrsim \lVert X\rVert_s$. This suffices to prove equivalence \eqref{eqn:equistep2}.

\smallskip

\noindent \emph{Step 3}: Conclusion. Without loss of generality we assume $\bar{G}$ is even, which can be removed in view of Step 4 of the proof of \cite[Lemma 13]{GO2015most}. We divide $\R^d$ into cubes of size $3r$, so that for a fixed $x\in[0,3r)^d$, $\{F(3z+x)\}_{z\in r\Z^d}$ are i.i.d. Therefore by stationarity, Step 2, passing to the limit $n\to\infty$, \begin{align*}
   \lVert \bar{G}*F \rVert_s=\lVert \bar{G}*F(0) \rVert_s= \lVert \int_{\R^d} \bar{G}(x)F(x) \ud x \rVert_s  & = \lVert  r^d \sum_{z\in r\Z^d}\fint_{[0,3r)^d} \bar{G}(x+3z)F(x+3z) \ud x \rVert_s  \\ &  \lesssim r^d\fint_{[0,3r)^d}\lVert   \sum_{z\in r\Z^d} \bar{G}(x+3z)F(x+3z)  \rVert_s \ud x \\ & \lesssim r^d\fint_{[0,3r)^d}(\sum_{z\in r\Z^d}\bar{G}^2(x+3z))^\frac{1}{2}\ud x \lVert F\rVert_s  \\ & \lesssim r^d(\fint_{[0,3r)^d}\sum_{z\in r\Z^d}\bar{G}^2(x+3z) \ud x )^\frac{1}{2}\lVert F\rVert_s  \\ & = r^\frac{d}{2}(\int_{\R^d}\bar{G}^2 )^\frac{1}{2}\lVert F\rVert_{s}. \qedhere
\end{align*}The rest of the proof for this lemma follows from Steps 4-6 of \cite[Lemma 13]{GO2015most}. \end{proof}

\section{Large scale regularity for stochastic homogenization}
 Below we present the crucial elliptic and parabolic large scale regularity estimates on scales larger than $\rr$ defined in \eqref{eqn:regminrad}. This slightly differs from \cite[Theorem 8.7]{armstrong2017quantitative}, \cite[Corollary 6]{GO2015most} since we use exponential averaging $\int \eta_R$ as opposed to the usual averaging $\fint_{B_R}$. In Theorem \ref{thm:mainthm}, we use the growth properties \eqref{intrpsiorig} of second-order corrector and flux $(\psi,\Psi)$  and the stochastic moment of $\rr$ to derive stochastic bounds of $\rss$. 

\begin{lemma}[$C^{0,1}$-estimate] \label{lem:parabolicmvp}
	The quantity $\rr$ satisfies the following: \begin{enumerate}[label=(\arabic*)]\item For every $R\ge \rr$ and $u$ satisfying \begin{equation*}
	\partial_t u-\nabla\cdot a \nabla u=0, \ \mbox{ in }(-R^2,0)\times \mathbb{R}^d
	\end{equation*}
	we have, for every $r\in [\rr,R]$, \begin{align}
	(\fint_{-r^2}^0 \ud t\int \eta_r |\nabla u(t)|^2)^\frac{1}{2} & \lesssim (\fint_{-R^2}^0\ud t\int \eta_R |\nabla u(t)|^2)^\frac{1}{2}, \label{eqn:mvpgrad} \\ 		\bigl(\fint_{-r^2}^0 \ud t\int \eta_r  u^2(t)\bigr)^\frac{1}{2} & \lesssim \bigl(\fint_{-R^2}^0\ud t\int \eta_R  u^2(t)\bigr)^\frac{1}{2}.\label{eqn:mvpnograd}
	\end{align}
	\item \cite[Corollary 6]{GO2015most},\cite[Theorem 1.2]{armstrong2016quantitative} 
 For every $R\ge \rr$, $u$ satisfying 
	\begin{equation*}
	-\nabla\cdot a \nabla u=0\text{ in } B_{R},
	\end{equation*}
	we have, for every $r\in [\rr,R]$, \begin{equation}\label{eqn:ellipticmvp}
	(\fint_{B_r} |\nabla u|^2)^\frac{1}{2} \lesssim (\fint_{B_R} |\nabla u|^2)^\frac{1}{2}.
	\end{equation} \end{enumerate}
\end{lemma} 
\begin{proof}
\noindent \emph{Step 1:} we prove \eqref{eqn:mvpgrad} from the usual version of $C^{0,1}$ estimate \cite[Remark 8.8]{armstrong2017quantitative}: for any $R\ge r\ge \rr(0)$, \begin{equation}\label{eqn:usualmvpgrad} (\fint_{B_r^{(P)}} |\nabla u|^2)^\frac{1}{2} \le C(\fint_{B_R^{(P)}} |\nabla u|^2)^\frac{1}{2} \ \mbox{ where }B_r^{(P)}:=(-r^2,0)\times Q_r.\end{equation}
Without loss of generality we assume $r=1$ and $R\ge 2$. We first claim that exponential average can be recovered from the ordinary one as follows:  \begin{equation}\label{eqn:recavg}
    \int \eta_R|\nabla u|^2 \sim \int_R^\infty\ud \rho \dfrac{\rho^d}{R^{d+1}} \exp(-\dfrac{\rho}{R})\fint_{B_\rho}|\nabla u|^2.
\end{equation}
Indeed, let $S(\rho)=\{x\in \R^d: |x|=\rho\}$ be the $d$-dimensional sphere of radius $r$. Then by Fubini's theorem \begin{align*}
    \int_R^\infty \ud \rho\dfrac{\rho^d}{R^{d+1}} \exp(-\dfrac{\rho}{R})\fint_{B_\rho}|\nabla u|^2 & \sim \int_R^\infty \dfrac{\ud \rho}{R^{d+1}} \exp(-\dfrac{\rho}{R})\int_0^\rho \ud \rho' \int_{S(\rho')}|\nabla u|^2 \\ & \sim \int_0^R \dfrac{\ud \rho'}{R^{d+1}}  \int_{S(\rho')}|\nabla u|^2 \int_R^\infty\ud \rho \exp(-\dfrac{\rho}{R}) \\ & \qquad +\int_R^\infty \dfrac{\ud \rho'}{R^{d+1}}  \int_{S(\rho')}|\nabla u|^2 \int_{\rho'}^\infty\ud \rho \exp(-\dfrac{\rho}{R}) \\ & \sim \int_0^R \dfrac{\ud \rho'}{R^d}  \int_{S(\rho')}|\nabla u|^2+ \int_R^\infty \dfrac{\ud \rho'}{R^d}   \exp(-\dfrac{\rho'}{R})\int_{S(\rho')}|\nabla u|^2 \\ & \sim \int \eta_R |\nabla u|^2.
\end{align*}

 We now prove \eqref{eqn:mvpgrad}. The strategy is to divide the $r$ integral into $(1,R)$ where we may appeal to \eqref{eqn:usualmvpgrad}, and $(R,\infty)$ where we use the faster decay of $\eta_1$ compared to $\eta_R$: \begin{align*}
    \fint_{-1}^0 & \ud t  \int \eta_1|\nabla u|^2 \\ & \leftstackrel{\eqref{eqn:recavg}}{\sim}  \int_{-1}^0 \ud t\int_1^R\ud \rho \rho^d \exp(-\rho)\fint_{B_\rho}|\nabla u|^2 + \int_{-1}^0 \ud t\int_R^\infty \ud \rho \rho^d \exp(-\rho)\fint_{B_\rho}|\nabla u|^2 \\ & \lesssim  \int_1^R\ud \rho \rho^{d+2}\exp(-\rho)\fint_{-\rho^2}^0 \ud t \fint_{B_\rho}|\nabla u|^2 + R^2 \fint_{-R^2}^0 \ud t\int_R^\infty \ud \rho \rho^d \exp(-\rho)\fint_{B_\rho}|\nabla u|^2\\ & \leftstackrel{\eqref{eqn:usualmvpgrad}}{\lesssim} \int_1^R\ud \rho \rho^{d+2}\exp(-\rho)\fint_{-R^2}^0 \ud t \fint_{B_R}|\nabla u|^2 + \fint_{-R^2}^0 \ud t\int_R^\infty \ud \rho \dfrac{\rho^d}{R^{d+1}} \exp(-\dfrac{\rho}{R})\fint_{B_\rho}|\nabla u|^2 \\ &\leftstackrel{\eqref{eqn:recavg}}{\lesssim} \fint_{-R^2}^0 \ud t \fint_{B_R}|\nabla u|^2 + \fint_{-R^2}^0 \ud t\int \eta_R |\nabla u|^2 \\ & \lesssim \fint_{-R^2}^0 \ud t \int \eta_R |\nabla u|^2.
\end{align*} 
Here the second term of the third inequality uses that $R^{d+3} \exp(-\rho)\lesssim \exp(-\frac{\rho}{R})$ when $\rho\ge R\ge 2$. This establishes \eqref{eqn:mvpgrad}.

\smallskip

\noindent \emph{Step 2:} proof of \eqref{eqn:mvpnograd}. The proof starts from the usual version of $C^{0,1}$ estimate without gradient \cite[Theorem 8.7]{armstrong2017quantitative}: for any $R\ge r\ge \rr(0)$, \begin{equation}\label{eqn:usualmvpnograd}
	(\fint_{B_r^{(P)}} |u-\fint_{B_r^{(P)}} u|^2)^\frac{1}{2} \lesssim \dfrac{r}{R} (\fint_{B_R^{(P)}} |u-\fint_{B_R^{(P)}} u|^2)^\frac{1}{2}.
\end{equation}
It suffices to prove \begin{equation}\label{eqn:midpsmvpnograd}
    (\fint_{B_r^{(P)}} |u|^2)^\frac{1}{2} \lesssim (\fint_{B_R^{(P)}} |u|^2)^\frac{1}{2}.
\end{equation} since we can then change spatial averaging from $\fint_{B_R}$ to $\int \eta_R$ using the same post-processing argument as above.

Let $u_r:=\fint_{B_r^{(P)}} u$. W.l.o.g. we may assume $R\ge 2r$. Now we can replace $r$ with $2r$ in \eqref{eqn:usualmvpnograd} and get  
\begin{equation*}
    (\fint_{B_{2r}^{(P)}} |u-u_{2r}|^2)^\frac{1}{2} \lesssim \dfrac{r}{R} (\fint_{B_R^{(P)}} |u-u_R|^2)^\frac{1}{2}.
\end{equation*}
Thus we can estimate 
\begin{align*}
    |u_r-u_{2r}|&  = |\fint_{B_r^{(P)}} (u-u_{2r})| \lesssim (\fint_{B_r^{(P)}} |u-u_{2r}|^2)^\frac{1}{2} \\ &\lesssim (\fint_{B_{2r}^{(P)}} |u-u_{2r}|^2)^\frac{1}{2}\lesssim \dfrac{r}{R} (\fint_{B_R^{(P)}} |u-u_R|^2)^\frac{1}{2}. \stepcounter{equation}\tag{\theequation}\label{eqn:mvptelescoping}
\end{align*}
We now obtain \eqref{eqn:midpsmvpnograd} by telescoping
\begin{align*}
    (\fint_{B_r^{(P)}} |u|^2)^\frac{1}{2} & \lesssim (\fint_{B_r^{(P)}} |u-u_r|^2)^\frac{1}{2}+|u_r| \\ & \leftstackrel{\eqref{eqn:usualmvpnograd}}{\lesssim} \dfrac{r}{R} (\fint_{B_R^{(P)}} |u-u_R|^2)^\frac{1}{2} + \sum_{\substack {\tilde{r}\in [r,R) \\ \text{dyadic}}} |u_{\tilde{r}}-u_{2\tilde{r}}| + |u_R| \\ & \leftstackrel{\eqref{eqn:mvptelescoping}}{\lesssim} \dfrac{r}{R} (\fint_{B_R^{(P)}} |u|^2)^\frac{1}{2} + \sum_{\substack {\tilde{r}\in [r,R) \\ \text{dyadic}}} \dfrac{\tilde{r}}{R} (\fint_{B_R^{(P)}} |u-u_R|^2)^\frac{1}{2} + |u_R| \\ & \lesssim   (\fint_{B_R^{(P)}} |u|^2)^\frac{1}{2}. 
\end{align*}
\end{proof}
\noindent The following corollary is a convenient post-processing of Lemma \ref{lem:parabolicmvp} from time intervals to single time slices.

\begin{corollary}\label{cor:upgrademvp} 
For any $\sqrt{t}\ge r \ge \rr$, we have\begin{align}
        &  (\int \eta_r |\nabla S(t)g|^2)^\frac{1}{2} \lesssim (\int \eta_{\sqrt{t}} |\nabla S(\dfrac{t}{2})g|^2)^\frac{1}{2},\label{eqn:detppsmvpgrad} \\
 &       (\int \eta_r | S(t)g|^2)^\frac{1}{2} \lesssim (\int \eta_{\sqrt{t}} | S(\dfrac{t}{2})g|^2)^\frac{1}{2}. \label{eqn:detppsmvpnograd} 
    \end{align} 
\end{corollary}
\begin{proof}
We will only prove \eqref{eqn:detppsmvpgrad} as the proof of \eqref{eqn:detppsmvpnograd} is identical. The estimate uses \eqref{eqn:mvpgrad} and the fact that the function $\tau\mapsto \int \eta_{\sqrt{t}} |\nabla S(\tau)g|^2$ is approximately non-increasing for any $g$, that is, for any $t_1 \le t_2 \le r^2$, $\int \eta_r |\nabla S(t_1)g|^2 \lesssim \int \eta_r |\nabla S(t_2)g|^2$ with the constant in $\lesssim$ independent of $r,t_1,t_2,g$:
\begin{align*}
    (\int \eta_r |\nabla S(t)g |^2)^\frac{1}{2} & \lesssim  (\int\eta_r \fint_{t-r^2}^t\ud \tau| \nabla S(\tau)g |^2)^\frac{1}{2}\\ & \leftstackrel{\eqref{eqn:mvpgrad}}{\lesssim } (\int \eta_{\sqrt{t}}\fint_{\frac{t}{2}}^t\ud \tau| \nabla S(\tau)g |^2)^\frac{1}{2} \\ & \lesssim(\int \eta_{\sqrt{t}}| \nabla S(\dfrac{t}{2})g |^2)^\frac{1}{2}. \end{align*}
\end{proof}
\begin{remark}\label{rmk:postproc}
     We further post-process Corollary \ref{cor:upgrademvp} into the following: for any $t,r\ge 0$,  \begin{align}
        & \|I(r\ge \rr) (\int \eta_r |\nabla S(t)g|^2)^\frac{1}{2}\|_s \lesssim \|I(r \vee \sqrt{t}\ge \rr)(\int \eta_{r\vee\sqrt{t}} |\nabla S(\dfrac{t}{2})g|^2)^\frac{1}{2}\|_s, \label{eqn:ppsmvpgrad} \\
 &       \|I(r\ge \rr)(\int \eta_r  | S(t)g|^2)^\frac{1}{2}\|_s \lesssim \|I(r \vee \sqrt{t}\ge \rr)(\int \eta_{r\vee\sqrt{t}} | S(\dfrac{t}{2})g|^2)^\frac{1}{2}\|_s. \label{eqn:ppsmvpnograd} 
    \end{align} The argument uses the fact that the function $ [0,r^2) \ni t \mapsto \int \eta_r |\nabla S(t)g|^2$ and $ [0,r^2) \ni t \mapsto \int \eta_r | S(t)g|^2$ are both approximately non-increasing for any $g$ (see the proof of \cite[Lemma 3]{GO2015most}).
\end{remark}

\bibliographystyle{alpha}
\bibliography{HomNumAlg}

\newcommand{\etalchar}[1]{$^{#1}$}
\begin{thebibliography}{AEEVE12}

\bibitem[ACB{\etalchar{+}}12]{anantharaman2012introduction}
Arnaud Anantharaman, Ronan Costaouec, Claude~Le Bris, Fr{\'e}d{\'e}ric Legoll,
  and Florian Thomines.
\newblock Introduction to numerical stochastic homogenization and the related
  computational challenges: some recent developments.
\newblock In {\em Multiscale modeling and analysis for materials simulation},
  pages 197--272. World Scientific, 2012.

\bibitem[AEEVE12]{abdulle2012heterogeneous}
Assyr Abdulle, Weinan E, Bj{\"o}rn Engquist, and Eric Vanden-Eijnden.
\newblock The heterogeneous multiscale method.
\newblock {\em Acta Numerica}, 21:1--87, 2012.

\bibitem[AHKM21]{armstrong2018iterative}
Scott Armstrong, Antti Hannukainen, Tuomo Kuusi, and Jean-Christophe Mourrat.
\newblock An iterative method for elliptic problems with rapidly oscillating
  coefficients.
\newblock {\em ESAIM: Mathematical Modelling and Numerical Analysis},
  55(1):37--55, 2021.

\bibitem[AKM16]{armstrong2016mesoscopic}
Scott Armstrong, Tuomo Kuusi, and Jean-Christophe Mourrat.
\newblock Mesoscopic higher regularity and subadditivity in elliptic
  homogenization.
\newblock {\em Communications in Mathematical Physics}, 347(2):315--361, 2016.

\bibitem[AKM17]{armstrong2017additive}
Scott Armstrong, Tuomo Kuusi, and Jean-Christophe Mourrat.
\newblock The additive structure of elliptic homogenization.
\newblock {\em Inventiones mathematicae}, 208(3):999--1154, 2017.

\bibitem[AKM19]{armstrong2017quantitative}
Scott Armstrong, Tuomo Kuusi, and Jean-Christophe Mourrat.
\newblock {\em Quantitative stochastic homogenization and large-scale
  regularity}, volume 352.
\newblock Springer, 2019.

\bibitem[AL87]{avellaneda1987compactness}
Marco Avellaneda and Fang-Hua Lin.
\newblock Compactness methods in the theory of homogenization.
\newblock {\em Communications on Pure and Applied Mathematics}, 40(6):803--847,
  1987.

\bibitem[AS16]{armstrong2016quantitative}
Scott~N Armstrong and Charles~K Smart.
\newblock Quantitative stochastic homogenization of convex integral
  functionals.
\newblock In {\em Annales scientifiques de l'Ecole normale sup{\'e}rieure},
  volume~49, pages 423--481. Societe Mathematique de France, 2016.

\bibitem[BCO94]{babuvska1994special}
Ivo Babu{\v{s}}ka, Gabriel Caloz, and John~E Osborn.
\newblock Special finite element methods for a class of second order elliptic
  problems with rough coefficients.
\newblock {\em SIAM Journal on Numerical Analysis}, 31(4):945--981, 1994.

\bibitem[BFFO17]{bella2017stochastic}
Peter Bella, Benjamin Fehrman, Julian Fischer, and Felix Otto.
\newblock Stochastic homogenization of linear elliptic equations: Higher-order
  error estimates in weak norms via second-order correctors.
\newblock {\em SIAM Journal on Mathematical Analysis}, 49(6):4658--4703, 2017.

\bibitem[BGO15]{bella2015quantitative}
Peter Bella, Arianna Giunti, and Felix Otto.
\newblock Quantitative stochastic homogenization: local control of
  homogenization error through corrector.
\newblock {\em arXiv preprint arXiv:1504.02487}, 2015.

\bibitem[BGO20]{bella2017effective}
Peter Bella, Arianna Giunti, and Felix Otto.
\newblock Effective multipoles in random media.
\newblock {\em Communications in Partial Differential Equations},
  45(6):561--640, 2020.

\bibitem[BLBL16]{blanc2016some}
Xavier Blanc, Claude Le~Bris, and Fr{\'e}d{\'e}ric Legoll.
\newblock Some variance reduction methods for numerical stochastic
  homogenization.
\newblock {\em Philosophical Transactions of the Royal Society A: Mathematical,
  Physical and Engineering Sciences}, 374(2066):20150168, 2016.

\bibitem[BO83]{babuvska1983generalized}
Ivo Babu{\v{s}}ka and John~E Osborn.
\newblock Generalized finite element methods: their performance and their
  relation to mixed methods.
\newblock {\em SIAM Journal on Numerical Analysis}, 20(3):510--536, 1983.

\bibitem[BO16]{bella2016corrector}
Peter Bella and Felix Otto.
\newblock Corrector estimates for elliptic systems with random periodic
  coefficients.
\newblock {\em Multiscale Modeling \& Simulation}, 14(4):1434--1462, 2016.

\bibitem[Bou18]{bourgain2018homogenization}
Jean Bourgain.
\newblock On a homogenization problem.
\newblock {\em Journal of Statistical Physics}, 172(2):314--320, 2018.

\bibitem[BP04]{bourgeat2004approximations}
Alain Bourgeat and Andrey Piatnitski.
\newblock Approximations of effective coefficients in stochastic
  homogenization.
\newblock In {\em Annales de l'IHP Probabilit{\'e}s et statistiques},
  volume~40, pages 153--165, 2004.

\bibitem[CELS15]{cances2015embedded}
Eric Cances, Virginie Ehrlacher, Fr{\'e}d{\'e}ric Legoll, and Benjamin Stamm.
\newblock An embedded corrector problem to approximate the homogenized
  coefficients of an elliptic equation.
\newblock {\em Comptes Rendus Mathematique}, 353(9):801--806, 2015.

\bibitem[DGL19]{duerinckx2019remark}
Mitia Duerinckx, Antoine Gloria, and Marius Lemm.
\newblock A remark on a surprising result by bourgain in homogenization.
\newblock {\em Communications in Partial Differential Equations},
  44(12):1345--1357, 2019.

\bibitem[DO20]{duerinckx2020higher}
Mitia Duerinckx and Felix Otto.
\newblock Higher-order pathwise theory of fluctuations in stochastic
  homogenization.
\newblock {\em Stochastics and Partial Differential Equations: Analysis and
  Computations}, 8(3):625--692, 2020.

\bibitem[Due21]{duerinckx2021non}
Mitia Duerinckx.
\newblock Non-perturbative approach to the {B}ourgain-{S}pencer conjecture in
  stochastic homogenization.
\newblock {\em arXiv preprint arXiv:2102.06319}, 2021.

\bibitem[EEH03]{weinan2003heterogeneous}
Weinan E, Bjorn Engquist, and Zhongyi Huang.
\newblock Heterogeneous multiscale method: a general methodology for multiscale
  modeling.
\newblock {\em Physical Review B}, 67(9):092101, 2003.

\bibitem[EEL{\etalchar{+}}07]{weinan2007heterogeneous}
Weinan E, Bjorn Engquist, Xiantao Li, Weiqing Ren, and Eric Vanden-Eijnden.
\newblock The heterogeneous multiscale method: A review.
\newblock In {\em Commun. Comput. Phys}. Citeseer, 2007.

\bibitem[EGMN15]{egloffe2015random}
A-C Egloffe, Antoine Gloria, J-C Mourrat, and Thahn~Nhan Nguyen.
\newblock Random walk in random environment, corrector equation and homogenized
  coefficients: from theory to numerics, back and forth.
\newblock {\em IMA journal of numerical analysis}, 35(2):499--545, 2015.

\bibitem[EH09]{efendiev2009multiscale}
Yalchin Efendiev and Thomas~Y Hou.
\newblock {\em Multiscale finite element methods: theory and applications},
  volume~4.
\newblock Springer Science \& Business Media, 2009.

\bibitem[EMZ{\etalchar{+}}05]{ming2005analysis}
Weinan E, Pingbing Ming, Pingwen Zhang, et~al.
\newblock Analysis of the heterogeneous multiscale method for elliptic
  homogenization problems.
\newblock {\em Journal of the American Mathematical Society}, 18(1):121--156,
  2005.

\bibitem[FGP21]{fischer2019priori}
Julian Fischer, Dietmar Gallistl, and Daniel Peterseim.
\newblock A priori error analysis of a numerical stochastic homogenization
  method.
\newblock {\em SIAM Journal on Numerical Analysis}, 59(2):660--674, 2021.

\bibitem[Fis19]{fischer2019choice}
Julian Fischer.
\newblock The choice of representative volumes in the approximation of
  effective properties of random materials.
\newblock {\em Archive for Rational Mechanics and Analysis}, 234(2):635--726,
  2019.

\bibitem[FO16]{fischer2016higher}
Julian Fischer and Felix Otto.
\newblock A higher-order large-scale regularity theory for random elliptic
  operators.
\newblock {\em Communications in partial differential equations},
  41(7):1108--1148, 2016.

\bibitem[FO17]{fischer2017sublinear}
Julian Fischer and Felix Otto.
\newblock Sublinear growth of the corrector in stochastic homogenization:
  optimal stochastic estimates for slowly decaying correlations.
\newblock {\em Stochastics and Partial Differential Equations: Analysis and
  Computations}, 5(2):220--255, 2017.

\bibitem[Glo12]{gloria2012numerical}
Antoine Gloria.
\newblock Numerical approximation of effective coefficients in stochastic
  homogenization of discrete elliptic equations.
\newblock {\em ESAIM: Mathematical Modelling and Numerical Analysis},
  46(1):1--38, 2012.

\bibitem[GM12]{gloria2012spectral}
Antoine Gloria and Jean-Christophe Mourrat.
\newblock Spectral measure and approximation of homogenized coefficients.
\newblock {\em Probability theory and related fields}, 154(1-2):287--326, 2012.

\bibitem[GM16]{gu2016scaling}
Yu~Gu and Jean-Christophe Mourrat.
\newblock Scaling limit of fluctuations in stochastic homogenization.
\newblock {\em Multiscale Modeling \& Simulation}, 14(1):452--481, 2016.

\bibitem[GN16]{gloria2016quantitative}
Antoine Gloria and James Nolen.
\newblock A quantitative central limit theorem for the effective conductance on
  the discrete torus.
\newblock {\em Communications on pure and applied mathematics},
  69(12):2304--2348, 2016.

\bibitem[GNO15]{gloria2015quantification}
Antoine Gloria, Stefan Neukamm, and Felix Otto.
\newblock Quantification of ergodicity in stochastic homogenization: optimal
  bounds via spectral gap on glauber dynamics.
\newblock {\em Inventiones mathematicae}, 199(2):455--515, 2015.

\bibitem[GNO20]{gloria2014regularity}
Antoine Gloria, Stefan Neukamm, and Felix Otto.
\newblock A regularity theory for random elliptic operators.
\newblock {\em Milan Journal of Mathematics}, 88(1):99--170, 2020.

\bibitem[GNO21]{gloria2019quantitative}
Antoine Gloria, Stefan Neukamm, and Felix Otto.
\newblock Quantitative estimates in stochastic homogenization for correlated
  coefficient fields.
\newblock {\em Analysis \& PDE}, 14(8):2497--2537, 2021.

\bibitem[GO11]{gloria2011optimal}
Antoine Gloria and Felix Otto.
\newblock An optimal variance estimate in stochastic homogenization of discrete
  elliptic equations.
\newblock {\em The annals of probability}, 39(3):779--856, 2011.

\bibitem[GO12]{gloria2012optimal}
Antoine Gloria and Felix Otto.
\newblock An optimal error estimate in stochastic homogenization of discrete
  elliptic equations.
\newblock {\em The annals of applied probability}, pages 1--28, 2012.

\bibitem[GO15]{GO2015most}
Antoine Gloria and Felix Otto.
\newblock The corrector in stochastic homogenization: optimal rates, stochastic
  integrability, and fluctuations.
\newblock {\em arXiv preprint arXiv:1510.08290}, 2015.

\bibitem[GO17]{gloria2017quantitative}
Antoine Gloria and Felix Otto.
\newblock Quantitative results on the corrector equation in stochastic
  homogenization.
\newblock {\em Journal of the European Mathematical Society},
  19(11):3489--3548, 2017.

\bibitem[Gu17]{gu2017high}
Yu~Gu.
\newblock High order correctors and two-scale expansions in stochastic
  homogenization.
\newblock {\em Probability Theory and Related Fields}, 169(3-4):1221--1259,
  2017.

\bibitem[HMS21]{hannukainen2019computing}
Antti Hannukainen, Jean-Christophe Mourrat, and Harmen~T Stoppels.
\newblock Computing homogenized coefficients via multiscale representation and
  hierarchical hybrid grids.
\newblock {\em ESAIM: Mathematical Modelling and Numerical Analysis},
  55:S149--S185, 2021.

\bibitem[HW97]{hou1997multiscale}
Thomas~Y Hou and Xiao-Hui Wu.
\newblock A multiscale finite element method for elliptic problems in composite
  materials and porous media.
\newblock {\em Journal of computational physics}, 134(1):169--189, 1997.

\bibitem[JKO94]{jikov2012homogenization}
Vasili~Vasilievitch Jikov, Sergei~M Kozlov, and Olga~Arsenievna Oleinik.
\newblock {\em Homogenization of differential operators and integral
  functionals}.
\newblock Springer Science \& Business Media, 1994.

\bibitem[KKO20]{khoromskaia2019numerical}
Venera Khoromskaia, Boris~N Khoromskij, and Felix Otto.
\newblock Numerical study in stochastic homogenization for elliptic partial
  differential equations: Convergence rate in the size of representative volume
  elements.
\newblock {\em Numerical Linear Algebra with Applications}, 27(3):e2296, 2020.

\bibitem[KL19]{kim2019averaged}
Jongchon Kim and Marius Lemm.
\newblock On the averaged {G}reen's function of an elliptic equation with
  random coefficients.
\newblock {\em Archive for Rational Mechanics and Analysis}, 234(3):1121--1166,
  2019.

\bibitem[Koz79]{kozlov1979averaging}
Sergei~Mikhailovich Kozlov.
\newblock Averaging of random operators.
\newblock {\em Matematicheskii Sbornik}, 151(2):188--202, 1979.

\bibitem[LBLM16]{le2016special}
Claude Le~Bris, Fr{\'e}d{\'e}ric Legoll, and William Minvielle.
\newblock Special quasirandom structures: a selection approach for stochastic
  homogenization.
\newblock {\em Monte Carlo Methods and Applications}, 22(1):25--54, 2016.

\bibitem[LO21]{lu2018optimal}
Jianfeng Lu and Felix Otto.
\newblock Optimal artificial boundary condition for random elliptic media.
\newblock {\em Foundations of Computational Mathematics}, pages 1--60, 2021.

\bibitem[LT13]{ledoux2013probability}
Michel Ledoux and Michel Talagrand.
\newblock {\em Probability in Banach Spaces: isoperimetry and processes}.
\newblock Springer Science \& Business Media, 2013.

\bibitem[MN17]{mourrat2017scaling}
Jean-Christophe Mourrat and James Nolen.
\newblock Scaling limit of the corrector in stochastic homogenization.
\newblock {\em The Annals of Applied Probability}, 27(2):944--959, 2017.

\bibitem[MO15]{marahrens2015annealed}
Daniel Marahrens and Felix Otto.
\newblock Annealed estimates on the green function.
\newblock {\em Probability theory and related fields}, 163(3-4):527--573, 2015.

\bibitem[MO16]{mourrat2016correlation}
Jean-Christophe Mourrat and Felix Otto.
\newblock Correlation structure of the corrector in stochastic homogenization.
\newblock {\em The annals of probability}, 44(5):3207--3233, 2016.

\bibitem[Mou19]{mourrat2019efficient}
J-C Mourrat.
\newblock Efficient methods for the estimation of homogenized coefficients.
\newblock {\em Foundations of Computational Mathematics}, 19(2):435--483, 2019.

\bibitem[MP14]{maalqvist2014localization}
Axel M{\aa}lqvist and Daniel Peterseim.
\newblock Localization of elliptic multiscale problems.
\newblock {\em Mathematics of Computation}, 83(290):2583--2603, 2014.

\bibitem[Nol14]{nolen2014normal}
James Nolen.
\newblock Normal approximation for a random elliptic equation.
\newblock {\em Probability Theory and Related Fields}, 159(3-4):661--700, 2014.

\bibitem[NS97]{naddaf1997homogenization}
Ali Naddaf and Thomas Spencer.
\newblock On homogenization and scaling limit of some gradient perturbations of
  a massless free field.
\newblock {\em Communications in mathematical physics}, 183(1):55--84, 1997.

\bibitem[NS98]{naddaf1998estimates}
Ali Naddaf and Thomas Spencer.
\newblock Estimates on the variance of some homogenization problems.
\newblock {\em preprint}, 1998.

\bibitem[OV00]{oden2000estimation}
J~Tinsley Oden and Kumar~S Vemaganti.
\newblock Estimation of local modeling error and goal-oriented adaptive
  modeling of heterogeneous materials: I. error estimates and adaptive
  algorithms.
\newblock {\em Journal of Computational Physics}, 164(1):22--47, 2000.

\bibitem[Owh15]{owhadi2015bayesian}
Houman Owhadi.
\newblock Bayesian numerical homogenization.
\newblock {\em Multiscale Modeling \& Simulation}, 13(3):812--828, 2015.

\bibitem[Owh17]{owhadi2017multigrid}
Houman Owhadi.
\newblock Multigrid with rough coefficients and multiresolution operator
  decomposition from hierarchical information games.
\newblock {\em SIAM Review}, 59(1):99--149, 2017.

\bibitem[PV79]{papanicolaou1979boundary}
George~C Papanicolaou and S.~R.~S. Varadhan.
\newblock Boundary value problems with rapidly oscillating random coefficients.
\newblock In {\em Colloquia Math. Soc., Janos Bolyai}, volume~27, pages
  853--873, 1979.

\bibitem[She18]{shen2018periodic}
Zhongwei Shen.
\newblock {\em Periodic Homogenization of Elliptic Systems}, volume 269.
\newblock Springer, 2018.

\bibitem[Su97]{su1997central}
Zhonggen Su.
\newblock Central limit theorems for random processes with sample paths in
  exponential orlicz spaces.
\newblock {\em Stochastic Processes and their applications}, 66(1):1--20, 1997.

\bibitem[Yur86]{yurinskii1986averaging}
VV~Yurinskii.
\newblock Averaging of symmetric diffusion in random medium.
\newblock {\em Siberian Mathematical Journal}, 27(4):603--613, 1986.

\end{thebibliography}

\end{document}